\documentclass[10pt, letterpaper,reqno]{amsart}
\usepackage[left=1in,right=1in,bottom=0.5in,top=0.8in]{geometry}
\usepackage{amsfonts}
\usepackage{amsmath, amssymb}
\usepackage{graphicx}
\usepackage[font=small,labelfont=bf]{caption}
\usepackage{epstopdf}
\usepackage{xcolor}
\usepackage{amsthm}
\usepackage{float}
\usepackage{enumitem}
\usepackage{pgfplots}
\usepackage{listings}
\usepackage{longtable}
\usepackage{mathrsfs}
\usepackage{dsfont}
\usepackage{mathtools}
\usepackage{dsfont} 
\usepackage{upgreek}
\usepackage{stmaryrd}
\usepackage[pdfpagelabels,hyperindex]{hyperref}
\hypersetup{linkbordercolor=green}

\newcommand*{\rom}[1]{\expandafter\@slowromancap\romannumeral #1@}
 \colorlet{lgray}{white!80!black}
\colorlet{lred}{white!85!red}
\colorlet{lgreen}{white!60!green}
\colorlet{dgreen}{black!30!green}
\colorlet{lpurple}{white!60!purple}
\colorlet{lblue}{white!60!blue}
\definecolor{green}{rgb}{0.1,0.8,0.1}
\definecolor{yellow}{rgb}{1.0,0.85,0.25}
\definecolor{purple}{rgb}{1.0, 0, 1.0}
\definecolor{blue}{rgb}{0, 0, 1.0}
\tikzstyle{unfused}=[lgray, line width=1.5pt, ->]
\tikzstyle{fused}=[lgray, line width=4pt, ->]
\tikzstyle{dual}=[black, line width=1pt, dashed]
\tikzstyle{lightdual}=[black, line width=0.5pt, dashed]
\tikzstyle{cut}=[black, line width=1.0pt]

\theoremstyle{plain}
\newtheorem{theorem}{Theorem}[section]
\newtheorem{lemma}[theorem]{Lemma}
\newtheorem{proposition}[theorem]{Proposition}

\theoremstyle{remark}
\newtheorem{definition}[theorem]{Definition}
\newtheorem{remark}[theorem]{Remark}

\colorlet{shadecolor}{gray!20}
\pgfplotsset{compat=1.9}
\usetikzlibrary{shapes.multipart}
\usetikzlibrary{patterns}
\usetikzlibrary{shapes.multipart}
\usetikzlibrary{arrows}
\usetikzlibrary{decorations.markings}
\usepgflibrary{decorations.shapes}
\usetikzlibrary{decorations.shapes}
\usepgflibrary{shapes.symbols}
\usetikzlibrary{shapes.symbols}
\usetikzlibrary{decorations.pathreplacing}
\tikzstyle{fleche}=[>=stealth', postaction={decorate}, thick]
\tikzstyle{axis}=[->, >=stealth', thick, gray]
\tikzstyle{paths}=[>->, >=stealth', thick]
\tikzstyle{path}=[->, >=stealth', thick]
\tikzstyle{grille}=[dotted, gray]
\DeclareFontFamily{U}{mathx}{}
\DeclareFontShape{U}{mathx}{m}{n}{<-> mathx10}{}
\DeclareSymbolFont{mathx}{U}{mathx}{m}{n}
\DeclareMathAccent{\widehat}{0}{mathx}{"70}
\DeclareMathAccent{\widecheck}{0}{mathx}{"71}

\usepackage{romannum}
\usepackage{upgreek}

\usepackage[scr=boondoxo]{mathalfa}

\newcommand{\kcr}{K^{\mathrm{hs}; \varpi}}

\newcommand{\icr}{I^{\mathrm{hs}; \varpi}}
\newcommand{\rcr}{R^{\mathrm{hs}; \varpi}}
\newcommand{\hsa}{\mathcal{A}^{\mathrm{hs}; \varpi}}
\newcommand{\hsai}{\mathcal{A}^{\mathrm{hs}; \infty}}
\newcommand{\hsail}{\mathcal{L}^{\mathrm{hs}; \infty}}

\newcommand{\hsm}{M^{\mathsf{S};\mathrm{hs};\varpi}}
\newcommand{\hsmi}{M^{\mathsf{S}; \mathrm{hs}; \infty}}

\newcommand{\hski}{K^{\mathrm{hs};\infty}}
\newcommand{\hsii}{I^{\mathrm{hs};\infty}}
\newcommand{\hsri}{R^{\mathrm{hs};\infty}}

\newcommand{\hsl}{\mathcal{L}^{\mathrm{hs}; \varpi}}

\newcommand{\ap}{\mathsf{a}}

\newcommand{\Ai}{\mathrm{Ai}}
\newcommand{\im}{\mathrm{i}} 

\usepgflibrary{fpu}
\makeatletter
\let\NAT@parse\undefined
\makeatother

\author{Jiaoyang Huang}
\address
{
Jiaoyang Huang, 
New York University 
}
\email{jh4427@nyu.edu} 
\author{Zongrui Yang}
\address
{
Zongrui Yang, 
University of California, Berkeley
}
\email{zongrui@berkeley.edu}

\title[Concentration of the half-space Airy line ensembles]{Concentration of the half-space Airy line ensembles}  
\begin{document}
\pagenumbering{arabic}

\begin{abstract}  
The half-space Airy line ensembles $\mathcal{A}^{\mathrm{hs};\varpi}$ ($\varpi\in\mathbb R\cup\{\infty\}$), introduced in \cite{dimitrov2026half,dimitrov2026pinned}, are natural analogues of the Airy line ensemble and are expected to arise as universal scaling limits of models in the half-space Kardar--Parisi--Zhang (KPZ) class. We prove that, at each fixed time $t\geq 0$, the height of the $k$th curve concentrates about $-(3\pi k/2)^{2/3}$ on the scale $k^{-1/3}\sqrt{\log(2+k)}$, with exponential tails.
\end{abstract}
\maketitle

\section{Introduction and the main theorem}

\subsection{Preface}

Over the past 20 years, the asymptotic behavior of half-space models in the Kardar--Parisi--Zhang (KPZ) universality class has received considerable attention. These models have local dynamics similar to those of their full-space counterparts, yet the interaction with a single boundary can change the scale and limiting law of their fluctuations. The half-space KPZ fixed point \cite{zhang2024tasep}, the half-space Airy line ensembles \cite{baik2018pfaffian,dimitrov2026half,dimitrov2026pinned}, and the half-space directed landscape \cite{dauvergne2026directed} are expected to arise as universal scaling limits in this class. They are the natural analogues of the full-space KPZ fixed point \cite{matetski2021kpz}, the Airy line ensemble \cite{prahofer2002scale,CorwinHammond2014}, and the directed landscape \cite{dauvergne2022directed}, respectively. 
For a non-exhaustive list of works on half-space KPZ models, we refer the reader to last-passage percolation (LPP) models \cite{BR01a, BR01b, BR01c,SI04,baik2018pfaffian,betea2018free,dimitrov2026half,das2026pinning,betea2020stationary,betea2022half,zeng2025stationary,occelli2026stationary}, log-gamma polymer \cite{o2014geometric,IMS22,barraquand2023identity,das2024half, BCD24,DasSerio25,barraquand2023stationary,zeng2026stationary}, stochastic six-vertex model and asymmetric simple exclusion process (ASEP) \cite{BBCW18,he2024boundary,he2026shift,barraquand2024markov,garbali2025symmetric,garbali2026transition}, and convergence to the KPZ equation \cite{Wu20,parekh2019positive,barraquand2023stationary,CS18, Par19,yang2025kpz,yang2026half}.

The study of Gibbsian line ensembles in half-space began with the log-gamma polymer \cite{BCD24}. The half-space Airy line ensembles $\hsa=\{\hsa_k\}_{k\geq 1}$ are ordered families of random continuous curves on $[0,\infty)$ indexed by a boundary parameter $\varpi\in\mathbb R\cup\{\infty\}$. The ensembles with finite $\varpi$ were introduced in \cite{dimitrov2026half}, and the pinned ensemble, corresponding to $\varpi=\infty$, was introduced in \cite{dimitrov2026pinned}. After a parabolic shift, these ensembles satisfy a Brownian Gibbs property. Away from the origin, the curves are locally described by avoiding Brownian bridges. Near the origin, the resampling law for finite $\varpi$ is given by avoiding reverse Brownian motions with alternating drifts determined by $\varpi$. In the pinned case, the boundary resampling law is instead described by Brownian motions pinned together in consecutive pairs at the origin.  
At finitely many positive times, the curves form a Pfaffian point process with the crossover kernel of \cite{baik2018pfaffian} for finite $\varpi$ and the symplectic--unitary kernel for $\varpi=\infty$. The latter kernel also appeared in \cite{ForresterNagaoHonner1999}.

A basic question for random matrices and point processes is how closely the location of a particle is determined by its rank. For the classical Airy line ensemble \cite{prahofer2002scale,CorwinHammond2014}, the joint law of the curves at finitely many times is described by the extended Airy kernel and is determinantal. At each fixed time, the height of the $k$th curve fluctuates about $-(3\pi k/2)^{2/3}$ on the scale $k^{-1/3}\sqrt{\log k}$ as $k\to\infty$ \cite[Theorem~3.3]{Hagg2007}. A concentration estimate, uniform over a bounded time interval \cite[Corollary~6.3]{DV2021}, enters the proof of the bridge representation of the Airy line ensemble, which approximates the curves by concatenating Brownian bridges between points on a fine time grid \cite{DV2021}. The bridge representation plays a key role in the construction of the directed landscape through last-passage percolation across the ensemble \cite{dauvergne2022directed}. For related results on counting fluctuations, tail estimates, and rigidity for the Airy point process, see \cite{Soshnikov2000,CorwinGhosal2020,ZhongAiry,BothnerBuckingham2018,ccGlobalRigidity}. Earlier work on fluctuations of particle counts in determinantal point processes includes \cite{costin1995gaussian,soshnikov1998level,forrester2014local}.

In this paper, we prove a concentration bound on the same scale for the half-space Airy line ensembles.
\begin{theorem}\label{thm:main}
Consider the half-space Airy line ensembles $\mathcal{A}^{\mathrm{hs};\varpi}$ for $\varpi\in\mathbb R\cup\{\infty\}$ introduced in \cite{dimitrov2026half,dimitrov2026pinned} (see Section \ref{subsec:def of half space Airy} below).
There exists an absolute constant $c>0$ with the following property. For
every $\varpi\in\mathbb R\cup\{\infty\}$ and $L>0$, there exists
$C_{\varpi,L}\in(0,\infty)$ such that, for every $k\in\mathbb N$,
$t\in[0,L]$, and $u\geq1$,
\[
 \mathbb P\!\left(
  \left|\mathcal A_k^{\mathrm{hs};\varpi}(t)
  +\left(\frac{3\pi k}{2}\right)^{2/3}\right|
  \geq u k^{-1/3}\sqrt{\log(2+k)}
 \right)
 \leq C_{\varpi,L}e^{-cu}.
\]  Moreover, there is an absolute
constant $C'\in(0,\infty)$ such that the same estimate, with $C'$ in
place of $C_{\varpi,L}$, holds simultaneously for all
$\varpi\in\mathbb R\cup\{\infty\}$, $t\geq0$, $k\geq2$, and $u\geq1$.
\end{theorem}   

The proof of Theorem~\ref{thm:main} will appear in Section~\ref{subsec:main inputs} and relies on several inputs. The subsequent sections establish these inputs. We mention that the concentration argument for the classical Airy line ensemble does not apply to the half-space setting because the determinantal point processes at finitely many times are replaced by Pfaffian point processes. Despite work on Janossy densities, fluctuations, and rigidity for Pfaffian point processes \cite{soshnikov2003janossy,kargin2014pfaffian,bufetov2018number,wang2025gaussian}, concentration in the Pfaffian setting remains less well understood than in the determinantal setting. We therefore develop a different argument. We first compare $\hsa$ for finite $\varpi$ with the pinned ensemble $\hsai$. An interlacing coupling of partitions sampled from Pfaffian Schur measures yields, in the Airy limit, a stochastic domination of the curve heights at each fixed time. Together with a separate upper-tail estimate for the top curve, this comparison reduces the problem to proving concentration for $\hsai$. For the pinned ensemble, we use the representation of its Pfaffian correlation kernel in terms of Airy-type functions and integral operators to obtain a tail bound of order $u^{-2}$. To strengthen this to an exponential bound, we use log-concavity, an approach inspired by \cite{wu2024applications}. At each fixed positive time, we realize the individual curve heights of the pinned ensemble as scaling limits of eigenvalues of matrix Ornstein--Uhlenbeck evolutions and use the log-concavity of the prelimit densities to prove that each curve of $\hsai$ at a fixed time has a log-concave law.

We expect that the concentration bound in Theorem~\ref{thm:main} can serve as the key input to establishing the Brownian bridge representation of \cite{DV2021} for the half-space Airy line ensembles. We also expect that the approach of \cite{dauvergne2022directed} can be adapted to realize the half-space directed landscape \cite{dauvergne2026directed} through last-passage percolation across the corresponding half-space Airy line ensemble. We leave these questions for future work.

\subsection{Definition of the half-space Airy line ensembles} 
\label{subsec:def of half space Airy}
In this subsection, we define the half-space Airy line ensembles $\hsa = \{\hsa_i\}_{i \geq 1}$ indexed by a boundary parameter $\varpi\in\mathbb R\cup\{\infty\}$, which were introduced in \cite{dimitrov2026half} in the case $\varpi\in\mathbb R$ and later in \cite{dimitrov2026pinned} in the case $\varpi=\infty$. 
The Pfaffian correlation kernels for these ensembles were introduced in \cite{baik2018pfaffian}. The symplectic--unitary transition corresponding to $\varpi=\infty$ was studied earlier in \cite{ForresterNagaoHonner1999}.

\begin{definition} \label{def:basic}
    We make the following conventions.
    \begin{enumerate}
        \item [(1)]   For a fixed $z \in \mathbb{C}$ and $\varphi \in (0, \pi)$, we denote by $\mathcal{C}_{z}^{\varphi}=\{z+|s|e^{\mathrm{sgn}(s)\im\varphi}: s\in \mathbb{R}\}$ the infinite contour oriented from $z+\infty e^{-\im\varphi}$ to $z+\infty e^{\im\varphi}$. 
        \item [(2)] For a finite set $\mathsf{S} = \{s_1, \dots, s_m\} \subset \mathbb{R}$, we let $\mu_{\mathsf{S}}$ denote the counting measure on $\mathbb{R}$, defined by $\mu_{\mathsf{S}}(A) = |A \cap \mathsf{S}|$. We also let $\mathrm{Leb}$ be the usual Lebesgue measure on $\mathbb{R}$ and $\mu_{\mathsf{S}} \times \mathrm{Leb}$ the product measure on $\mathbb{R}^2$.
        \item [(3)] Let $\mathcal{L} = \{\mathcal{L}_{i}\}_{i \geq 1}$ be a line ensemble on an interval $\Lambda \subseteq \mathbb{R}$. We say that $\mathcal{L}$ is {\em non-intersecting} if almost surely $\mathcal{L}_i(t) > \mathcal{L}_j(t)$ for all $t \in \Lambda$ and $1 \leq i < j$. We say that $\mathcal{L}$ is {\em ordered} if almost surely $\mathcal{L}_i(t) \geq \mathcal{L}_j(t)$ for all $t \in \Lambda$ and $1 \leq i < j$.   
    \end{enumerate}
\end{definition}

We define the Pfaffian correlation kernels $\kcr$ for the half-space Airy line ensembles $\hsa$, $\varpi\in\mathbb R$. 
\begin{definition} \label{def:kcr}
Fix $\varpi \in \mathbb{R}$, $m \in \mathbb{N}$ and $\mathsf{S} = \{s_1, \dots, s_m\}$, where $0 \leq s_1 < s_2 < \cdots < s_m$, and set $\ap_i = |\varpi| + 3i$ for $i = 1,2,3$. For $s,t \in \mathsf{S}$ and $x, y \in \mathbb{R}$, we define the kernel
\begin{equation}\label{Eq.S1DefKcross}
\begin{split}
&\kcr(s,x; t,y) = \begin{bmatrix}
    \kcr_{11}(s,x;t,y) & \kcr_{12}(s,x;t,y)\\
    \kcr_{21}(s,x;t,y) & \kcr_{22}(s,x;t,y) 
\end{bmatrix} \\
&= \begin{bmatrix}
    \icr_{11}(s,x;t,y) & \icr_{12}(s,x;t,y) + \rcr_{12}(s,x;t,y) \\
    -\icr_{12}(t,y;s,x) - \rcr_{12}(t,y;s,x) & \icr_{22}(s,x;t,y) + \rcr_{22}(s,x;t,y) 
\end{bmatrix} ,
\end{split}
\end{equation}
where the kernels $\icr_{ij}, \rcr_{ij}$ are defined as follows. We have
\begin{equation}\label{Eq.S1DefIcross}
\begin{split}
\icr_{11}(s,x;t,y) = &\frac{1}{(2\pi \im)^2} \int_{\mathcal{C}_{\ap_1 - s}^{\pi/3}}dz \int_{\mathcal{C}_{\ap_3 - t}^{\pi/3}} dw \frac{z + s - w - t}{(z + s + w + t)(z+s)(w+ t)} \\
& \times (z + s+ \varpi )(w + t + \varpi) \cdot e^{z^3/3 + w^3/3 - xz - y w}, \\
\icr_{12}(s,x;t,y) = &\frac{1}{(2\pi \im)^2} \int_{\mathcal{C}_{\ap_3 - s}^{\pi/3}}dz \int_{\mathcal{C}_{\ap_1 - t}^{2\pi/3}} dw \frac{z + s + w + t}{2(z+ s)(z + s - w -t)} \\
& \times \frac{z + \varpi + s}{w + \varpi + t} \cdot e^{z^3/3 - w^3/3 - xz + y w},\\
\icr_{22}(s,x;t,y) = &\frac{1}{(2\pi \im)^2} \int_{\mathcal{C}_{-\ap_2 - s}^{2\pi/3}}dz \int_{\mathcal{C}_{-\ap_2 - t}^{2\pi/3}} dw \frac{z + s - w - t}{4 (z + s + w + t)} \\
& \times \frac{1}{(z + s  + \varpi )(w + t + \varpi)} \cdot e^{-z^3/3 - w^3/3 + xz + y w},
\end{split}
\end{equation}
with contours as in Definition \ref{def:basic} \textup{(1)}. We also have \[\rcr_{12}(s,x;t,y) = 
 \frac{-{\bf1}{\{s<t\}}}{\sqrt{4\pi(t-s)}}\exp\left(
\dfrac{-(s-t)^4+6(x+y)(s-t)^2+3(x-y)^2}{12(s-t)}
\right),\]   $\rcr_{22}(s,x; t,y) = - \rcr_{22}(t,y; s,x) $, and when $x- s^2 > y -t^2$ we have
\begin{equation}\label{Eq.S1DefRcross2}
\begin{split}
\rcr_{22}(s,x;t,y) = & \frac{1}{2\pi \im} \int_{\mathcal{C}_{\ap_1 }^{2\pi/3}}dz \frac{e^{(-z + s)^3/3 + (\varpi + t)^3/3 - x (-z+s) - y (\varpi + t)}}{4(z - \varpi)}\\
&- \frac{1}{2\pi \im}\int_{\mathcal{C}_{-\ap_2 }^{2\pi/3}}dw \frac{e^{(-w + t)^3/3 + (\varpi + s)^3/3 - y (-w+t) - x (\varpi + s)}}{4(w - \varpi)}\\
& + \frac{{\bf 1}\{s + t > 0\} }{2\pi \im} \int_{\mathcal{C}^{2\pi/3}_{-\ap_{2}}} dw \frac{we^{ (-w + t)^3/3 + (w + s)^3/3   - y (-w + t) - x (w+s)}   }{2(w- \varpi)(w + \varpi) }.
\end{split}
\end{equation}
We set $\rcr_{22}(s,x;t,y)=0$ on $x-s^2=y-t^2$.
\end{definition}

We next recall the main theorem in \cite{dimitrov2026half}. See also the reformulation in \cite[Proposition 1.10]{dimitrov2026pinned}
\begin{proposition} \label{Prop.HSAiry}  Fix $\varpi \in \mathbb{R}$. There exists a line ensemble $\hsa = \{\hsa_i\}_{i \geq 1}$ on $[0, \infty)$ that satisfies the following properties. The ensemble is non-intersecting in the sense of Definition \ref{def:basic} \textup{(3)}. In addition, for any finite set $\mathsf{S} = \{s_1, \dots, s_m\} \subset [0,\infty)$ with $s_1 < s_2 < \cdots < s_m$, the random measure on $\mathbb{R}^2$ given by
\begin{equation}\label{eq:HSA point process}
\hsm(\omega, A) = \sum_{i \geq 1} \sum_{j = 1}^m {\bf 1}\{(s_j, \hsa_i(s_j,\omega) ) \in A\}
\end{equation}
is a Pfaffian point process with correlation kernel $\kcr$ as in Definition \ref{def:kcr} and reference measure $\mu_{\mathsf{S}} \times \mathrm{Leb}$ as in Definition \ref{def:basic} \textup{(2)}. 
Moreover, if $\hsl = \{\hsl_i\}_{i \geq 1}$ is the line ensemble defined through 
\begin{equation}\label{eq:Parabolic HSA}
\sqrt{2} \cdot \hsl_i(t) + t^2 = \hsa_i(t) \mbox{ for } i \geq 1, t \in [0,\infty),
\end{equation}
then $\hsl$ satisfies the half-space Brownian Gibbs property from \cite[Definition 2.7]{dimitrov2026half} with parameters $\mu_i = (-1)^{i} \sqrt{2} \cdot \varpi$.
\end{proposition}

We next define the Pfaffian correlation kernels $\hski$ for the pinned half-space Airy line ensemble $\hsai$.
\begin{definition}\label{def:kcrinfty}
    Fix $m \in \mathbb{N}$ and $\mathsf{S} = \{s_1, \dots, s_m\}$, where $0 < s_1 < s_2 < \cdots < s_m$. For $s,t \in \mathsf{S}$ and $x, y \in \mathbb{R}$, we define the kernel
    \begin{equation}\label{Eq.DefK}
\begin{split}
&\hski(s,x; t,y) = \begin{bmatrix}
    \hski_{11}(s,x;t,y) & \hski_{12}(s,x;t,y)\\
    \hski_{21}(s,x;t,y) & \hski_{22}(s,x;t,y) 
\end{bmatrix} \\
&= \begin{bmatrix}
    \hsii_{11}(s,x;t,y) & \hsii_{12}(s,x;t,y) + \hsri_{12}(s,x;t,y) \\
    -\hsii_{12}(t,y;s,x) - \hsri_{12}(t,y;s,x) & \hsii_{22}(s,x;t,y) + \hsri_{22}(s,x;t,y) 
\end{bmatrix},
\end{split}
\end{equation}
where the kernels $\hsii_{ij}, \hsri_{ij}$ are defined as follows. We have
\begin{equation}\label{Eq.DefII}
\begin{split}
\hsii_{11}(s,x;t,y) = &\frac{1}{(2\pi \im)^2} \int_{\mathcal{C}_{1+s}^{\pi/3}}dz \int_{\mathcal{C}_{1+t}^{\pi/3}} dw \frac{(z + s - w - t)H(z,x;w,y)}{4(z + s + w + t)(z+s)(w+ t)}, \\
\hsii_{12}(s,x;t,y) = &\frac{1}{(2\pi \im)^2} \int_{\mathcal{C}_{1+s}^{\pi/3}}dz \int_{\mathcal{C}_{1+t}^{\pi/3}} dw \frac{(z + s - w + t) H(z,x;w,y)}{2(z+ s)(z + s + w -t)}, \\
\hsii_{22}(s,x;t,y) = &\frac{1}{(2\pi \im)^2} \int_{\mathcal{C}_{1+ s}^{\pi/3}}dz \int_{\mathcal{C}_{1+t}^{\pi/3}} dw \frac{(z - s - w + t)H(z,x;w,y)}{z - s + w - t},
\end{split}
\end{equation}
where we have set 
\begin{equation}\label{Eq.DefH}
H(z,x;w,y) = e^{z^3/3 + w^3/3 - xz - y w},
\end{equation}
and the contours $\mathcal{C}_{z}^{\varphi}$ are as in Definition \ref{def:basic} \textup{(1)}. We also have that 
 \[R^{\mathrm{hs};\infty}_{12}(s,x;t,y) = 
 \frac{-{\bf1}{\{s<t\}}}{\sqrt{4\pi(t-s)}}\exp\left(
\dfrac{-(s-t)^4+6(x+y)(s-t)^2+3(x-y)^2}{12(s-t)}
\right),\]  
and
\begin{equation}\label{Eq.DefR22I}
\begin{split}
\hsri_{22}(s,x;t,y) = \frac{H(s,x;t,y) (y-t^2 -x + s^2) }{2\pi^{1/2} (t+s)^{3/2} } \cdot \exp \left(-\frac{(y - t^2 - x +s^2)^2}{4(t+s)}\right).
\end{split}
\end{equation}
\end{definition}

The next result is a combination of \cite[Theorems 1.13, 1.16, 1.20,  1.26, and Remarks  1.14, 1.15]{dimitrov2026pinned}.
\begin{proposition} \label{Prop.HSAiryinfty}
    Let $\hsa = \{\hsa_i\}_{i \geq 1}$ be as in Proposition \ref{Prop.HSAiry}. Then, the line ensembles $\hsa$ converge weakly to a line ensemble $\hsai$ as $\varpi \rightarrow \infty$.  This convergence is of random elements in $C(\mathbb{N} \times [0,\infty))$, where the latter space is endowed with the topology of uniform convergence over compact sets. The ensemble $\hsai$ is ordered on $[0,\infty)$ and non-intersecting on $(0,\infty)$ as in Definition \ref{def:basic} \textup{(3)}. At the origin, 
\begin{equation}\label{Eq.OriginCollision}
\hsai_1(0) = \hsai_2(0) > \hsai_3(0) = \hsai_4(0) > \hsai_5(0) = \hsai_6(0) > \cdots.
\end{equation}
Let  $\hsmi$ be as in (\ref{eq:HSA point process}) with $\varpi = \infty$ for $0< s_1 < s_{2} < \cdots < s_m$. Then, $\hsmi$ is a Pfaffian point process on $\mathbb{R}^2$ with reference measure $\mu_{\mathsf{S}} \times \mathrm{Leb}$ as in Definition \ref{def:basic} \textup{(2)} and correlation kernel $\hski$ as in Definition \ref{def:kcrinfty}. 
Furthermore,  let $\{\Lambda_k\}_{k \geq 0}$ denote, in ascending order, the eigenvalues of the stochastic Airy operator
  $\mathcal{H}_{\beta}$ from \cite{ramirez2011beta} with $\beta = 4$.  Then we have the following equality in law as random vectors in $\mathbb{R}^{\infty}$
\begin{equation}\label{Eq.Equality in law}
(\hsai_1(0), \hsai_2(0), \hsai_3(0), \hsai_4(0), \dots) \overset{d}{=}  (-2^{2/3}\Lambda_0, -2^{2/3}\Lambda_0, -2^{2/3}\Lambda_1, -2^{2/3}\Lambda_1, \dots).
\end{equation}
Moreover,  define $\hsail = \{\hsail_{i}\}_{i \geq 1}$ through
\begin{equation}\label{Eq.ParabolicHSALE}
\sqrt{2} \cdot \hsail_i(t) + t^2 = \hsai_i(t) \mbox{ for } i \geq 1, t \in [0, \infty).
\end{equation}
Then the restriction of $\hsail$ to $(0,\infty)$, viewed as an $\mathbb{N}$-indexed line ensemble on $(0,\infty)$, satisfies the Brownian Gibbs property from \cite[Definition 4.12]{dimitrov2026pinned}. On $[0,\infty)$ the ensemble $\hsail$ satisfies the pinned half-space Brownian Gibbs property from  \cite[Definition 4.10]{dimitrov2026pinned}.
\end{proposition}

\subsection{Main inputs and the proof of Theorem \ref{thm:main}}\label{subsec:main inputs}
In this subsection, we first state all the inputs needed to prove Theorem \ref{thm:main}, and then prove this theorem modulo these inputs.  The proofs of these inputs will be the focus of later sections of the paper.  

Suppose $X$ and $Y$ are real-valued random variables. We write $X\leq_{\mathrm{st}} Y$ or $Y\geq_{\mathrm{st}} X$ if $F_X(t)\geq F_Y(t)$ for each $t\in\mathbb R$, where $F_X$ and $F_Y$ are the cumulative distribution functions of $X$ and $Y$. 
\begin{proposition}\label{prop:domination}
    For each $\varpi\in\mathbb R$,  $t\geq0$, $k\in\mathbb N$, we have 
    $ \mathcal A_{k+1}^{\mathrm{hs};\infty}(t)\leq_{\mathrm{st}}\mathcal A_k^{\mathrm{hs};\varpi}(t)$
    and $ \mathcal A_{k+1}^{\mathrm{hs};\varpi}(t)\leq_{\mathrm{st}}\mathcal A_{k}^{\mathrm{hs};\infty}(t)$.
\end{proposition}

The proof of Proposition \ref{prop:domination} will be the focus of Section \ref{sec:comparison}.

\begin{proposition}\label{prop:bounding variance}
For \(t>0\) and \(r\geq1\), set
\[
 N_{t,r}=M^{t;\mathrm{hs};\infty}((-r,\infty))
 =\#\{k\in\mathbb N:\hsai_k(t)>-r\}.
\]
Then \(N_{t,r}<\infty\) almost surely and \(N_{t,r}\in L^2\).
Moreover, there exists an absolute constant \(C<\infty\) such that
\[
 \operatorname{Var}(N_{t,r})\leq C\log(2+r)
 \qquad\mbox{ for } t>0,\, r\geq1.
\]
\end{proposition}

We will prove Proposition \ref{prop:bounding variance}  in Section \ref{subsec:variance control}.

\begin{proposition}\label{prop:uniform upper tail}
    We have the following upper tail estimates of the first curve in the half-space Airy line ensembles:
    \begin{enumerate}
        \item [(1)]  For each $p>0$ there exists $C_p\in(0,\infty)$ depending only on $p$ such that 
    \[
     \mathbb P\left(\hsai_1(t)\geq x\right)\le C_p(1+x)^{-p}
 \qquad\mbox{ for } t>0,\,x\ge0. \]
 \item [(2)]    Fix $\varpi\in\mathbb R$ and $L\in(0,\infty)$. There exists  a constant $C_{\varpi,L}\in(0,\infty)$ such that
    \begin{equation*}
 \sup_{t\in[0,L]}
 \mathbb P\left(\mathcal A_1^{\mathrm{hs};\varpi}(t)\ge x\right)
 \le C_{\varpi,L}e^{-6x} 
 \qquad\mbox{ for }x\ge0. 
\end{equation*} 
    \end{enumerate}
\end{proposition} 

Proposition \ref{prop:uniform upper tail}  will be proved in Section \ref{subsec:uniform tail}.

We next define log-concavity of measures and random variables. As shown in Lemma \ref{lem:log-concavity improve concentration} below, this property improves a polynomial concentration bound of order $u^{-2}$ to an exponential bound.
\begin{definition}\label{def:log concave intro} A finite Borel measure $\mu$ on
$\mathbb R^m$ is called log-concave if 
\[
 \mu(\theta E+(1-\theta)F)
 \geq\mu(E)^\theta\mu(F)^{1-\theta}\qquad\mbox{ for all nonempty compact sets }E,F\mbox{ and }\theta\in(0,1).
\]
Here
$\theta E+(1-\theta)F=\{\theta x+(1-\theta)y:x\in E,\ y\in F\}$. 

A random vector is called log-concave when its law is log-concave. 
\end{definition}

\begin{lemma}\label{lem:log-concavity improve concentration} 
Let $Z$ have a log-concave law on $\mathbb R$.  If $a\in\mathbb R$,
$\sigma>0$, $C\ge1$, and
\[
 \mathbb P(|Z-a|\ge u\sigma)\le Cu^{-2}\qquad\mbox{ for each } u\ge1,
\]
then
\begin{equation}
 \mathbb P(|Z-a|\ge u\sigma)
 \le\sqrt3\exp\left(-\frac{\log3}{4\sqrt C}u\right)\qquad\mbox{ for each } u\ge1.
 \label{eq:logconcave-upgrade}
\end{equation} 
\end{lemma}

The proof of Lemma \ref{lem:log-concavity improve concentration} will appear in Section \ref{subsec:log concavity}.

\begin{proposition}\label{prop:log concave input}
    For each $k\in\mathbb N$, $t\geq0$, the random variable $\hsai_k(t)$ is log-concave.
\end{proposition}

The proof of Proposition \ref{prop:log concave input} will be the focus of Section \ref{sec: proof of log concavity input} and will appear in Section \ref{subsec: proof of log concavity input}.

We next prove Theorem \ref{thm:main} assuming the inputs stated above. 

\begin{proof}[Proof of Theorem \ref{thm:main}]
We begin by giving an outline of the proof. We first prove the result for the ensemble $\hsai$ (i.e., case $\varpi=\infty$). 
The mean of the random variable $\#\{k\in\mathbb N:\hsai_k(t)>-r\}$ can be estimated by controlling the integration of the Pfaffian kernel, written in terms of Airy functions. 
Combining this mean estimate with the variance bound given by Proposition \ref{prop:bounding variance} and the upper tail estimate of $\hsai_1(t)$ given by Proposition \ref{prop:uniform upper tail}(1), we obtain the $u^{-2}$-concentration of the curve $\hsai_k(t)$, which can then be upgraded to exponential concentration using Lemma \ref{lem:log-concavity improve concentration}. Then the theorem for $\hsa$, $\varpi\in\mathbb R$ can be obtained using the comparison with $\hsai$ as in Proposition \ref{prop:domination} and the upper tail bound for the first curve of $\hsa$ given by Proposition \ref{prop:uniform upper tail}(2).

Throughout the proof, we use $C,c\in(0,\infty)$ to denote absolute constants, and let 
\[
 \gamma_k=\left( 3\pi k/2\right)^{2/3},\qquad
 \ell_k=\log(2+k),\qquad
 s_k=k^{-1/3}\sqrt{\ell_k} .
\]
We first prove the theorem for the pinned half-space Airy line ensemble.  
For each $t>0$ and $r\geq1$, let $N_{t,r}=M^{t;\mathrm{hs};\infty}((-r,\infty))$ be as in
Proposition~\ref{prop:bounding variance}, and denote
$m_t(r)=\mathbb E[N_{t,r}]$.  
We next prove  
\begin{equation}\label{eq:main-count-estimates}
 \left|m_t(r)- 2/(3\pi) r^{3/2}\right|\leq C ,
 \qquad
 \operatorname{Var}(N_{t,r})\leq C \log(2+r) 
 \qquad\mbox{ for } t>0,\,r\geq1.
\end{equation}
Indeed, by \eqref{eq:alt12} from Lemma \ref{lem:alternative kernels}, \cite[(5.3) and (5.12)]{dimitrov2026half} and monotone convergence, 
\[
 m_t(r)=\int_{-r}^{\infty}K_{\Ai}(x,x)\,dx
 -\frac12\int_{-r}^{\infty}h_t(x)\Ai(x)\,dx.
\]
In view of \cite[(1.16), (1.17)]{Soshnikov2000},
\[
 K_{\Ai}(-s,-s)=\frac{s^{1/2}}{\pi}
 -\frac{\cos(4s^{3/2}/3)}{4\pi s}+O(s^{-5/2})
 \qquad\mbox{ as } s\to\infty.
\]
The oscillatory term above has a bounded primitive on $[1,\infty)$, as follows
from the change of variables $z=s^{3/2}$ and Dirichlet's test, while the
remainder is integrable. The integral of $K_{\Ai}(x,x)$ on $[-1,\infty)$ can be controlled using $K_{\Ai}(x,x)=\Ai'(x)^2-x\Ai(x)^2$ and the decay provided by
Lemma~\ref{lem:Airy decay}. Therefore,
\[
 \int_{-r}^{\infty}K_{\Ai}(x,x)\,dx
 =2/(3\pi) r^{3/2}+O(1)\qquad\mbox{ as }r\to\infty,
\]
where the error is bounded on all $r\geq1$. 
Using \eqref{eq:second line ht} in Lemma \ref{lem:ht bounds} and
Proposition~\ref{prop:bounding variance}, we conclude  \eqref{eq:main-count-estimates}.

We next prove that there is an absolute constant $C\geq1$ such that
\begin{equation}\label{eq:main-pinned-polynomial}
 \mathbb P\left(|\hsai_k(t)+\gamma_k|\geq us_k\right)
 \leq C u^{-2}
 \qquad\mbox{ for } t\geq0,\, k\in\mathbb N,\, u\geq1.
\end{equation}
 Indeed, it suffices to prove this for $t>0$ and $u\geq u_0$, where $u_0$ is a sufficiently large absolute constant.
 We first consider the lower tail.
 Denote $a=us_k/\gamma_k$ and $r_+=\gamma_k+us_k=\gamma_k(1+a)$. Since $\gamma_k\geq\gamma_1>1$, $r_+\geq1$. Moreover, $(1+a)^{3/2}-1\geq a\sqrt{1+a}$ and $ka= \left( 2/(3\pi)\right)^{2/3}u\sqrt{\ell_k}$. Thus \eqref{eq:main-count-estimates} implies, after fixing $u_0$ sufficiently large, 
\[
 m_t(r_+)-(k-1)
 \geq c u\sqrt{\ell_k}\sqrt{1+a}.
\] Furthermore, since $\ell_k\geq\log3>1$ and $\log(1+a)\leq a$,
\[  \log(2+r_+)
 \leq\log(2+\gamma_k)+\log(1+a)
 \leq C \ell_k(1+a).\]  
Since
$ \{\hsai_k(t)\leq-r_+\}=\{N_{t,r_+}\leq k-1\}$, by
Chebyshev's inequality, \eqref{eq:main-count-estimates}, and the last two displays, 
\begin{equation}\label{eq:main-pinned-lower-polynomial}
 \mathbb P(\hsai_k(t)\leq-\gamma_k-us_k)\leq C u^{-2}.
\end{equation}
We next consider the upper tail for \eqref{eq:main-pinned-polynomial}.  Let $r_-=\gamma_k-us_k$.
If $r_-\geq1$, then
$0\leq a<1$ and $1-(1-a)^{3/2}\geq a$. Thus, after increasing $u_0$
if necessary,
\[
 k-m_t(r_-)
 \geq k (1-(1-a)^{3/2} )-C \geq ka-C \geq \left( 2/(3\pi)\right)^{2/3}u\sqrt{\ell_k}-C 
 \geq c u\sqrt{\ell_k}.
\]
Also $\log(2+r_-)\leq\log(2+\gamma_k)\leq C\ell_k$. Since \(M^{t;\mathrm{hs};\infty}\) has reference measure $\mathrm{Leb}$,
\cite[(5.3)]{dimitrov2026half} gives
\(\mathbb E[M^{t;\mathrm{hs};\infty}(\{-r_-\})]=0\). Hence
$ \mathbb P(\hsai_k(t)\geq-r_-)
 =\mathbb P(N_{t,r_-}\geq k)$.
By Chebyshev's inequality and \eqref{eq:main-count-estimates},
\begin{equation}\label{eq:main-pinned-upper-regular}
 \mathbb P(\hsai_k(t)\geq-\gamma_k+us_k)\leq Cu^{-2}
 \qquad\text{if }r_-\geq1.
\end{equation}
Suppose instead that $r_-<1$.
Then
 $\{\hsai_k(t)\geq -r_-\}\subseteq\{N_{t,1}\geq k\}$,
and  by   $\sup_{t>0}\mathbb E[N_{t,1}^2]<\infty$ from \eqref{eq:main-count-estimates},
\begin{equation}\label{eq:main-exceptional-Q}
 \mathbb P(\hsai_k(t)\geq -r_-)\leq Ck^{-2}.
\end{equation}
If $-r_-\leq k^{2/3}$, then
$ u= (\gamma_k -r_-)/s_k
 \leq C k/\sqrt{\ell_k} \leq Ck$,
so \eqref{eq:main-exceptional-Q} is at most $Cu^{-2}$. If $-r_->k^{2/3}$, then, since
$\hsai_k(t)\leq\hsai_1(t)$ by orderedness,
Proposition \ref{prop:uniform upper tail} (1) with $p=3$ gives
\[
 \mathbb P(\hsai_k(t)\geq -r_-)\leq -Cr_-^{-3}.
\]
In this case $\gamma_k\leq C(-r_-)$ and hence
$ u\leq -C r_-k^{1/3}/\sqrt{\ell_k }$, and
$ u^2\leq C r_-^2k^{2/3}/\ell_k \leq -C r_-^3$;
thus $-r_-^{-3}\leq Cu^{-2}$. This proves the upper-tail bound also when
$r_-<1$,  and hence establishes \eqref{eq:main-pinned-polynomial} at positive times. The bound \eqref{eq:main-pinned-polynomial} can then be extended to $t=0$ in view of the almost sure continuity of the curve $\hsai_k(t)$. 

In view of Proposition \ref{prop:log concave input} and Lemma \ref{lem:log-concavity improve concentration},  the polynomial decay \eqref{eq:main-pinned-polynomial} can be strengthened  into exponential decay: for an absolute $c_{\mathrm{pin}}>0$,
\begin{equation}\label{eq:main-pinned-exponential}
 \mathbb P\left(|\hsai_k(t)+\gamma_k|\geq us_k\right)
 \leq\sqrt3e^{-c_{\mathrm{pin}
 }u} 
 \qquad\mbox{ for } t\geq0,\, k\in\mathbb N,\, u\geq1.
\end{equation} 

In the remainder of the proof, we consider $\hsa$ for finite $\varpi$.  We first observe that there is an absolute $K\geq1$
such that, for every $k\geq2$,
\begin{equation}\label{eq:main-neighbor-scales}
 K^{-1}s_k\leq s_{k-1},s_{k+1}\leq Ks_k,
 \qquad
 0\leq\gamma_{k+1}-\gamma_k,\ \gamma_k-\gamma_{k-1}\leq Ks_k.
\end{equation}
Indeed, for the first assertion, divide
$s_{k\pm1}=(k\pm1)^{-1/3}\sqrt{\log(k+2\pm1)}$ by $s_k$; the resulting
ratios are bounded above and away from zero uniformly for $k\geq2$,
since
\[
 1\leq\frac{k}{k-1}\leq2,\qquad
 \frac23\leq\frac{k}{k+1}\leq1,\qquad
 \frac12\leq\frac{\log(k+1)}{\log(k+2)}\leq1,
 \qquad
 1\leq\frac{\log(k+3)}{\log(k+2)}\leq\frac32.
\]
The second assertion follows from the mean-value theorem applied to
$x\mapsto(3\pi x/2)^{2/3}$, since
 $(k-1)^{-1/3}\leq2^{1/3}k^{-1/3}$,
and
$ k^{-1/3}\leq(\log4)^{-1/2}s_k$.
Fix $\varpi\in\mathbb R$. Proposition~\ref{prop:domination} gives
\[
 \hsai_{k+1}(t)\leq_{\mathrm{st}}
 \mathcal A_k^{\mathrm{hs};\varpi}(t)
 \leq_{\mathrm{st}}\hsai_{k-1}(t) 
 \qquad \mbox{ for }k\geq2,\, t\geq0.
\]Hence, for $u\geq2K$,
\begin{align*}
 \mathbb P\!\left(
  \mathcal A_k^{\mathrm{hs};\varpi}(t)\leq-\gamma_k-us_k
 \right)
 &\leq
 \mathbb P\!\left(
  |\hsai_{k+1}(t)+\gamma_{k+1}|
  \geq u/(2K)\cdot s_{k+1}
 \right),\\
 \mathbb P\!\left(
  \mathcal A_k^{\mathrm{hs};\varpi}(t)\geq-\gamma_k+us_k
 \right)
 &\leq
 \mathbb P\!\left(
  |\hsai_{k-1}(t)+\gamma_{k-1}|
  \geq u/(2K)\cdot s_{k-1}
 \right).
\end{align*}
By \eqref{eq:main-pinned-exponential} and the trivial probability bound
for $1\leq u<2K$, there are absolute constants $C_2,c_2>0$ such that
\begin{equation}\label{eq:main-finite-k-geq-two}
 \mathbb P\!\left(
  |\mathcal A_k^{\mathrm{hs};\varpi}(t)+\gamma_k|\geq us_k
 \right)
 \leq C_2e^{-c_2u}\qquad\mbox{ for all }\varpi\in\mathbb R,\,t\geq0,\,k\geq2,\,u\geq1.
\end{equation}
We next consider $k=1$. Set $u_-=\max (1, 2(\gamma_2-\gamma_1)/s_1, 2s_2/s_1 )$.
For $u\geq u_-$, by Proposition~\ref{prop:domination} and
\eqref{eq:main-pinned-exponential},
\begin{align*}
 \mathbb P\!\left(
  \mathcal A_1^{\mathrm{hs};\varpi}(t)\leq-\gamma_1-us_1
 \right)
 &\leq
 \mathbb P\!\left(
  |\hsai_2(t)+\gamma_2|\geq us_1/2 
 \right) \leq\sqrt{3}
 \exp\!\left(- c_{\mathrm{pin}}s_1/(2s_2)\cdot u\right).
\end{align*}
Thus this bound holds for every $u\geq1$ after replacing the prefactor 
by $C_-=\max (\sqrt{3},e^{c_-u_-})$, where $c_-= c_{\mathrm{pin}}s_1/(2s_2)$.
For the upper tail, put
 $u_+=\max\left(1, 2\gamma_1/s_1\right)$,
and let $D_{\varpi,L}$ be the constant from Proposition~\ref{prop:uniform upper tail}\textup{(2)}. If $u\geq u_+$, then
$-\gamma_1+us_1\geq us_1/2\geq0$, and that proposition gives
\[
 \sup_{t\in[0,L]}
 \mathbb P\!\left(
  \mathcal A_1^{\mathrm{hs};\varpi}(t)\geq-\gamma_1+us_1
 \right)
 \leq D_{\varpi,L}e^{-3s_1u}.
\]
Consequently this upper-tail bound holds for all $u\geq1$ with
$C_+e^{-c_+u}$ on the right, where $c_+=3s_1$ and $ C_+=\max (D_{\varpi,L},e^{c_+u_+})$. 
All of $c_{\mathrm{pin}},c_2,c_-,c_+$ are absolute and positive. Taking
$ c=\min(c_{\mathrm{pin}},c_2,c_-,c_+)$
and combining the last two bounds with
\eqref{eq:main-finite-k-geq-two} proves the assertion for finite
$\varpi$, for example with
$ C_{\varpi,L}=\max(C_2,C_-+C_+)$.
The case $\varpi=\infty$ is \eqref{eq:main-pinned-exponential}, and one
    may take $C_{\infty,L}=\sqrt{3}$. Finally,
 $C'=\max(\sqrt{3},C_2)$
proves the uniform assertion for $k\geq2$.
\end{proof}

\subsection{Outline of the rest of the paper} Section \ref{sec:comparison} proves Proposition~\ref{prop:domination} by constructing the coupling of partitions sampled from Pfaffian Schur measures and taking the Airy scaling limit. Section \ref{sec:3} proves Propositions \ref{prop:bounding variance} and \ref{prop:uniform upper tail} through an intricate analysis of the Pfaffian correlation kernels expressed as Airy-type functions and operators. Section \ref{sec: proof of log concavity input} proves the log-concavity (Proposition \ref{prop:log concave input}) by proving log-concavity of the spectrum of certain matrix Ornstein–Uhlenbeck evolutions and establishing their scaling convergence.

\subsection*{Acknowledgements} We thank Ivan Corwin for many helpful discussions throughout this project and, in particular, for helping us clarify the exposition of the coupling in Section~\ref{sec:comparison}. We also thank Milind Hegde and Leonid Petrov for helpful conversations. 
J.H. was supported partially by NSF Grant DMS-2337795 and a Sloan Research Fellowship.
Z.Y. was partially supported by the Miller Institute for Basic Research in Science at the University of California, Berkeley, Ivan Corwin's NSF grants DMS-1811143, DMS-2246576, Simons Foundation Grant 929852, and the Fernholz Foundation's ``Summer Minerva Fellows'' program.

\subsection*{AI declaration}
This project began in July 2024, and the proofs were written by the end of August 2026. ChatGPT 5.6 Sol suggested the monotone coupling in Proposition~\ref{prop:domination} and the use of log-concavity to improve the concentration bound in Lemma~\ref{lem:log-concavity improve concentration}. The authors also used AI tools to help edit the manuscript, perform some computations, and check arguments.

\section{Monotone couplings: Proof of Proposition \ref{prop:domination}}
\label{sec:comparison}

In this section, we prove Proposition~\ref{prop:domination}, which gives a stochastic domination between curves of half-space Airy line ensembles with different boundary parameters at each fixed time. We first construct an interlacing coupling of two partitions sampled from Pfaffian Schur measures \cite{borodin2005eynard,baik2018pfaffian}. We then take the Airy scaling limit of \cite{dimitrov2026half} to obtain the desired stochastic domination of continuous ensembles.

A {\em partition} is a non-increasing sequence of non-negative integers $\lambda=(\lambda_1\geq\lambda_2\geq\dots\geq\lambda_k \geq \cdots)$ with finitely many non-zero elements. We denote the space of partitions by $\mathbb Y$. For any partition $\lambda\in\mathbb Y$, we denote its {\em weight} by $\vert\lambda\vert=\sum_{i=1}^{\infty}\lambda_i$. There is a single partition of weight $0$, which we denote by $\emptyset$.  Given two partitions $\lambda$ and $\mu$, we write $\mu\preceq\lambda$ or $\lambda \succeq \mu$ and say that $\mu$ and $\lambda$ {\em interlace} if $\lambda_1\geq\mu_1\geq\lambda_2\geq\mu_2\geq\dots$.
For a partition $\lambda\in\mathbb Y$, its dual, denoted $\lambda'$, is the partition given by $\lambda'_i=\#\{j:\lambda_j\geq i\}$. We say that a partition $\lambda$ is even if all its components are even integers. Note that $\lambda'$ is even if and only if $\lambda_{2i-1}=\lambda_{2i}$ for each $i\in\mathbb N$.
Given finitely many variables $x_1, \dots, x_n$, we define the {\em skew Schur polynomials} via
\begin{equation}\label{def: skew Schur}
s_{\lambda/ \mu}(x_1, \dots, x_n) = \sum_{\mu = \lambda^{0} \preceq  \lambda^{1} \preceq \cdots \preceq \lambda^{n} = \lambda} \prod_{i = 1}^n x_i^{|\lambda^{i}| - |\lambda^{i-1}|}.
\end{equation}
When $\mu = \emptyset$ we drop it from the notation, and write $s_{\lambda}$, which is then the Schur polynomial indexed by $\lambda$. We also define the {\em boundary monomial} in a single variable $c$ by
\begin{equation}\label{def:DefTau}
\tau_\lambda(c) =c^{\sum_{j=1}^{\infty}(-1)^{j-1}\lambda_j}=c^{\lambda_1-\lambda_2+\lambda_3-\lambda_4+\dots}.
\end{equation}
At $c=0$, the right side is understood with $0^0=1$, so that 
$\tau_\lambda(0)=\mathbf 1_{\{\lambda'\text{ is even}\}}$.

\begin{lemma}\label{lem:interlace}
    Fix $M\in\mathbb Z_{\geq0}$, $N\in\mathbb N$, $q\in(0,1)$, and $c\in[0,q^{-1})$.  Write $q^r=(q,\ldots,q)$ for the specialization
consisting of $r$ copies of $q$.  Then the following is a probability measure 
\[
\mathsf P_c^{M,N}(\lambda)
 =   \frac{1}{Z_{M,N}(c)}    s_\lambda(q^N) 
   \sum_{\rho\in\mathbb Y} \tau_\rho(c)s_{\lambda/\rho}(q^M),\qquad\lambda\in\mathbb Y,
\]
where $ Z_{M,N}(c)=(1-cq)^{-N}
 (1-q^2)^{-MN-N(N-1)/2}$ is the normalization constant.
Suppose $\lambda\sim\mathsf P_c^{M,N}$ and $\nu\sim\mathsf P_0^{M,N}$. Then there exists a coupling between random partitions $\nu$ and $\lambda$, such that $\nu\preceq\lambda$ almost surely.
\end{lemma}
\begin{remark}
   As a corollary, the joint law of $(\lambda_2,\lambda_4,\ldots)$ under $\mathsf P_c^{0,N}$ is independent of $c$. Indeed, for $\nu\sim\mathsf P_0^{0,N}$, we have $\nu_{2j}=\nu_{2j-1}$ for $j\in\mathbb N$. Using the coupling given by Lemma \ref{lem:interlace}, we have $\nu_{2j-1}\geq\lambda_{2j}\geq\nu_{2j}=\nu_{2j-1}$. Thus all the even-indexed parts of $\lambda\sim\mathsf P_c^{0,N}$ agree simultaneously with those of $\nu\sim\mathsf P_0^{0,N}$.
\end{remark}
\begin{proof}  Define the nonnegative measure
    \[
     \mathsf Q_c^{M,N}(\nu,\lambda)
 =\frac{1}{Z_{M,N}(c)}
 U_M(\nu)s_{\lambda/\nu}(c)s_\lambda(q^N),\qquad\nu,\lambda\in\mathbb Y,
    \]
    where 
     $U_M(\nu)=\sum_{\mu\in\mathbb Y:\,\mu'\text{ is even}}s_{\nu/\mu}(q^M)$. 
The skew Cauchy identity \cite[(7)]{baik2018pfaffian} and the
Littlewood identity \cite[(9)]{baik2018pfaffian} give
\begin{align*}
\sum_{\nu,\lambda\in\mathbb Y}\mathsf Q_c^{M,N}(\nu,\lambda)
&=\frac{(1-cq)^{-N}}{Z_{M,N}(c)}
  \sum_{\nu\in\mathbb Y}U_M(\nu)s_\nu(q^N) =\frac{(1-cq)^{-N}(1-q^2)^{-MN}}{Z_{M,N}(c)}
  \sum_{ \mu\in\mathbb Y,\, \mu'\text{ is even}} s_\mu(q^N) =1.
\end{align*}
Thus $\mathsf Q_c^{M,N}$ is a probability measure.
Here we use Tonelli's theorem since all summands above are nonnegative;  and since $cq<1$ and $q^2<1$, the specialized Cauchy and Littlewood products are finite.

   Suppose $(\nu,\lambda)\sim\mathsf Q_c^{M,N}$, we will show that $\lambda\sim\mathsf P_c^{M,N}$ and $\nu\sim\mathsf P_0^{M,N}$. Indeed, by \eqref{def: skew Schur}, $s_{\rho/\mu}(c)$ vanishes unless
$\mu\preceq\rho$, and, when $\mu\preceq\rho$, it equals
$c^{|\rho|-|\mu|}$. If also $\mu'$ is even, then
$\mu_{2i-1}=\mu_{2i}$, and $\mu\preceq\rho$ forces $\mu_{2i-1}=\mu_{2i}=\rho_{2i}$ for $i\in\mathbb N$. 
Conversely, the partition
$\mu^{(\rho)}=(\rho_2,\rho_2,\rho_4,\rho_4,\ldots)$ has even conjugate
and interlaces $\rho$. Consequently, 
\[
\sum_{\mu\in\mathbb Y:\,\mu'\text{ is even}}s_{\rho/\mu}(c)
=s_{\rho/\mu^{(\rho)}}(c)
=c^{\sum_{i\geq1}(\rho_{2i-1}-\rho_{2i})}
=\tau_\rho(c).
\]
Using this identity and the branching rule \cite[(8)]{baik2018pfaffian}, the claim regarding
the $\lambda$-marginal follows from
\[ 
 \sum_{\rho\in\mathbb Y}\tau_\rho(c)s_{\lambda/\rho}(q^M)
  =\sum_{ \rho,\mu\in\mathbb Y,\,\mu'\text{ even}} 
       s_{\rho/\mu}(c)s_{\lambda/\rho}(q^M) 
  =\sum_{\mu\in\mathbb Y:\,\mu'\text{ even}}s_{\lambda/\mu}(c,q^M) 
  =\sum_{\nu\in\mathbb Y} U_M(\nu)s_{\lambda/\nu}(c). \]
For the $\nu$-marginal, the skew Cauchy identity gives
\[
\sum_{\lambda\in\mathbb Y}s_{\lambda/\nu}(c)s_\lambda(q^N)
 =H(c;q^N)s_\nu(q^N),
 \qquad
 H(c;q^N)=\prod_{j=1}^N(1-cq)^{-1}=(1-cq)^{-N}.
\]
Since $Z_{M,N}(c)=H(c;q^N)Z_{M,N}(0)$ and
$\tau_\rho(0)=\mathbf 1_{\{\rho'\text{ is even}\}}$, it follows that
\begin{align*}
\sum_{\lambda\in\mathbb Y}\mathsf Q_c^{M,N}(\nu,\lambda)
&=\frac{H(c;q^N)}{Z_{M,N}(c)}
  U_M(\nu)s_\nu(q^N) =\frac{1}{Z_{M,N}(0)}s_\nu(q^N)
  \sum_{\rho\in\mathbb Y}\tau_\rho(0)s_{\nu/\rho}(q^M)
 =\mathsf P_0^{M,N}(\nu).
\end{align*}
Since $\nu\preceq\lambda$ almost surely under $(\nu,\lambda)\sim\mathsf Q_c^{M,N}$, we conclude the proof.
\end{proof}

\begin{remark}
    The coupling $\mathsf Q_c^{M,N}$ constructed in the proof of Lemma \ref{lem:interlace} can be written as 
    $$\mathsf Q_c^{M,N}(\nu,\lambda)=\mathsf P_0^{M,N}(\nu)\cdot \mathcal{U}^{\llcorner}_{q^N,(c)}(\lambda | \nu,\emptyset),$$
    where $\mathcal{U}^{\llcorner}$, defined by \cite[(14)]{baik2018pfaffian}, is a transition operator between Pfaffian Schur processes along down-right paths in the half-quadrant.
    This gives the coupling $\mathsf Q_c^{M,N}$ an interpretation in the half-space (symmetric) geometric last-passage percolation (LPP) model: starting with the half-space geometric LPP with zero weights along the diagonal and entries below the diagonal having i.i.d.    weights with distribution $\operatorname{Geom}(q^2)$, we   append a column at $(M+N+1,j)$, $1\leq j\leq N$, consisting of i.i.d.  weights with distribution $\operatorname{Geom}(cq)$. Then the joint law of the partitions $\nu$ at $(M+N,N)$ and $\lambda$ at $(M+N+1,N)$ is given by $\mathsf Q_c^{M,N}$. 
    
    We also mention that the monotone coupling of partitions constructed in Lemma \ref{lem:interlace} can be generalized to the half-space Macdonald processes \cite{barraquand2020half}, which may admit future applications.
\end{remark}

\begin{proof}[Proof of Proposition \ref{prop:domination}]
Fix $q\in(0,1)$ and denote
\[
 \sigma_q=\frac{q^{1/3}(1+q)^{1/3}}{1-q},\qquad
 f_q=\frac{q^{1/3}}{2(1+q)^{2/3}},\qquad
 h_q=\frac{2q}{1-q},\qquad p_q=\frac{q}{1-q}.
\]
For $t\geq0$ let $M_N=\lfloor f_q^{-1}tN^{2/3}\rfloor$.  For each $\theta\in\mathbb R$, choose $c_N(\theta)\in[0,q^{-1})$ so that
\[ 
 c_N(\theta)=1-\theta\sigma_q^{-1}N^{-1/3}+o(N^{-1/3}).\]
By Lemma \ref{lem:interlace}, for each  $N\in\mathbb N$, there exists a coupling 
\[
 \lambda^{N,\theta}\sim\mathsf P_{c_N(\theta)}^{M_N,N},
 \qquad
 \nu^{N}\sim\mathsf P_0^{M_N,N},
 \qquad
 \nu^{N}\preceq\lambda^{N,\theta}.
\]
Let 
\begin{align*}
    X_i^{N,\theta}(t)&=\sigma_q^{-1}N^{-1/3}
 \left(\lambda_i^{N,\theta}-h_qN-p_qM_N-i\right)+t^2,\\
  V_i^{N}(t)&=\sigma_q^{-1}N^{-1/3}
 \left(\nu_i^{N}-h_qN-p_qM_N-i\right)+t^2.
\end{align*}
For all sufficiently large $N$, we have
$c_N(\theta)\in(q,q^{-1})$ and $M_N\leq N$. Let
$(\Lambda^0,\ldots,\Lambda^N)$ have the Pfaffian Schur process from
\cite[Definition~5.16]{dimitrov2026half} with
$a_1=\cdots=a_N=b_1=\cdots=b_N=q$ and $c=c_N(\theta)$.
The branching and skew Cauchy identities applied to
\cite[(5.20)]{dimitrov2026half} give
\[
 \mathbb P(\Lambda^{M_N}=\lambda)
 =\frac{1}{Z_{M_N,N}(c_N(\theta))}
 s_\lambda(q^N)\sum_{\rho\in\mathbb Y}
 \tau_\rho(c_N(\theta))s_{\lambda/\rho}(q^{M_N}).
\]
Indeed, the branching and skew Cauchy identities show that summing over the partitions after time $M_N$ gives
$ (1-q^2)^{-(N-M_N)N}s_\lambda(q^N)$,
whereas summing over the partitions before time $M_N$ gives the displayed sum. The factor $(1-q^2)^{-(N-M_N)N}$ is independent of $\lambda$ and is absorbed into the normalization, since
$Z_{N,N}(c_N(\theta))
 =(1-q^2)^{-(N-M_N)N}Z_{M_N,N}(c_N(\theta))$.
  Thus
$\Lambda^{M_N}$ has the same law as $\lambda^{N,\theta}$. Let $s=f_q^{-1}t$, so that the quantity
$T_s(N)$ in \cite[Definition~5.20]{dimitrov2026half} equals $M_N$.
The variable in \cite[(5.25)]{dimitrov2026half} at time $s$ is then
$X_i^{N,\theta}(t)-t^2$. Hence
\cite[Proposition~5.21 and the proof of Theorem~1.4, Steps~2--3,
in particular (5.33)]{dimitrov2026half} yields, for each fixed
$i\in\mathbb N$ and $t\geq0$,
\begin{equation}\label{eq:hbrteyhn}
    X_i^{N,\theta}(t)\Rightarrow
    \mathcal A_i^{\mathrm{hs};\theta}(t)
    \qquad\mbox{as }N\to\infty.
\end{equation}
 For every $i\in\mathbb N$,
\begin{align*}
 V_i^{N}(t)-X_{i+1}^{N,\theta}(t)
 =\sigma_q^{-1}N^{-1/3}
   \left(\nu_i^{N}-\lambda_{i+1}^{N,\theta}+1\right)>0,\qquad
 X_i^{N,\theta}(t)-V_i^{N}(t)
 =\sigma_q^{-1}N^{-1/3}
   \left(\lambda_i^{N,\theta}-\nu_i^{N}\right)\geq0.
\end{align*}
Consequently,
\begin{equation*}
 X_{i+1}^{N,\theta}(t)\leq V_i^{N}(t)\leq X_i^{N,\theta}(t). 
\end{equation*}
Fix $\theta,\varpi\in\mathbb R$ and a bounded, continuous,
nondecreasing $\phi:\mathbb R\to\mathbb R$. Apply the above
display separately with parameter $\eta=\theta$ and with
parameter $\eta=\varpi$, and write $V_k^{N,\eta}(t)$ for the
corresponding intermediate variable. These two couplings need not be
simultaneous, but $V_k^{N,\theta}(t)$ and $V_k^{N,\varpi}(t)$ have the
same law because the underlying partitions both have law
$\mathsf P_0^{M_N,N}$. Hence, for each $k\in\mathbb N$,
\begin{align*}
 \mathbb E\phi(X_{k+1}^{N,\theta}(t))
 &\leq \mathbb E\phi(V_k^{N,\theta}(t))
  =\mathbb E\phi(V_k^{N,\varpi}(t))
  \leq \mathbb E\phi(X_k^{N,\varpi}(t)),\\
 \mathbb E\phi(X_{k+1}^{N,\varpi}(t))
 &\leq \mathbb E\phi(V_k^{N,\varpi}(t))
  =\mathbb E\phi(V_k^{N,\theta}(t))
  \leq \mathbb E\phi(X_k^{N,\theta}(t)).
\end{align*}
Letting $N\to\infty$ and using \eqref{eq:hbrteyhn}, for each $k\in\mathbb N$,
\begin{align*}
 \mathbb E\phi(\mathcal A_{k+1}^{\mathrm{hs};\theta}(t))
  \leq\mathbb E\phi(\mathcal A_k^{\mathrm{hs};\varpi}(t)),
  \qquad
 \mathbb E\phi(\mathcal A_{k+1}^{\mathrm{hs};\varpi}(t))
 \leq\mathbb E\phi(\mathcal A_{k}^{\mathrm{hs};\theta}(t)).
\end{align*}
By Proposition \ref{Prop.HSAiryinfty},
$\mathcal A^{\mathrm{hs};\theta}\Rightarrow
\mathcal A^{\mathrm{hs};\infty}$ in
$C(\mathbb N\times[0,\infty))$ as $\theta\to\infty$.
Since coordinate evaluation at $(i,t)$ is continuous in this topology,
including when $t=0$, letting $\theta\to\infty$ gives
\[
\mathbb E\phi(\mathcal A_{k+1}^{\mathrm{hs};\infty}(t))
\leq\mathbb E\phi(\mathcal A_k^{\mathrm{hs};\varpi}(t)),
\qquad
\mathbb E\phi(\mathcal A_{k+1}^{\mathrm{hs};\varpi}(t))
\leq\mathbb E\phi(\mathcal A_k^{\mathrm{hs};\infty}(t)).
\]
For each $a\in\mathbb R$, apply these to the bounded,
continuous, nondecreasing functions
$\phi_{n,a}(x)= (n(x-a)_+ )\wedge1$.
By bounded convergence as $n\to\infty$, 
\[
\mathbb P\left(\mathcal A_{k+1}^{\mathrm{hs};\infty}(t)>a\right)
\leq
\mathbb P\left(\mathcal A_k^{\mathrm{hs};\varpi}(t)>a\right),
\qquad
\mathbb P\left(\mathcal A_{k+1}^{\mathrm{hs};\varpi}(t)>a\right)
\leq
\mathbb P\left(\mathcal A_k^{\mathrm{hs};\infty}(t)>a\right).
\]
We conclude the proof.
\end{proof}

\section{Analysis of the fixed-time Pfaffian point processes: Proofs of Propositions \ref{prop:bounding variance} and \ref{prop:uniform upper tail} } \label{sec:3}

We begin in Section~\ref{subsec:3.1} by collecting estimates for the Airy function and identities for the associated integral operators. In Section~\ref{subsec:3.2}, we derive alternative formulas for the entries of the Pfaffian kernel $\hski$ at equal times $s=t>0$ in terms of Airy-type integrals. In Section~\ref{subsec:uniform tail}, we use estimates of the one-point densities to prove the upper tail bounds for the top curves of the ensembles given by Proposition~\ref{prop:uniform upper tail}. Finally, in Section~\ref{subsec:variance control} we prove the bound of the variance of $\#\{k\in\mathbb N:\hsai_k(t)>-r\}$ given by Proposition~\ref{prop:bounding variance}.

\subsection{Properties of the Airy functions and operators} \label{subsec:3.1} In this subsection, we collect a few properties of the Airy functions and operators which will be useful in the later subsections.

\begin{lemma}[cf. \cite{DLMF} 9.7.5--9.7.10]\label{lem:Airy decay}
    Let $\Ai(z)$ denote the Airy function. Then there is an absolute constant $C\in(0,\infty)$ such that, for $s\ge0$,
   \[  \begin{aligned}
 |\Ai(-s)|&\le C(1+s)^{-1/4},&
 |\Ai'(-s)|&\le C(1+s)^{1/4},\\
 |\Ai(s)|&\le C(1+s)^{-1/4}e^{-2s^{3/2}/3},&
 |\Ai'(s)|&\le C(1+s)^{1/4}e^{-2s^{3/2}/3}.
\end{aligned}\]
Furthermore, 
\[
\left|\Ai(-s)- \pi^{-1/2}s^{-1/4}\sin\left(\frac23s^{3/2}+\frac\pi4\right)\right|\leq Cs^{-7/4}\qquad\mbox{ for }s\geq1.
\] 
\end{lemma}

\begin{definition}\label{def:Airy kernel}
For $x,y\in\mathbb R$, we define the Airy kernel
\[
 K_{\Ai}(x,y)=\int_0^\infty\Ai(x+s)\Ai(y+s)\,ds.
\]
It is well-known (see, for example, \cite[(2)]{BasorWidom1999}) that 
\begin{equation}\label{eq:Darboux KAi}
  K_{\Ai}(x,y)
 =\frac{\Ai(x)\Ai'(y)-\Ai'(x)\Ai(y)}{x-y}\qquad\mbox{ when }x\ne y.
\end{equation}
\end{definition}

\begin{lemma}\label{lem:operators Airy}
 For $f\in C_c^\infty(\mathbb R_+)$, define
\[
 (\mathsf A_-f)(x)=\int_0^\infty\Ai(x+s)f(s)\,ds,
 \qquad
 (\mathsf A_+f)(x)=\int_0^\infty\Ai(x-v)f(v)\,dv,
 \qquad x\in\mathbb R.
\]
They extend uniquely to linear isometries $ \mathsf A_-,\mathsf A_+:L^2(\mathbb R_+)\to L^2(\mathbb R)$. Their ranges give an orthogonal decomposition $L^2(\mathbb R)
 =\operatorname{Ran}(\mathsf A_-) \oplus \operatorname{Ran}(\mathsf A_+)$. Moreover, the integral operator with kernel $K_{\Ai}$ from
Definition~\ref{def:Airy kernel}, defined on $C_c^\infty(\mathbb R)$ by $(K_{\Ai}g)(x)=\int_{\mathbb R}K_{\Ai}(x,y)g(y)\,dy$, has a unique bounded extension to $L^2(\mathbb R)$,  and we have
\[
 \mathsf A_-^*\mathsf A_-
 =\mathsf A_+^*\mathsf A_+
 =I_{L^2(\mathbb R_+)} ,
 \qquad
 \mathsf A_-^*\mathsf A_+
 =\mathsf A_+^*\mathsf A_-=0,\qquad K_{\Ai}=\mathsf A_-\mathsf A_-^*,
 \qquad
 \mathsf A_-\mathsf A_-^*+\mathsf A_+\mathsf A_+^*
 =I_{L^2(\mathbb R)} .
\]In particular, $K_{\Ai}$ is the orthogonal projection onto
$\operatorname{Ran}(\mathsf A_-)$.  All operator equalities above are equalities of bounded operators on
the indicated $L^2$-spaces.
\end{lemma}

\begin{proof}
    Let $\mathcal A$ be the Airy transform on $L^2(\mathbb R)$ from \cite[Section 2]{BasorWidom1999}. By \cite[Lemma 2.1 and the preceding paragraph]{BasorWidom1999}, we have $\mathcal A^*=\mathcal A$, $\mathcal A^2=I_{L^2(\mathbb R)}$, and  for each $h\in L^1(\mathbb R)\cap L^2(\mathbb R)$,
\[
 (\mathcal A h)(x)=\int_{\mathbb R}\Ai(x+y)h(y)\,dy.
\]
 Define isometries $\iota_-,\iota_+:L^2(\mathbb R_+)\to L^2(\mathbb R)$ by
\[
 (\iota_-f)(y)=\mathbf 1_{(0,\infty)}(y)f(y),
 \qquad
 (\iota_+f)(y)=\mathbf 1_{(-\infty,0)}(y)f(-y).
\]For each $f\in C_c^\infty(\mathbb R_+)$, the integral representation of
$\mathcal A$ and the change of variables $y=-v$ give
\[
 \mathsf A_-f=\mathcal A\iota_-f,
 \qquad
 \mathsf A_+f=\mathcal A\iota_+f.
\]
Hence $\mathsf A_-$ and $\mathsf A_+$ take values in $L^2(\mathbb R)$ and are isometries there.  Since
$C_c^\infty(\mathbb R_+)$ is dense in $L^2(\mathbb R_+)$, their continuous extensions exist and are unique. 
 The ranges of $\iota_-$ and $\iota_+$ are the mutually orthogonal subspaces
of functions supported, respectively, on $(0,\infty)$ and $(-\infty,0)$;
their orthogonal direct sum is $L^2(\mathbb R)$.  Hence
\[
 \iota_-^*\iota_-=\iota_+^*\iota_+=I_{L^2(\mathbb R_+)},
 \qquad
 \iota_-^*\iota_+=\iota_+^*\iota_-=0,
 \qquad
 \iota_-\iota_-^*+\iota_+\iota_+^*=I_{L^2(\mathbb R)}.
\]
Combining these identities with
$\mathsf A_-=\mathcal A\iota_-$,
$\mathsf A_+=\mathcal A\iota_+$, and $\mathcal A^*=\mathcal A$, $\mathcal A^2=I_{L^2(\mathbb R)}$, we conclude the claimed identities  except for $K_{\Ai}=\mathsf A_-\mathsf A_-^*$. For $\phi,\psi\in C_c^\infty(\mathbb R)$, using $\mathsf A_-^*=\iota_-^*\mathcal A$, we have
\begin{equation}\label{eq:AGRE}\begin{aligned}
 &\langle\phi,\mathsf A_-\mathsf A_-^*\psi\rangle_{L^2(\mathbb R)}
 =\int_0^\infty
   \overline{\left(\int_{\mathbb R}\Ai(x+s)\phi(x)\,dx\right)}
   \left(\int_{\mathbb R}\Ai(y+s)\psi(y)\,dy\right)ds\\
 &=\int_{\mathbb R}\int_{\mathbb R}
   \overline{\phi(x)}
   \left(\int_0^\infty\Ai(x+s)\Ai(y+s)\,ds\right)
   \psi(y)\,dy\,dx =\int_{\mathbb R}\int_{\mathbb R}
   \overline{\phi(x)}K_{\Ai}(x,y)\psi(y)\,dy\,dx.
\end{aligned}   \end{equation}
Here, the use of Fubini's theorem is justified because the decay of Airy functions as in Lemma~\ref{lem:Airy decay} gives
\[
 \int_0^\infty\int_{\mathbb R}\int_{\mathbb R}
 |\phi(x)\psi(y)\Ai(x+s)\Ai(y+s)|\,dy\,dx\,ds<\infty.
\] 
For  $\psi\in C_c^\infty(\mathbb R)$, set  $F_\psi(x)=\int_{\mathbb R}K_{\Ai}(x,y)\psi(y)\,dy$.
By the bound above,  $F_\psi\in L^1_{\mathrm{loc}}(\mathbb R)$.
By \eqref{eq:AGRE},  $F_\psi=\mathsf A_-\mathsf A_-^*\psi$
as distributions.  Since 
$\mathsf A_-\mathsf A_-^*\psi\in L^2(\mathbb R)$, this equality holds almost everywhere and
$F_\psi\in L^2(\mathbb R)$.  Hence the integral operator with kernel
$K_{\Ai}$ agrees on $C_c^\infty(\mathbb R)$ with the bounded operator
$\mathsf A_-\mathsf A_-^*$.  Since $C_c^\infty(\mathbb R)$ is dense in
$L^2(\mathbb R)$, it has the asserted unique bounded extension, and
that extension equals $\mathsf A_-\mathsf A_-^*$.
\end{proof}

\begin{lemma}\label{lem:Airy leakage}
Define, for each $u,v\geq0$ and $r\geq1$,
\[
 G_r(u,v) =\int_{-r}^\infty\Ai(x+u)\Ai(x-v)\,dx.
\] 
Then there is an absolute constant  $C\in(0,\infty)$ such that  
\begin{equation*}
 \int_0^\infty\!\!\int_0^\infty |G_r(u,v)|^2\,du\,dv
 \le C\log(2+r)\qquad\mbox{ for }r\geq1. 
\end{equation*} 
\end{lemma}
\begin{proof}
By Lemma \ref{lem:Airy decay}, for each nonempty compact set $Q\subset[0,\infty)^2$, 
\[
 \sup_{(u,v)\in Q}|\Ai(x+u)\Ai(x-v)|
 \leq C_1e^{-c_1x^{3/2}},\qquad x\geq2\max\{v:(u,v)\in Q\}+2,
\]
where $C_1,c_1\in(0,\infty)$ depend only on $Q$. Hence the defining integral for $G_r$ is
absolutely convergent, and dominated convergence shows that $G_r$ is
continuous on $[0,\infty)^2$.

In the rest of the proof, we denote $H=L^2(\mathbb R)$ and $H_+=L^2(\mathbb R_+).$
    Let $ \mathsf A_-,\mathsf A_+:H_+\to H$ and $K_{\Ai}: H\to H$ be the operators from Lemma \ref{lem:operators Airy}, and let $R_r$ be multiplication on $H$ by $\mathbf 1_{(-r,\infty)}$.
    Then $X_r=R_r\mathsf A_-: H_+\to H$ is Hilbert--Schmidt, since it has kernel $(x,u)\mapsto
 \mathbf 1_{\{x>-r\}}\Ai(x+u)$ in $L^2(\mathbb R\times\mathbb R_+)$. Hence $A_r=R_rK_{\Ai}R_r= X_rX_r^*$ is positive and trace class. It is also a contraction, since for $h\in H$, $0\leq\langle h,A_rh\rangle_{H}
 =\langle R_rh,K_{\Ai}R_rh\rangle_{H}
 =\|K_{\Ai}R_rh\|_{H}^2
 \leq\|h\|_{H}^2$.  We next claim that 
\begin{equation}\label{eq:claim leakage}\int_0^\infty\!\!\int_0^\infty |G_r(u,v)|^2\,du\,dv 
 =\operatorname{Tr}_{ H}(A_r-A_r^2).\end{equation} 
Indeed, let $B_r=\mathsf A_-^*R_r\mathsf A_+$. For each $f,g\in C_c^\infty(\mathbb R_+)$, we have
\begin{equation}
    \label{eq:htrynyre}\begin{aligned}
 &\langle f,B_rg\rangle_{H_+}
  =\langle \mathsf A_-f,R_r\mathsf A_+g\rangle_H\\
 &=\int_0^\infty\int_0^\infty
   \overline{f(u)}
   \left(
     \int_{-r}^{\infty}\Ai(x+u)\Ai(x-v)\,dx
   \right)
   g(v)\,dv\,du =\int_0^\infty\int_0^\infty
   \overline{f(u)}\,G_r(u,v)\,g(v)\,dv\,du,
\end{aligned}\end{equation}
where the use of Fubini's theorem is justified by 
\[
 \int_{-r}^{\infty}\int_K\int_K
 |\Ai(x+u)\Ai(x-v)f(u)g(v)|\,dv\,du\,dx<\infty,
\]
from Lemma \ref{lem:Airy decay}. Here $K=\mathrm{supp}(f)\cup\mathrm{supp}(g)\subset\mathbb R_+$ is compact.
Since $R_r=R_r^*$, we have $B_r=X_r^*\mathsf A_+$.  As $X_r$ is
Hilbert--Schmidt and $\mathsf A_+$ is bounded, the Hilbert--Schmidt
ideal property shows that $B_r$ is Hilbert--Schmidt.  In view of \eqref{eq:htrynyre}, the Hilbert--Schmidt kernel of $B_r$ equals $G_r$ almost everywhere on
$\mathbb R_+^2$.
Hence 
\begin{equation}\label{eq:leakage proof1}
    \int_0^\infty\!\!\int_0^\infty |G_r(u,v)|^2\,du\,dv
 =\|B_r\|_{\mathrm{HS}}^2.
\end{equation}
On the other hand, define $S_r=X_r^*X_r=\mathsf A_-^*R_r\mathsf A_-$, which is positive and trace class. Then 
\begin{align*}
 B_rB_r^* =\mathsf A_-^*R_r\mathsf A_+\mathsf A_+^*R_r\mathsf A_- =\mathsf A_-^*R_r(I_{H}-K_{\Ai})R_r\mathsf A_- =S_r-\mathsf A_-^*R_r\mathsf A_-\mathsf A_-^*R_r\mathsf A_-
 =S_r-S_r^2.
\end{align*}
Moreover, 
\[
 \operatorname{Tr}_{H}(A_r)
 =\operatorname{Tr}_{H}(X_rX_r^*)
 =\operatorname{Tr}_{H_+}(X_r^*X_r)
 =\operatorname{Tr}_{H_+}(S_r),
\]
and, using $A_r^2=X_rS_rX_r^*$ and  cyclicity,
\[
 \operatorname{Tr}_{H}(A_r^2)
 =\operatorname{Tr}_{H_+}(S_rX_r^*X_r)
 =\operatorname{Tr}_{H_+}(S_r^2).
\]
Combining the above three displays, $ \operatorname{Tr}_{H}(A_r-A_r^2)=\|B_r\|_{\mathrm{HS}}^2$. Combining with \eqref{eq:leakage proof1}, we conclude  \eqref{eq:claim leakage}. 

Let $\mathcal X_{\Ai}$ denote the Airy  point process with correlation kernel $K_{\Ai}$, and $N_r=\#(\mathcal X_{\Ai}\cap (-r,\infty))$. Then 
\begin{equation}\label{eq:bhtrw}
 \operatorname{Var}(N_r)
 =\operatorname{Tr}_{H}(A_r-A_r^2).
\end{equation}
Indeed, 
\[
 \mathbb E[N_r]
 =\int_{(-r,\infty)}K_{\Ai}(x,x)\,dx
 =\operatorname{Tr}_{H}(A_r)<\infty,
\]
and $\rho_2(x,y)
 =K_{\Ai}(x,x)K_{\Ai}(y,y)-|K_{\Ai}(x,y)|^2$, hence  
\[
\mathbb E[N_r(N_r-1)] =\int_{(-r,\infty)^2}\rho_2(x,y)\,dx\,dy =(\operatorname{Tr}_{H} A_r)^2-\operatorname{Tr}_{H}(A_r^2),
\]
and we conclude \eqref{eq:bhtrw}. 
In view of \cite[(1.1) and the variance formula on p.~1303]{CharlierClaeys2020} with $x=-r$,  there exist absolute constants $r_2>1$ and
$C_2<\infty$ such that
\[
 \operatorname{Tr}_{H}(A_r-A_r^2)
 =\operatorname{Var}(N_r)
 \leq C_2\log r,
 \qquad r\geq r_2.
\]
It remains to treat $1\leq r\leq r_2$. Since $A_r$ is a positive contraction, $0\leq\operatorname{Tr}_{H}(A_r-A_r^2)
 \leq\operatorname{Tr}_{H}(A_r)$. Hence  
\begin{align*}
 \operatorname{Tr}_{H}(A_r)
  =\int_{-r}^{\infty}K_{\Ai}(x,x)\,dx 
  \leq\int_{-r_2}^{\infty}K_{\Ai}(x,x)\,dx 
  =\int_{-r_2}^{\infty}(z+r_2)\Ai(z)^2\,dz
 =C_3<\infty.
\end{align*}
We conclude the proof of the lemma with $C=\max(C_2,C_3/\log3)$.
\end{proof}

\subsection{An alternative formula for the Pfaffian kernels} \label{subsec:3.2}
In this subsection, we derive an alternative formula for the Pfaffian correlation kernel $\hski$ at equal times $s=t>0$ in terms of Airy-type integrals.

\begin{lemma}\label{lem:alternative kernels}
    For $t>0$ and $x\in\mathbb R$, define
\[
 h_t(x)=\int_0^\infty e^{-ts}\Ai(x+s)\,ds.
\] 
Recall the Pfaffian kernel $\hski$ from Definition \ref{def:kcrinfty}. For $t>0$, let $K^{t;\mathrm{hs};\infty}(x,y)=K^{\mathrm{hs};\infty}(t,x;t,y)$. Then
\begin{align}
 \label{eq:alt11}  K_{11}^{t;\mathrm{hs};\infty}(x,y)
 &=\frac14\int_0^\infty e^{-2ts}
 \left(h_t(x+s)h_t'(y+s)-h_t'(x+s)h_t(y+s)\right)\,ds,\\
\label{eq:alt12}   K_{12}^{t;\mathrm{hs};\infty}(x,y)
 &=K_{\Ai}(x,y)-\frac12h_t(x)\Ai(y),\\
\label{eq:alt22}  K_{22}^{t;\mathrm{hs};\infty}(x,y)
 &=\int_0^\infty e^{-2tv}
 \left(\Ai(x-v)\partial_v\Ai(y-v)
 -\partial_v\Ai(x-v)\Ai(y-v)\right)\,dv.
\end{align}
We also have an alternative formula \begin{equation}\label{eq:hski operator} K_{11}^{t;\mathrm{hs};\infty}(x,y) =\left(\mathsf A_-\mathsf J_t\mathsf A_-^*\right)(x,y),\end{equation} for $t>0$, where
$\mathsf A_-:L^2(\mathbb R_+)\to L^2(\mathbb R)$ is the operator from Lemma \ref{lem:operators Airy}, and
$\mathsf J_t$ is the Hilbert--Schmidt operator on
$L^2(\mathbb R_+)$ with kernel
\[  J_t(s,u)=\frac14e^{-t(s+u)}\operatorname{sgn}(u-s).\] The right side of \eqref{eq:hski operator} denotes the continuous
representative of the Hilbert--Schmidt integral kernel of
$\mathsf A_-\mathsf J_t\mathsf A_-^*$. 
\end{lemma}
\begin{proof}
    We take $s=t>0$ in Definition \ref{def:kcrinfty}. With $\Gamma_t=\mathcal C_{1+t}^{\pi/3}$, we have
    \begin{align}
 \label{eq:altproof11} K_{11}^{t;\mathrm{hs};\infty}(x,y)
 &=\frac1{(2\pi\im)^2}\int_{\Gamma_t}\!dz\int_{\Gamma_t}\!dw\,
   \frac{z-w }
        {4(z+w+2t)(z+t)(w+t)}e^{z^3/3+w^3/3-xz-yw},\\
\label{eq:altproof12} K_{12}^{t;\mathrm{hs};\infty}(x,y)
 &=\frac1{(2\pi\im)^2}\int_{\Gamma_t}\!dz\int_{\Gamma_t}\!dw\,
   \frac{z-w+2t}{2(z+t)(z+w)}e^{z^3/3+w^3/3-xz-yw},\\
\label{eq:altproof22} K_{22}^{t;\mathrm{hs};\infty}(x,y)
 &=I_{22}^{t;\mathrm{hs};\infty}(x,y)+R_{22}^{t;\mathrm{hs};\infty}(x,y),\\
 \label{eq:altproofI22}I_{22}^{t;\mathrm{hs};\infty}(x,y)
 &=\frac1{(2\pi\im)^2}\int_{\Gamma_t}\!dz\int_{\Gamma_t}\!dw\,
   \frac{z-w}{z+w-2t}e^{z^3/3+w^3/3-xz-yw},\\
 \label{eq:altproofR22}R_{22}^{t;\mathrm{hs};\infty}(x,y)
 &=\frac{y-x}{2^{5/2}\sqrt\pi\,t^{3/2}}
   e^{-(x-y)^2/(8t)-t(x+y)+2t^3/3}.
\end{align}
By \cite[(9.5.4)]{DLMF} and deforming the contour to $\Gamma_t$ in view of the exponential cubic decay of the integrand, 
\begin{equation}\label{eq:Airy contour}
 \Ai(a)=\frac{1}{2\pi\im}\int_{\Gamma_t}
 e^{z^3/3-az}\,dz\qquad\mbox{ for each }a\in\mathbb R.
\end{equation} By Fubini's theorem and differentiation under the contour integral,
\[
 h_t(a)=\frac1{2\pi\im}\int_{\Gamma_t}
 \frac{e^{z^3/3-az}}{z+t}\,dz,
 \qquad
 h_t'(a)=-\frac1{2\pi\im}\int_{\Gamma_t}
 \frac{z\,e^{z^3/3-az}}{z+t}\,dz.
\]
For $z,w\in\Gamma_t$, inserting
$(z+w+2t)^{-1}=\int_0^\infty e^{-r(z+w+2t)}\,dr$ into \eqref{eq:altproof11} gives
\[
 K_{11}^{t;\mathrm{hs};\infty}(x,y)
 =\frac14\int_0^\infty e^{-2tr}
 \left(h_t(x+r)h_t'(y+r)-h_t'(x+r)h_t(y+r)\right)\,dr,
\]
which is \eqref{eq:alt11}.  
For $a\in\{x,y\}$, set
 $f_a(r)=e^{-tr}\Ai(a+r) $
and
 $F_a(r)=\int_r^\infty f_a(v)\,dv
       =e^{-tr}h_t(a+r)$.
We have 
\begin{align*}
&\int_0^\infty\!\!\int_0^\infty
 e^{-t(s+u)}\operatorname{sgn}(u-s)
 \Ai(x+s)\Ai(y+u)\,ds\,du=\int_0^\infty\!\!\int_0^\infty
 f_x(s)f_y(u)\operatorname{sgn}(u-s)\,ds\,du \\&=\int_0^\infty
 \left(f_x(r)F_y(r)-F_x(r)f_y(r)\right)\,dr =\int_0^\infty
 \left(-F_x'(r)F_y(r)+F_x(r)F_y'(r)\right)\,dr\\
& =\int_0^\infty e^{-2tr}
 \left(h_t(x+r)h_t'(y+r)-h_t'(x+r)h_t(y+r)\right)\,dr  =4K_{11}^{t;\mathrm{hs};\infty}(x,y),
\end{align*}
where the use of Fubini's theorem is justified by Lemma \ref{lem:Airy decay}. By absolute convergence, $K_{11}^{t;\mathrm{hs};\infty}(x,y)$ is continuous in $x,y$ for fixed $t>0$. For each
$\phi,\psi\in C_c^\infty(\mathbb R)$, by  the
definition of $\mathsf A_-$  and the above display,
\[
 \left\langle
 \phi,\mathsf A_-\mathsf J_t\mathsf A_-^*\psi
 \right\rangle_{L^2(\mathbb R)}
 =
 \int_{\mathbb R}\int_{\mathbb R}
 \overline{\phi(x)}\,K_{11}^{t;\mathrm{hs};\infty}(x,y)\psi(y)\,dy\,dx.
\]  
Since $\mathsf J_t$ is Hilbert--Schmidt and $\mathsf A_-$ is bounded,
$\mathsf A_-\mathsf J_t\mathsf A_-^*$ is Hilbert--Schmidt. We conclude \eqref{eq:hski operator}. 

For the $12$ entry, use
\[
 \frac{z-w+2t}{2(z+t)(z+w)}
 =\frac1{z+w}-\frac1{2(z+t)}
\]
in \eqref{eq:altproof12}. We then use $(z+w)^{-1}=\int_0^\infty e^{-r(z+w)}\,dr$  and $(z+t)^{-1}=\int_0^\infty e^{-r(z+t)}\,dr$   for $z,w\in\Gamma_t$. We get
\begin{align*}
K_{12}^{t;\mathrm{hs};\infty}(x,y)
 &=\int_0^\infty\Ai(x+u)\Ai(y+u)\,du -\frac12
 \Ai(y)\int_0^\infty e^{-tu}\Ai(x+u)\,du  =K_{\Ai}(x,y)-\frac12h_t(x)\Ai(y),
\end{align*}
which is \eqref{eq:alt12}.

For the $22$ entry, inserting
$(z+w-2t)^{-1}=\int_0^\infty e^{-s(z+w-2t)}\,ds$ for $z,w\in\Gamma_t$ into \eqref{eq:altproofI22} gives 
\begin{equation}\label{eq:I22-sector}
I_{22}^{t;\mathrm{hs};\infty}(x,y)=\int_0^\infty e^{2ts}W_s(x,y)\,ds,\qquad\mbox{where } W_s(x,y)=\Ai(x+s)\Ai'(y+s)-\Ai'(x+s)\Ai(y+s).
\end{equation}
We claim that 
\[
 P_t(x,y)=\int_{\mathbb R}e^{2ts}\Ai(x+s)\Ai(y+s)\,ds
 =\frac{e^{-(x-y)^2/(8t)-t(x+y)+2t^3/3}}{\sqrt{8\pi t}}.
\]
Indeed, by deforming the Airy contour in \eqref{eq:Airy contour} to
$\Re z=t$ and parameterizing the latter by $z=t-\im p$, 
\begin{equation}\label{eq:whyrthsdf}
 \Ai(u)=\frac{1}{2\pi}\int_{\mathbb R}
 \exp\left(
 \frac{\im}{3}(p+\im t)^3+\im(p+\im t)u
 \right)\,dp.
\end{equation}
By Lemma \ref{lem:Airy decay},
$f_{x,t}(s)=e^{ts}\Ai(x+s)$ belongs to $L^2(\mathbb R)$.
Use the Fourier convention
\[
 \widehat f(p)=\int_{\mathbb R}e^{-\im ps}f(s)\,ds,
 \qquad
 f(s)=\frac1{2\pi}\int_{\mathbb R}e^{\im ps}\widehat f(p)\,dp.
\]
Set
$ g_{x,t}(p)=
 \exp (
 \frac{\im}{3}(p+\im t)^3+\im(p+\im t)x )$.
Since
 $|g_{x,t}(p)|=e^{-tp^2+t^3/3-tx}$,
we have $g_{x,t}\in L^1(\mathbb R)\cap L^2(\mathbb R)$.
Multiplying \eqref{eq:whyrthsdf} by $e^{ts}$ gives the
absolutely convergent inverse Fourier representation
\[
 f_{x,t}(s)=\frac1{2\pi}\int_{\mathbb R}
 e^{\im ps}g_{x,t}(p)\,dp.
\]
The $L^2$ Fourier-inversion theorem therefore yields
\[
 \widehat f_{x,t}(p)=g_{x,t}(p)
 =
 \exp\left(
 \frac{\im}{3}(p+\im t)^3+\im(p+\im t)x
 \right).
\]
Since $f_{y,t}$ is real-valued, Plancherel's theorem gives
\begin{align*}
&\int_{\mathbb R}e^{2ts}\Ai(x+s)\Ai(y+s)\,ds
 =\frac{1}{2\pi}\int_{\mathbb R}
 \widehat f_{x,t}(p)\,
 \overline{\widehat f_{y,t}(p)}\,dp\\
&=\frac{1}{2\pi}\int_{\mathbb R}
 \exp\left(
 \frac{\im}{3}\left((p+\im t)^3+(-p+\im t)^3\right)
 +\im(p+\im t)x+\im(-p+\im t)y
 \right)\,dp\\
&=\frac{e^{-t(x+y)+2t^3/3}}{2\pi}
 \int_{\mathbb R}e^{-2tp^2+\im p(x-y)}\,dp =\frac{e^{-(x-y)^2/(8t)-t(x+y)+2t^3/3}}
 {\sqrt{8\pi t}}.
\end{align*}
Moreover, Lemma \ref{lem:Airy decay} shows, locally uniformly in
$(x,y)$, that
\[
 s\longmapsto e^{2ts}\left(
 |\Ai'(x+s)\Ai(y+s)|
 +|\Ai(x+s)\Ai'(y+s)|
 \right)
\]
is integrable on $\mathbb R$.  Hence the first derivatives of $P_t$
may be taken under its defining integral.
Therefore,
\[
 \int_{\mathbb R}e^{2ts}W_s(x,y)\,ds
 =(\partial_y-\partial_x)P_t(x,y)
 =\frac{x-y}{2t}P_t(x,y)=-R_{22}^{t;\mathrm{hs};\infty}(x,y).
\]
Combining this  with \eqref{eq:I22-sector}, then substituting
$s=-v$ on $(-\infty,0)$, we conclude the proof of \eqref{eq:alt22}.
\end{proof}

\subsection{Uniform upper tail bound of the ensembles: Proof of Proposition \ref{prop:uniform upper tail}}
\label{subsec:uniform tail}

\begin{proof}[Proof of Proposition \ref{prop:uniform upper tail}]
    We begin by proving the first assertion.    By Proposition~\ref{Prop.HSAiryinfty},
\cite[(5.3) and (5.12)]{dimitrov2026half}, and monotone convergence,
\begin{equation}\label{eq:mdyrtsau}
\mathbb P\left(\hsai_1(t)\geq x\right)
\leq 
\int_x^{\infty}
K_{12}^{t;\mathrm{hs};\infty}(u,u)\,du.
\end{equation}
Recall from \eqref{eq:alt12} in Lemma \ref{lem:alternative kernels} that
\begin{equation*} 
    K_{12}^{t;\mathrm{hs};\infty}(u,u)=K_{\Ai}(u,u)-\frac12h_t(u)\Ai(u).\end{equation*}
For $u\ge0$, by $e^{-ts}\le1$ in $h_t$ and Lemma \ref{lem:Airy decay}, there exist absolute constants $C,c\in(0,\infty)$ such that $0\leq K_{12}^{t;\mathrm{hs};\infty}(u,u)  \leq Ce^{-cu^{3/2}}$ for each $t>0$.
Plugging in \eqref{eq:mdyrtsau}, we conclude the proof of the first assertion.    

We next prove the second assertion.       By Proposition~\ref{Prop.HSAiry},
\cite[(5.3) and (5.12)]{dimitrov2026half}, and monotone convergence,
\begin{equation*} 
\mathbb P\left(\hsa_1(t)\geq x\right)
\leq 
\int_x^{\infty}
K_{12}^{t;\mathrm{hs};\varpi}(u,u)\,du,
\end{equation*}
where $K^{t;\mathrm{hs};\varpi}(x,y)=K^{\mathrm{hs};\varpi}(t,x;t,y)$.
Therefore, it suffices to show that
\begin{equation}\label{eq:nsrythg}
 0\le K_{12}^{t;\mathrm{hs};\varpi}(u,u)
 \le C_{\varpi,L}e^{-6u}
 \qquad\mbox{ for } t\in[0,L],\, u\ge0.
\end{equation}
By Definition \ref{def:kcr},
\begin{equation}
K_{12}^{t;\mathrm{hs};\varpi}(u,u)=\frac1{(2\pi\im)^2}
 \int_{\mathcal C_{|\varpi|+9-t}^{\pi/3}}\!dz
 \int_{\mathcal C_{|\varpi|+3-t}^{2\pi/3}}\!dw\,
 \frac{(p+q)(p+\varpi)}{2p(p-q)(q+\varpi)}
 e^{z^3/3-w^3/3-u(z-w)},
 \label{eq:finite-density-contour}
\end{equation}
where $p=z+t$ and $q=w+t$. We parameterize the rays by radial variables $r_z,r_w\ge0$. Then
\[ 
 \operatorname{Re} (p)\ge|\varpi|+9,\qquad |p-q|\ge6,\qquad|q+\varpi|^2
\geq3(|\varpi|+\varpi+3)^2/4\geq 27/4,
\] 
 and hence there exists $C_1\in(0,\infty)$ depending only on $\varpi$ such that, on each pair of rays,
\[
\left|
\frac{(p+q)(p+\varpi)}
     {2p(p-q)(q+\varpi)}
\right|
\le C_1(1+r_z+r_w)^2.
\] 
For the exponent, we have 
$\operatorname{Re}(z-w)=6+(r_z+r_w)/2$, 
\begin{align*}
 \operatorname{Re}(z^3/3) =-r_z^3/3-(a/2)r_z^2+(a^2/2)r_z+a^3/3,\qquad
 \operatorname{Re}(-w^3/3) =-r_w^3/3+(b/2)r_w^2+(b^2/2)r_w-b^3/3,
\end{align*}
where $a=|\varpi|+9-t$ and $b=|\varpi|+3-t$.  Therefore, uniformly for $t\in[0,L]$, 
\[
\operatorname{Re}(z^3/3-w^3/3-u(z-w))\leq - (r_z^3+r_w^3)/3
 +C_2(1+r_z^2+r_w^2)-6u- (r_z+r_w)u/2\leq   C_3- (r_z^3+r_w^3)/6-6u,
\]
where $C_2,C_3\in(0,\infty)$ depend only on $\varpi,L$. We conclude the proof of \eqref{eq:nsrythg}.
\end{proof}

\subsection{Controlling the variance: Proof of Proposition \ref{prop:bounding variance}} \label{subsec:variance control}
For $t>0$, consider the point process
\[
 M^{t;\mathrm{hs};\infty}(\omega,A)
  =\sum_{k=1}^{\infty}\mathbf 1\{\hsai_k(t,\omega)\in A\}.
\]  
By Proposition \ref{Prop.HSAiryinfty} and \cite[Lemma 5.13]{dimitrov2026half}, this is a Pfaffian point process on $\mathbb R$ with reference measure $\mathrm{Leb}$ and with correlation kernel $K^{t;\mathrm{hs};\infty}(x,y)=K^{\mathrm{hs};\infty}(t,x;t,y)$.

\begin{lemma}\label{lem:ht bounds} Let $\mathsf A_-$ and $K_{\Ai}$ be as in  Lemma \ref{lem:operators Airy}. For $t>0$ and $x\in\mathbb R$, let 
    \[
    h_t(x)=\int_0^\infty e^{-ts}\Ai(x+s)\,ds.
    \] 
    Then there exists an absolute constant $C\in(0,\infty)$ such that for each $t>0$ and $r\geq1$,
\begin{align}
\label{eq:first ht bounds}& h_t=\mathsf A_-(e^{-t\,\cdot}),\qquad K_{\Ai}h_t=h_t,\qquad
 \|h_t\|_2^2= 1/(2t), \qquad \sup_{x\in\mathbb R}|h_t(x)|\in[0,C],\\
\label{eq:second line ht}    &  h_t'=t h_t-\Ai,\qquad \int_{-r}^\infty h_t(x)\Ai(x)\,dx
 =t\int_{-r}^\infty h_t(x)^2\,dx+\frac12h_t(-r)^2\in[0,C].
\end{align} 
\end{lemma}
\begin{proof}
    In this proof, we denote $f_t(s)=e^{-ts}$ for $s\ge0$. Choose $\chi_n\in C_c^\infty(\mathbb R_+)$ such that
$0\le\chi_n\le1$ and $\chi_n(s)\to1$ for every $s>0$. By dominated convergence, $\chi_nf_t\to f_t$ in $L^2(\mathbb R_+)$. Hence,
by Lemma~\ref{lem:operators Airy}, $\mathsf A_-(\chi_nf_t)\to\mathsf A_-f_t$ in $L^2(\mathbb R)$. For each fixed $x\in\mathbb R$, by the  dominated convergence theorem with dominating function $s\mapsto e^{-ts}|\Ai(x+s)|$,  we have
$\mathsf A_-(\chi_nf_t)(x)
=\int_0^\infty\chi_n(s)e^{-ts}\Ai(x+s)\,ds
\to h_t(x)$.
 Passing to an almost everywhere convergent subsequence of $\mathsf A_-(\chi_nf_t)$, we have $h_t=\mathsf A_-f_t$ in $L^2(\mathbb R)$.  Therefore, by Lemma~\ref{lem:operators Airy},
 \[
 K_{\Ai}h_t
=\mathsf A_-\mathsf A_-^*\mathsf A_-f_t
=\mathsf A_-f_t=h_t,\qquad\qquad \|h_t\|_{L^2(\mathbb R)}^2
=\|f_t\|_{L^2(\mathbb R_+)}^2
=\int_0^\infty e^{-2ts}\,ds
= 1/(2t).
 \] Define $F(x)=\int_x^\infty\Ai(u)\,du$ for $x\in\mathbb R$. By Lemma~\ref{lem:Airy decay}, this integral is well-defined, and $F$ is bounded on $[0,\infty)$. To control the negative
axis, put $\theta(s)= 2/3\cdot s^{3/2}+ \pi/4$.
 For $s\geq1$, recall from Lemma~\ref{lem:Airy decay} that
 \[
\Ai(-s)=\pi^{-1/2}s^{-1/4}\sin\theta(s)+\rho(s),
\qquad |\rho(s)|\le C_1s^{-7/4},
\]
where $C_1$ is an absolute constant. Note that $\theta'(s)=s^{1/2}$. Using integration by parts, we have, for $a\geq1$,
\[
\int_1^a s^{-1/4}\sin\theta(s)\,ds
=
-\left[s^{-3/4}\cos\theta(s)\right]_{1}^{a}
-\frac34\int_1^a s^{-7/4}\cos\theta(s)\,ds.
\] The right side above and $\int_1^a\rho(s)\,ds$ are bounded uniformly in $a\ge1$. Since
$F(-a)=F(-1)+\int_1^a\Ai(-s)\,ds$, 
and $F$ is continuous on $[-1,0]$, we conclude that
$C_2=\sup_{x\in\mathbb R}|F(x)|<\infty$. Note that $F'=-\Ai$. Using integration by parts first on $[0,R]$ and then
letting $R\to\infty$, we get
 \[
h_t(x)
=F(x)-t\int_0^\infty e^{-ts}F(x+s)\,ds.
\]
Indeed, the boundary term $e^{-tR}F(x+R)\to0$, and the remaining integral is absolutely convergent.
 Thus
\[
|h_t(x)|
\le C_2+tC_2\int_0^\infty e^{-ts}\,ds
=2C_2
\qquad\mbox{ for each } t>0,\, x\in\mathbb R.
\]

In the remainder of the proof, we show \eqref{eq:second line ht}. By a change of variables, 
$h_t(x)=e^{tx}\int_x^\infty e^{-tu}\Ai(u)\,du$. 
Hence $h_t\in C^1(\mathbb R)$ and $h_t'(x)=t h_t(x)-\Ai(x)$. By Lemma~\ref{lem:Airy decay},  for $x\ge0$, $|h_t(x)|
\le\int_x^\infty|\Ai(u)|\,du\to0$ as $x\to\infty$. Again by Lemma~\ref{lem:Airy decay},  $\Ai\in L^2((-r,\infty))$. Since $h_t\in L^2(\mathbb R)$, Cauchy--Schwarz shows that $h_t\Ai\in L^1((-r,\infty))$. Using
$\Ai=t h_t-h_t'$ and integrating on $[-r,R]$, we obtain
\[
\int_{-r}^{R}h_t(x)\Ai(x)\,dx
=
t\int_{-r}^{R}h_t(x)^2\,dx
+\frac12h_t(-r)^2-\frac12h_t(R)^2.
\] Taking $R\to\infty$,
\[
\int_{-r}^{\infty}h_t(x)\Ai(x)\,dx
=
t\int_{-r}^{\infty}h_t(x)^2\,dx+\frac12h_t(-r)^2.
\]
Note that we have
\[0\leq
t\int_{-r}^{\infty}h_t(x)^2\,dx+\frac12h_t(-r)^2
\le
t\|h_t\|_{L^2(\mathbb R)}^2+\frac12\|h_t\|_\infty^2
\le\frac12+2C_2^2.
\] Taking $C=\max(2C_2,1/2+2C_2^2)$, we conclude the proof.
\end{proof}

\begin{lemma}\label{lem:exact identity variance}
    For \(t>0\) and \(r\geq1\), let $N_{t,r}$ be as in Proposition \ref{prop:bounding variance}. For $u,v\geq0$, let
    \[
 G_r(u,v) =\int_{-r}^\infty\Ai(x+u)\Ai(x-v)\,dx .
\] 
Let $h_t(x)$ and $K_{\Ai}$ be as in Lemma \ref{lem:ht bounds}. Let $R_r$ be the multiplication on $L^2(\mathbb R)$ by $\mathbf 1_{(-r,\infty)}$.  We use \(R_r\Ai\) to denote the function
$x\mapsto \mathbf 1_{\{x>-r\}}\Ai(x)$,
which belongs to \(L^2(\mathbb R)\).
 Then \(N_{t,r}<\infty\) almost surely and \(N_{t,r}\in L^2\),
all quantities below are finite, and
\begin{equation}\label{eq:etyjnaed}
    \begin{split}
         \operatorname{Var}(N_{t,r})=& \int_0^\infty\!\!\int_0^\infty
 \left(e^{-2t(u+v)}+1\right)G_r(u,v)^2\,du\,dv\\
 &+\frac{\langle R_r\Ai,R_rh_t\rangle}{2}-\frac{\langle R_r\Ai,R_rh_t\rangle^2}{2}-\langle R_r\Ai,K_{\Ai}(I-R_r)h_t\rangle\\
 & -\int_0^\infty\!\!\int_0^\infty
 e^{-2t(u+v)}h_t(-r+u)\Ai(-r-v)G_r(u,v)\,du\,dv.
    \end{split}
\end{equation}
\end{lemma}
\begin{proof}
    Throughout the proof, we denote  $q=\langle R_r\Ai,R_rh_t\rangle$, $A_r=R_rK_{\Ai}R_r$, and the kernels \(K_{ij}^{t;\mathrm{hs};\infty}\) by \(K_{ij}\). By \eqref{eq:claim leakage} in the proof of Lemma~\ref{lem:Airy leakage}, and Lemma~\ref{lem:ht bounds}, 
    \begin{align}
  \label{eq:variance-use-leakage}&A_r\mbox{ is trace class, }\qquad       \operatorname{Tr}(A_r-A_r^2)
 =\int_{\mathbb R_+^2}G_r(u,v)^2\,du\,dv<\infty,\\
\label{eq:variance-h-derivative} &h_t\in L^2(\mathbb R),\qquad K_{\Ai}h_t=h_t,
 \qquad h_t'=th_t-\Ai.
    \end{align}
    
    For clarity, we split the rest of the proof into two steps.\\

\noindent\textbf{Step 1. } In this step, we show that \(N_{t,r}<\infty\) almost surely and \(N_{t,r}\in L^2\).
We begin by proving the following estimates: there exist \(C_1,c_1\in(0,\infty)\), depending only on \(r,t\), such that, for each
\(u,v\geq0\), 
\begin{align}
 \|R_r\Ai(\,\cdot-v)\|_2&\leq C_1(1+v)^{1/4},\qquad 
 \|R_r\Ai'(\,\cdot-v)\|_2\leq C_1(1+v)^{3/4},
 \label{eq:variance-shift-minus}\\
 \|R_rf(\,\cdot+u)\|_2&\leq C_1 e^{-c_1(u-r)_+^{3/2}}\qquad\mbox{for }
 f\in\{\Ai,\Ai',h_t,h_t'\}.
 \label{eq:variance-shift-plus}
\end{align}
Indeed, by a change of variables
and Lemma~\ref{lem:Airy decay}, for \(v\geq0\),
\[
 \|R_r\Ai(\,\cdot-v)\|_2^2
 =\int_{-r-v}^{\infty}\Ai(y)^2\,dy
 \leq C_2(1+v)^{1/2},\qquad 
 \|R_r\Ai'(\,\cdot-v)\|_2^2
 =\int_{-r-v}^{\infty}\Ai'(y)^2\,dy
 \leq C_2(1+v)^{3/2}.
\]
By Lemma~\ref{lem:Airy decay}, for \(y\geq0\), we have
$ |h_t(y)|
 \leq C_3\int_0^\infty e^{-c_3(y+s)^{3/2}}\,ds
 \leq C_4 e^{-c_4 y^{3/2}}$,
and hence, by \eqref{eq:variance-h-derivative},
$ |h_t'(y)|
 \leq t|h_t(y)|+|\Ai(y)|
 \leq C_5 e^{-c_5 y^{3/2}}$.
Combining with Lemma \ref{lem:Airy decay}, we conclude  \eqref{eq:variance-shift-minus} and \eqref{eq:variance-shift-plus}. 

We next show that the operators  \(R_rK_{ij}R_r\) are Hilbert--Schmidt for $i,j\in\{1,2\}$. Indeed, by \eqref{eq:hski operator}, \(R_rK_{11}R_r\) is Hilbert--Schmidt. Put $\phi_v=R_r\Ai(\,\cdot-v)$ and $\psi_v=R_r\partial_v\Ai(\,\cdot-v)$. Then by  \eqref{eq:alt22}, we  have
\[
 R_rK_{22}R_r
 =\int_0^\infty e^{-2tv}
   (\phi_v\otimes\psi_v-\psi_v\otimes\phi_v)\,dv,
\]
where \((f\otimes g)(x,y)=f(x)g(y)\).
This integral converges absolutely in Hilbert--Schmidt norm because
\[
 \left\|e^{-2tv}
   (\phi_v\otimes\psi_v-\psi_v\otimes\phi_v)\right\|_{\mathrm{HS}}
 \leq 2e^{-2tv}\|\phi_v\|_2\|\psi_v\|_2
 \leq Ce^{-2tv}(1+v),
\]
which is integrable since \(t>0\).  Furthermore, \eqref{eq:alt12} gives
\begin{equation}\label{eq:srfg} R_rK_{12}R_r=A_r- 1/2\cdot (R_rh_t)\otimes(R_r\Ai),\end{equation}
so \(R_rK_{12}R_r\) is trace class and therefore Hilbert--Schmidt. 
 
By Proposition~\ref{Prop.HSAiryinfty}, the first two correlation functions of the point process at time \(t\) are
\begin{equation}\label{eq:nsrnbrmtue}
 \rho_1(x)=K_{12}(x,x),\qquad
 \rho_2(x,y)=\rho_1(x)\rho_1(y)
 -K_{12}(x,y)K_{12}(y,x)-K_{11}(x,y)K_{22}(x,y).
\end{equation}
We have \(\rho_1\in L^1((-r,\infty))\) in view of \eqref{eq:srfg}, since
\[
\int_{(-r,\infty)}K_{\Ai}(x,x)\,dx=\operatorname{Tr}(A_r)<\infty,\qquad  \int_{(-r,\infty)}|h_t(x)\Ai(x)|\,dx
 \leq\|R_rh_t\|_2\|R_r\Ai\|_2<\infty.
\]
By the Cauchy--Schwarz inequality,
\[\begin{aligned}
& \int_{(-r,\infty)^2}
 \left|K_{12}(x,y)K_{12}(y,x)\right|\,dx\,dy
 \leq\|R_rK_{12}R_r\|_{\mathrm{HS}}^2,\\
& \int_{(-r,\infty)^2}
 \left|K_{11}(x,y)K_{22}(x,y)\right|\,dx\,dy
 \leq
 \|R_rK_{11}R_r\|_{\mathrm{HS}}\|R_rK_{22}R_r\|_{\mathrm{HS}}.
\end{aligned}\]
 Consequently,
\[
 \|\rho_2\|_{L^1((-r,\infty)^2)}
 \leq \|\rho_1\|_{L^1((-r,\infty))}^2
 +\|R_rK_{12}R_r\|_{\mathrm{HS}}^2
 +\|R_rK_{11}R_r\|_{\mathrm{HS}}\|R_rK_{22}R_r\|_{\mathrm{HS}}
 <\infty.
\]

For \(m\in\mathbb N\), let
\(N_m=M^{t;\mathrm{hs};\infty}((-r,m))\). Then 
\[
 \mathbb E[N_m]=\int_{(-r,m)}\rho_1(x)\,dx,\qquad
 \mathbb E[N_m(N_m-1)]
 =\int_{(-r,m)^2}\rho_2(x,y)\,dx\,dy.
\]
Since \(N_m\uparrow N_{t,r}\), monotone convergence and the preceding integrability bounds imply 
 $\mathbb E[N_{t,r}]<\infty$ 
and
 $\mathbb E[N_{t,r}(N_{t,r}-1)]<\infty$. 
Thus \(N_{t,r}<\infty\) almost surely and \(N_{t,r}\in L^2\).\\

\noindent\textbf{Step 2. } In this step, we prove the exact identity \eqref{eq:etyjnaed}. 
We first observe that every term therein is finite. Indeed, the first line being finite follows from Lemma \ref{lem:Airy leakage}; the terms in the second line being finite follows from $|q|\leq\|R_r\Ai\|_2\|R_rh_t\|_2$ and $|\langle R_r\Ai,K_{\Ai}(I-R_r)h_t\rangle| \leq\|R_r\Ai\|_2\|(I-R_r)h_t\|_2$, since $K_{\Ai}$ is an orthogonal projection; and the third line being finite  follows from
Cauchy--Schwarz, \eqref{eq:variance-use-leakage}, \(h_t\in L^2(\mathbb R)\), and Lemma~\ref{lem:Airy decay}.

Substituting \eqref{eq:nsrnbrmtue} into
$ \operatorname{Var}(N_{t,r})
 =\mathbb E[N_{t,r}]
 +\mathbb E[N_{t,r}(N_{t,r}-1)]
 -\mathbb E[N_{t,r}]^2$ gives
\begin{equation}\label{eq:variance-pf-short}
 \operatorname{Var}(N_{t,r})
 =\int_{(-r,\infty)}K_{12}(x,x)\,dx
 -\int_{(-r,\infty)^2}\!
 \left(K_{12}(x,y)K_{12}(y,x)
       +K_{11}(x,y)K_{22}(x,y)\right)dx\,dy.
\end{equation}
Expanding \eqref{eq:alt12} in Lemma \ref{lem:alternative kernels}, using \(K_{\Ai}h_t=h_t\) and
\eqref{eq:variance-use-leakage}, we obtain
\begin{equation}\label{eq:variance-12-short}
    \begin{split}
         & \int_{(-r,\infty)}K_{12}(x,x)\,dx
 -\int_{(-r,\infty)^2}K_{12}(x,y)K_{12}(y,x)\,dx\,dy   =\operatorname{Tr}(A_r-A_r^2)-\frac q2
 +\langle R_r\Ai,K_{\Ai}R_rh_t\rangle-\frac{q^2}{4} \\
 & =\int_{\mathbb R_+^2}G_r(u,v)^2\,du\,dv
 +\frac q2-\frac{q^2}{4}
 -\langle R_r\Ai,K_{\Ai}(I-R_r)h_t\rangle.
    \end{split}
\end{equation}
Indeed, 
$ \langle R_r\Ai,K_{\Ai}R_rh_t\rangle
 =\langle R_r\Ai,K_{\Ai}h_t\rangle
  -\langle R_r\Ai,K_{\Ai}(I-R_r)h_t\rangle
 =q-\langle R_r\Ai,K_{\Ai}(I-R_r)h_t\rangle$.
 
We denote
\[
 J=\int_{(-r,\infty)^2}K_{11}(x,y)K_{22}(x,y)\,dx\,dy.
\]
In view of \eqref{eq:variance-12-short}, to prove \eqref{eq:etyjnaed},
it suffices to show  
\begin{align}
 J= &-\int_{\mathbb R_+^2}e^{-2t(u+v)}
 G_r(u,v)^2\,du\,dv+\frac{q^2}{4} +\int_{\mathbb R_+^2}e^{-2t(u+v)}
 h_t(-r+u)\Ai(-r-v)G_r(u,v)\,du\,dv.
 \label{eq:variance-J-short}
\end{align}
Define
\[
 H_r(u,v)=\int_{(-r,\infty)} h_t(x+u)\Ai(x-v)\,dx,
 \qquad
 L_r(u,v)=h_t(-r+u)\Ai(-r-v).
\]
We begin by proving that
\begin{equation}\label{eq:heyrtarweg}
     J=\frac12\int_{\mathbb R_+^2}e^{-2t(u+v)}
 \left(
 G_r(u,v)\partial_vH_r(u,v)
 -\partial_vG_r(u,v)H_r(u,v)
 \right)\,du\,dv.
\end{equation}
Indeed, the bounds \eqref{eq:variance-shift-minus} and \eqref{eq:variance-shift-plus} justify the use of all differentiation under integrals, integration by parts, and Fubini's theorem below.  
Using \(h_t'=th_t-\Ai\), we have $\partial_uH_r =tH_r -G_r $. 
Integration by parts in \(x\) gives
\begin{equation}\label{eq:wrytABS}
\begin{aligned}
 \partial_vH_r(u,v) &=-\int_{(-r,\infty)} h_t(x+u)\Ai'(x-v)\,dx =L_r(u,v)+\int_{(-r,\infty)} h_t'(x+u)\Ai(x-v)\,dx \\&=L_r(u,v)+tH_r(u,v)-G_r(u,v).
\end{aligned}
\end{equation}
Moreover, by \eqref{eq:variance-h-derivative} and
\eqref{eq:alt11}, 
\[
 K_{11}(x,y)=\frac14\int_0^\infty e^{-2tu}
 \left(
 \Ai(x+u)h_t(y+u)-h_t(x+u)\Ai(y+u)
 \right)\,du.
\]
Expanding this formula together with \eqref{eq:alt22} for $K_{22}$, using Fubini's theorem, we conclude \eqref{eq:heyrtarweg}. 
By \eqref{eq:wrytABS},  
\[
\begin{aligned}
 e^{-2t(u+v)}
 \left(G_r\partial_vH_r-(\partial_vG_r)H_r\right)
 ={}&-\partial_v\left(e^{-2t(u+v)}G_rH_r\right) -2e^{-2t(u+v)}G_r^2
 +2e^{-2t(u+v)}L_rG_r.
\end{aligned}
\]
By \eqref{eq:variance-shift-minus} and \eqref{eq:variance-shift-plus}, $e^{-2t(u+v)}G_r(u,v)H_r(u,v)\to0$ as  $v\to\infty$.  In view of \eqref{eq:heyrtarweg}, we have
\begin{equation*} 
 J={}-\int_{\mathbb R_+^2}e^{-2t(u+v)}G_r(u,v)^2\,du\,dv +\int_{\mathbb R_+^2}e^{-2t(u+v)}
 L_r(u,v)G_r(u,v)\,du\,dv +\frac12\int_0^\infty
 e^{-2tu}G_r(u,0)H_r(u,0)\,du.
\end{equation*}
Let \(Q(u)=H_r(u,0)\).  Then $Q'(u)=tQ(u)-G_r(u,0)$ and $Q(0)=q$.
Since \(e^{-2tu}Q(u)^2\to0\) as \(u\to\infty\),
\[
\begin{aligned}
 \int_0^\infty e^{-2tu}G_r(u,0)H_r(u,0)\,du
 &=-\frac12\int_0^\infty
 \frac{d}{du}\left(e^{-2tu}Q(u)^2\right)\,du =\frac{q^2}{2}.
\end{aligned}
\]
Since \(L_r(u,v)=h_t(-r+u)\Ai(-r-v)\), substitution into the last formula for $J$ proves \eqref{eq:variance-J-short}  and hence
\eqref{eq:etyjnaed}.
\end{proof}

\begin{lemma}\label{lem:technical bounds variance}
We use the same notation as in Lemma~\ref{lem:exact identity variance}.
Then there exists an absolute constant \(C\in(0,\infty)\) such that,
for every \(t>0\) and \(r\geq1\),
\begin{align}
 \left|\langle R_r\Ai,K_{\Ai}(I-R_r)h_t\rangle\right|
 &\leq C,
 \label{eq:technical-variance-first}\\
 \left|\int_0^\infty\!\!\int_0^\infty
 e^{-2t(u+v)}h_t(-r+u)\Ai(-r-v)G_r(u,v)\,du\,dv\right|
 &\leq C\log(2+r).
 \label{eq:technical-variance-second}
\end{align}
\end{lemma}

\begin{proof}
Throughout this proof, we use $C\in(0,\infty)$ to denote an absolute constant which may change from line to line. 
We also denote
\[\begin{aligned}
  a(u)&=\Ai(-r+u),\qquad b(v)=\Ai(-r-v),\qquad
 \eta(u)=h_t(-r+u),\\  p&=u+v,\qquad\qquad
 \Omega=\{(u,v)\in\mathbb R_+^2:p\geq1\}.\end{aligned}
\]
By Lemmas~\ref{lem:Airy decay} and~\ref{lem:ht bounds},  
\begin{gather}
 |a(0)|\leq Cr^{-1/4},\qquad
 b(0)^2\leq Cr^{-1/2},\qquad
 b(v)^2\leq C(r+v)^{-1/2},\qquad
 \|\eta\|_\infty\leq C,
 \label{eq:technical-basic-pointwise}\\
 \eta'=t\eta-a,\qquad
 \int_0^\infty a(u)^2\,du\leq Cr^{1/2},\qquad
 \int_0^\infty\frac{|a(u)|}{1+u}\,du\leq Cr^{1/4}.
 \label{eq:technical-basic-integral}
\end{gather}
Here the last estimate follows from Cauchy--Schwarz and the preceding
one. Therefore, we have
\begin{align}
 \int_{(1-u)_+}^\infty\frac{b(v)^2}{u+v}\,dv
 &\leq
 Cr^{-1/2}\int_{(1-u)_+}^r\frac{dv}{u+v}
 +C\int_r^\infty v^{-3/2}\,dv
 \leq Cr^{-1/2}\log(2+r),
 \label{eq:technical-one-dimensional-first}\\
 \int_{(1-u)_+}^\infty\frac{b(v)^2}{(u+v)^2}\,dv
 &\leq
 Cr^{-1/2}\int_{(1-u)_+}^\infty\frac{dv}{(u+v)^2}
 \leq \frac{Cr^{-1/2}}{1+u},
 \label{eq:technical-one-dimensional-second}
\end{align}
and hence,
\begin{align}
 \iint_\Omega\frac{a(u)^2b(v)^2}{p}\,du\,dv
 &\leq C\log(2+r),
 \label{eq:technical-bulk-log}\\
 \iint_\Omega\frac{|a(u)|b(v)^2}{p^2}\,du\,dv
 &\leq C.
 \label{eq:technical-bulk-uniform}
\end{align}
The estimates \eqref{eq:technical-basic-pointwise} and \eqref{eq:technical-basic-integral} also imply
\begin{equation}
 \begin{split}
 &|a(0)|\int_1^\infty\frac{b(v)^2}{v}\,dv
 +b(0)^2\int_1^\infty\frac{|a(u)|}{u}\,du 
 +\int_0^1|a(u)|b(1-u)^2\,du
 \leq C.
 \end{split}
 \label{eq:technical-boundary-uniform}
\end{equation}

Define $\zeta_1(u,v)=e^{-tu}$ and $\zeta_2(u,v)=e^{-2tp}\eta(u)$.
We next prove the claim that 
\begin{equation}\label{eq:technical-master-short}
    \begin{split}
        & \iint_\Omega\zeta_i bG_r\,du\,dv=a(0)\int_1^\infty
       \frac{\zeta_i(0,v)b(v)^2}{v}\,dv
 +\frac{b(0)^2}{2}\int_1^\infty
       \frac{\zeta_i(u,0)a(u)}{u}\,du \\
       &+\frac32\int_0^1
       \zeta_i(u,1-u)a(u)b(1-u)^2\,du +\iint_\Omega\frac{a(u)b(v)^2}{p}
       \left(\partial_u\zeta_i+
             \frac12\partial_v\zeta_i\right)\,du\,dv-\frac32\iint_\Omega
       \frac{\zeta_i a(u)b(v)^2}{p^2}\,du\,dv,
    \end{split}
\end{equation}
where we note that for \(i=1\), the left side denotes the iterated improper integral
determined by this order of limits. Indeed, using  \eqref{eq:Darboux KAi}, for
\(p>0\),
\[
 G_r(u,v)
 =K_{\Ai}(-r+u,-r-v)
 =-\left(a(u)b'(v)+a'(u)b(v)\right)/p,
\]
 and hence
\begin{equation}
 b(v)G_r(u,v)
 =-\frac{a(u)}{2p}\,\partial_v\!\left(b(v)^2\right)
  -\frac{b(v)^2}{p}\,\partial_u a(u).
 \label{eq:technical-darboux-short}
\end{equation}
Using \(\eta'=t\eta-a\), we have
\begin{align*}
 \partial_u\zeta_1+\frac12\partial_v\zeta_1
  =-te^{-tu},\qquad
 \partial_u\zeta_2+\frac12\partial_v\zeta_2
  =-e^{-2tp}\left(2t\eta(u)+a(u)\right).
\end{align*}
Since \(\Ai\) is bounded and \(te^{-2tp}\leq C/p\) on \(\Omega\),
\eqref{eq:technical-one-dimensional-first},
\eqref{eq:technical-bulk-log}, and \eqref{eq:technical-bulk-uniform}
give
\begin{align}
 \iint_\Omega\frac{|a|b^2}{p}
 \left|\partial_u\zeta_1+\frac12\partial_v\zeta_1\right|\,du\,dv
 &\leq
 Cr^{-1/2}\log(2+r)
 \int_0^\infty te^{-tu}|a(u)|\,du
 \leq C,
 \label{eq:technical-zeta-one-derivative}\\
 \iint_\Omega\frac{|a|b^2}{p}
 \left|\partial_u\zeta_2+\frac12\partial_v\zeta_2\right|\,du\,dv
 &\leq
 C\iint_\Omega\frac{|a|b^2}{p^2}\,du\,dv
 +C\iint_\Omega\frac{a^2b^2}{p}\,du\,dv  \leq C\log(2+r).
 \label{eq:technical-zeta-two-derivative}
\end{align} Integrating \eqref{eq:technical-darboux-short} by parts on
\(\Omega\cap([0,M]\times[0,V])\), then sending first \(M\to\infty\)
and then \(V\to\infty\), we conclude \eqref{eq:technical-master-short}. Indeed, for fixed \(V\), the boundary at
\(u=M\) vanishes by the decay of the Airy function along the positive axis. After
\(M\to\infty\), the boundary at \(v=V\) is bounded, for \(i=1\), by
 $Cb(V)^2/V\cdot
 \int_0^\infty e^{-tu}|a(u)|\,du\to0$.
For \(i=2\), it is bounded by
 $ Ce^{-2tV}b(V)^2/V\cdot
 \int_0^\infty e^{-2tu}|a(u)|\,du\to0$,
where we used \(\|\eta\|_\infty\leq C\).
The bulk terms are absolutely convergent by
\eqref{eq:technical-bulk-uniform},
\eqref{eq:technical-zeta-one-derivative},  
\eqref{eq:technical-zeta-two-derivative}, while the boundary terms
are absolutely convergent by
\eqref{eq:technical-boundary-uniform}.  

We next prove \eqref{eq:technical-variance-second}. 
Since \(\sup_\Omega|\zeta_i|\leq C\),
\eqref{eq:technical-boundary-uniform}, \eqref{eq:technical-master-short}, \eqref{eq:technical-bulk-uniform}, \eqref{eq:technical-zeta-one-derivative}, and \eqref{eq:technical-zeta-two-derivative} imply
\begin{equation}
 \left|\iint_\Omega\zeta_1bG_r\,du\,dv\right|\leq C,
 \qquad
 \left|\iint_\Omega\zeta_2bG_r\,du\,dv\right|
 \leq C\log(2+r).
 \label{eq:technical-two-weights-omega}
\end{equation}
On \(\{p<1\}\), Lemma~\ref{lem:Airy leakage} and
\eqref{eq:technical-basic-pointwise} give, for \(i=1,2\),
\[
 \iint_{\{p<1\}}|\zeta_i bG_r|\,du\,dv
 \leq
 Cr^{-1/4}\|G_r\|_{L^2(\mathbb R_+^2)}
 \leq Cr^{-1/4}\sqrt{\log(2+r)}
 \leq C.
\]
It follows that the iterated improper integral
\begin{equation}
 J_{t,r}=
 \lim_{V\to\infty}\int_0^V b(v)
 \left(\int_0^\infty e^{-tu}G_r(u,v)\,du\right)dv
 \label{eq:technical-J-short}
\end{equation}
exists and satisfies \(|J_{t,r}|\leq C\). For finite \(V\), the interchange of the \(u\)- and \(v\)-integrals is justified by
Cauchy--Schwarz, since \(G_r\in L^2(\mathbb R_+^2)\), \(e^{-t\,\cdot}\in L^2(\mathbb R_+)\), and \(b\in L^2((0,V))\).

The integral corresponding to \(\zeta_2\) is absolutely convergent. Indeed, 
\[
\begin{aligned}
 &\iint_{\mathbb R_+^2}
 e^{-2tp}|\eta(u)b(v)G_r(u,v)|\,du\,dv \leq
 \|G_r\|_2
 \left(\int_0^\infty e^{-4tu}\eta(u)^2\,du\right)^{1/2}
 \left(\int_0^\infty e^{-4tv}b(v)^2\,dv\right)^{1/2}
 <\infty.
\end{aligned}
\]
Therefore, by \eqref{eq:technical-two-weights-omega} and the estimate on \(\{p<1\}\), we conclude the proof of
\eqref{eq:technical-variance-second}.

In the remainder of the proof, we show \eqref{eq:technical-variance-first}. Define
\[
 \mathcal H(v)=\int_0^\infty e^{-tu}G_r(u,v)\,du,
 \qquad
 \alpha=h_t(-r),
 \qquad
q=\langle R_r\Ai,R_rh_t\rangle.
\]
For each fixed \(v\geq0\), Fubini's theorem applies because
$ \int_0^\infty e^{-tu}
 \|R_r\Ai(\,\cdot+u)\|_2
 \|R_r\Ai(\,\cdot-v)\|_2\,du<\infty$
by Lemma~\ref{lem:Airy decay}, and hence
\begin{equation}
 \mathcal H(v)
 =\int_{-r}^\infty h_t(x)\Ai(x-v)\,dx.
 \label{eq:technical-H-short}
\end{equation}
Lemma~\ref{lem:operators Airy} and \eqref{eq:technical-H-short} therefore give
\[
 \mathcal H=\mathsf A_+^*R_rh_t,
 \qquad
 G_r(0,\cdot)=\mathsf A_+^*R_r\Ai,
 \qquad
 \|\mathcal H\|_2\leq\|h_t\|_2,
 \qquad
 \mathcal H(0)=q.
\]
The Airy decay and \(h_t'=th_t-\Ai\) justify differentiation under
\eqref{eq:technical-H-short}, integration by parts on \([-r,R]\),
and passage to \(R\to\infty\). They give
\begin{equation}
 \mathcal H'(v)
 =t\mathcal H(v)+\alpha b(v)-G_r(0,v).
 \label{eq:technical-H-ode-short}
\end{equation}
Using \(K_{\Ai}h_t=h_t\) and
\(\mathsf A_+\mathsf A_+^*=I-K_{\Ai}\), we obtain
\begin{align}
 \langle R_r\Ai,K_{\Ai}(I-R_r)h_t\rangle
 &=\langle R_r\Ai,(I-K_{\Ai})R_rh_t\rangle  =\int_0^\infty G_r(0,v)\mathcal H(v)\,dv.
 \label{eq:technical-inner-product-short}
\end{align}
The last integral is absolutely convergent because both factors
belong to \(L^2(\mathbb R_+)\).

Since \(\mathcal H\in L^2(\mathbb R_+)\), choose \(V_n\to\infty\)
such that \(\mathcal H(V_n)\to0\). Multiplying
\eqref{eq:technical-H-ode-short} by \(\mathcal H\), integrating on
\([0,V_n]\), and using \eqref{eq:technical-J-short} gives, as
\(n\to\infty\),
\[
 \langle R_r\Ai,K_{\Ai}(I-R_r)h_t\rangle
 =\frac{q^2}{2}+t\|\mathcal H\|_2^2+\alpha J_{t,r}.
\]
Finally, Lemma~\ref{lem:ht bounds} and \(|J_{t,r}|\leq C\) imply $ |q|+|\alpha|\leq C$ and $t\|\mathcal H\|_2^2
 \leq t\|h_t\|_2^2
 =1/2$. We conclude \eqref{eq:technical-variance-first}.
\end{proof}

\begin{proof}[Proof of Proposition \ref{prop:bounding variance}]
By Lemma~\ref{lem:exact identity variance}, \(N_{t,r}\) is almost surely finite and  
\(N_{t,r}\in L^2\). Let
$ q =\langle R_r\Ai,R_rh_t\rangle$.
By Lemma~\ref{lem:ht bounds}, $q =\int_{-r}^{\infty}h_t(x)\Ai(x)\,dx$ is finite, and hence  
 $q/2-q^2/2 
 \leq 1/8$. Applying Lemmas~\ref{lem:Airy leakage} and
\ref{lem:technical bounds variance} to  
Lemma~\ref{lem:exact identity variance}, we obtain absolute constants
\(C_1,C_2,C_3\in(0,\infty)\) such that, for \(t>0\) and \(r\geq1\),
\[
 \operatorname{Var}(N_{t,r})
 \leq
 2C_1\log(2+r)+ 1/8+C_2+C_3\log(2+r)\leq
 \left(
 2C_1+C_3+ (C_2+1/8)/\log3
 \right)\log(2+r).
\]  We conclude the proof.
\end{proof}

\section{Log-concavity of $\hsai_k(t)$: Proof of Proposition \ref{prop:log concave input}} 
\label{sec: proof of log concavity input}

In this section, we prove Proposition \ref{prop:log concave input}, which states that $\hsai_k(t)$ is   log-concave   for each $t\geq0$ and $k\in\mathbb N$. To prove this, we realize this  variable as the weak limit of a sequence of   variables, and show that each variable in this sequence is log-concave. We will define the prelimit variables in Section \ref{subsec: proof of log concavity input}, which are (rescaled) spectrums of certain matrix Ornstein–Uhlenbeck processes. The log-concavity of the spectrum will be shown in Section \ref{subsec:log concavity} and the weak convergence to $\hsai_k(t)$ will be shown in Section \ref{subsec: proof of scaling limit}.

\subsection{Proof of Proposition \ref{prop:log concave input}}
\label{subsec: proof of log concavity input}
In this subsection, we define the prelimit variables mentioned above and prove Proposition \ref{prop:log concave input} modulo the log-concavity of the spectrum (Proposition \ref{prop:log concave prelimit}) and the weak (scaling) convergence of these variables to $\hsai_k(t)$ (Proposition \ref{prop:scaling limit PP}).
 For each $m\in\mathbb N$ and $x\in\mathbb R^m$, we write
\[
 \Delta_m(x)=\prod_{1\le i<j\le m}(x_i-x_j).
\]

\begin{definition} 
\label{def:GUE-point-process}
For $m\in\mathbb N$, let $\mathsf{Herm}_m=\{H\in\mathbb C^{m\times m}:H=H^*\}$.
An $m\times m$ GUE matrix is a $\mathsf{Herm}_m$-valued random matrix $B_m$ with density
\begin{equation*}
 \frac{2^{m(m-1)/2}}{\pi^{m^2/2}}
 e^{-\operatorname{Tr}(H^2)}\,dH,\qquad\mbox{ where } dH=\prod_{i=1}^m dH_{ii}
 \prod_{1\leq i<j\leq m}
 d\operatorname{Re}H_{ij}\,d\operatorname{Im}H_{ij}.
\end{equation*}
The ordered eigenvalues $b_1>\cdots>b_m$ of $B_m$ then have joint density
\begin{equation*}
 \frac{1}{Z_m^{\rm GUE}}
 \mathbf1_{\{b_1>\cdots>b_m\}}
 e^{-\sum_{i=1}^m b_i^2}\Delta_m(b)^2,
 \qquad
 Z_m^{\rm GUE}
 =2^{-m(m-1)/2}\pi^{m/2}\prod_{j=1}^{m-1}j!. 
\end{equation*}
\end{definition}

\begin{definition}\label{def:prelimit} 
Let $n\in\mathbb N$ and $N=2n$. Let $A_N$ be an $N\times N$ Hermitian random matrix with the following GSE spectral law. Its spectrum is
$(a_1,a_1,\ldots,a_n,a_n)$, where $a_1>\cdots>a_n$ has density
\[
 \frac{1}{Z^{\rm GSE}_n}
 \mathbf1_{\{a_1>\cdots>a_n\}}
 e^{-\sum_{i=1}^n a_i^2}\Delta_n(a)^4,
 \qquad
 Z^{\rm GSE}_n
 =\frac{2^{-n^2}\pi^{n/2}}{n!}\prod_{j=1}^n(2j)!.
\]
Independently of $A_N$, let $B_N$ be an $N\times N$  GUE matrix as in Definition \ref{def:GUE-point-process}.
We now define
\[
 A_N(\tau)=e^{-\tau}A_N+ \sqrt{1-e^{-2\tau}}\, B_N,\qquad \tau>0,
\]
which is a matrix Ornstein–Uhlenbeck evolution started from $A_N$.
As  shown in Proposition \ref{prop:log concave prelimit}, $ A_N(\tau)$ has simple spectrum almost surely.
We write $Y_1^{(N)}(\tau)>\cdots>Y_N^{(N)}(\tau)$ to denote the ordered eigenvalues of $A_N(\tau)$.
\end{definition}
 \begin{proposition}
\label{prop:log concave prelimit}
Assume the same notation as in Definition \ref{def:prelimit}. Then $A_N(\tau)$ has simple spectrum almost surely, and the joint law on $\mathbb R^N$ of its ordered eigenvalues $Y_1^{(N)}(\tau)>\cdots>Y_N^{(N)}(\tau)$ is log-concave. 
\end{proposition}

The proof of Proposition \ref{prop:log concave prelimit} will be the focus of Section \ref{subsec:log concavity}.

\begin{proposition}\label{prop:scaling limit PP}
Let $\hsai$ be the pinned half-space Airy line ensemble as in Proposition \ref{Prop.HSAiryinfty}. 
For even $N\in\mathbb N$, let $Y_1^{(N)}(\tau)>\cdots>Y_N^{(N)}(\tau)$ be as in Definition \ref{def:prelimit}. Then for each fixed $t>0$ and $k\in\mathbb N$, 
\[Z_k^N(t)=\sqrt2N^{1/6}
 \left(Y_k^{(N)}(tN^{-1/3})-\sqrt{2N}\right)
 \Longrightarrow \hsai_k(t)\qquad\mbox{ as }N\to\infty.\]
\end{proposition}

The proof of Proposition \ref{prop:scaling limit PP} will be the focus of Section \ref{subsec: proof of scaling limit}.
 
\begin{proof}[Proof of Proposition \ref{prop:log concave input} assuming Propositions \ref{prop:log concave prelimit}  and \ref{prop:scaling limit PP} ]
We first assume $t>0$.  
For every even $N\geq k$ by Proposition~\ref{prop:log concave prelimit} with $\tau=tN^{-1/3}>0$ and Lemma~\ref{lem:log-concavity-properties}\textup{(vi)},  
$Y_k^{(N)}(tN^{-1/3})$ has a log-concave law.
Lemma~\ref{lem:log-concavity-properties}\textup{(v)} then shows that
$Z_k^N(t)$ has a log-concave law. Proposition~\ref{prop:scaling limit PP} and
Lemma~\ref{lem:log-concavity-properties}\textup{(vi)} therefore show that
$\hsai_k(t)$ has a log-concave law. If $t_n\downarrow0$, the path continuity
in Proposition~\ref{Prop.HSAiryinfty} gives
$\hsai_k(t_n)\to \hsai_k(0)$ almost surely, so  
Lemma~\ref{lem:log-concavity-properties}\textup{(vi)} proves
log-concavity at $t=0$.
\end{proof}

\subsection{Log-concavity of the prelimiting process: Proof of Proposition \ref{prop:log concave prelimit}} \label{subsec:log concavity}

We recall that a finite Borel measure $\mu$ on
$\mathbb R^m$ is called log-concave if 
\[
 \mu(\theta E+(1-\theta)F)
 \geq\mu(E)^\theta\mu(F)^{1-\theta}\qquad\mbox{ for all nonempty compact sets }E,F\mbox{ and }\theta\in(0,1).
\]
Here
$\theta E+(1-\theta)F=\{\theta x+(1-\theta)y:x\in E,\ y\in F\}$. 
A random vector is called log-concave when its law is log-concave. 
We call a function $f:\mathbb R^m\to[0,\infty)$   log-concave if
\[
 f(\theta x+(1-\theta)y)
 \geq f(x)^\theta f(y)^{1-\theta}\qquad\mbox{ for all }x,y\in\mathbb R^m,\,\theta\in(0,1).
\]

\begin{lemma}\label{lem:log-concavity-properties}
    The following statements hold.
    \begin{enumerate}[label=\textup{(\roman*)}] 
\item Finite products and nonnegative scalar multiples of log-concave functions are log-concave. If $f$ is log-concave and $T$ is affine,
then $f\circ T$ is log-concave. Any pointwise limit of log-concave functions that takes values in $[0,\infty)$ is log-concave.
Moreover, let  $C\subset\mathbb R^m$ be convex and let $\phi:C\to\mathbb R$ be a concave function. The zero extension to
$\mathbb R^m$ of $e^\phi$ is log-concave. In particular, indicators of convex sets are log-concave, as is $e^{-V}$ for every convex
$V:\mathbb R^m\to\mathbb R$. 
\item If $f:\mathbb R^m\to[0,\infty)$ is measurable and
log-concave, $g:\mathbb R^m\to[0,\infty)$ is continuous, and
$f=g$ Lebesgue-almost everywhere, then $g$ is log-concave.

\item Let $r,s\in\mathbb N$, let
$H:\mathbb R^{r+s}\to[0,\infty)$ be measurable and log-concave, and
suppose that
\[
 G(x)=\int_{\mathbb R^s}H(x,u)\,du<\infty,
 \qquad x\in\mathbb R^r.
\]
Then $G$ is log-concave. In particular, the Lebesgue marginals of every log-concave $H$ satisfying
$0<\int_{\mathbb R^{r+s}}H<\infty$ are finite and log-concave.

\item If $f:\mathbb R^m\to[0,\infty)$ is integrable and
log-concave, then the finite Borel measure
\[
 \mu(B)=\int_B f(x)\,dx
\]
is log-concave.

\item If a real random variable $X$ has a log-concave law, then
$aX+b$ has a log-concave law for $a,b\in\mathbb R$.

\item If finite log-concave Borel measures $\mu_j$ on $\mathbb R^m$
converge weakly to a finite Borel measure $\mu$, then $\mu$ is
log-concave. Moreover, every coordinate marginal of a finite log-concave Borel measure is
log-concave. 
    \end{enumerate}
\end{lemma}
\begin{proof}
The assertions in~\textup{(i)} follow immediately from the defining inequality.

For~\textup{(ii)}, let $D=\{f=g\}$, which has full Lebesgue measure.
For fixed $\theta\in(0,1)$, the set $\{(x,y):x,y,\theta x+(1-\theta)y\in D \}$
has full Lebesgue measure in $\mathbb R^{2m}$, hence is dense: its
complement is a union of inverse images of a null set under three
surjective linear maps, and these inverse images are null after a
linear change of variables and Fubini's theorem. Approximate an
arbitrary $(x,y)$ by points in this dense set, apply log-concavity of
$f$, and use continuity of $g$.

Part~\textup{(iii)} is the fiberwise form of Pr\'ekopa's marginal theorem; it follows from the Pr\'ekopa--Leindler inequality and is
recorded in \cite[Theorem~3.3 and Corollary~3.5]{BrascampLieb1976}. Its proof requires only the displayed finiteness, not joint
integrability. To justify the last sentence, let $0<\int H<\infty$ and retain the notation $G$. Tonelli's theorem shows that $G<\infty$ almost everywhere and that $A=\{x:0<G(x)<\infty\}$ has positive Lebesgue measure. If $G(x_0)=\infty$, set
\[
 h_n(u)=(H(x_0,u)\wedge n)\mathbf 1_{\{|u|\leq n\}}.
\]
For every $x\in A$, the Pr\'ekopa--Leindler inequality gives
\[
 G((x_0+x)/2)\geq
 \left(\int h_n\right)^{1/2}G(x)^{1/2}\longrightarrow\infty.
\]
This contradicts the almost-everywhere finiteness of $G$, since $(x_0+A)/2$ has positive Lebesgue measure.

For~\textup{(iv)}, fix nonempty compact sets $E,F$ and $\theta\in(0,1)$, put $C=\theta E+(1-\theta)F$, and apply the
Pr\'ekopa--Leindler inequality \cite[Theorem~3.3]{BrascampLieb1976} to
$f\mathbf 1_E$, $f\mathbf 1_F$, and $f\mathbf 1_C$. This gives $\mu(C)\geq \mu(E)^\theta\mu(F)^{1-\theta}$.

For~\textup{(v)}, let $T(x)=ax+b$. If $a\neq0$, then $T$ is an affine bijection and
$ T^{-1} (\theta E+(1-\theta)F )
 =\theta T^{-1}(E)+(1-\theta)T^{-1}(F)$
for compact $E,F\subset\mathbb R$, so the defining measure inequality for the law of $X$ proves~\textup{(v)}. If $a=0$, the image law is
the point mass $\delta_b$, which is log-concave directly from the definition.

For the first claim in~\textup{(vi)}, fix $E,F,C$ as above and write
 $K^\varepsilon
 =\{x:\operatorname{dist}(x,K)\leq\varepsilon\}$.
Choose $\varepsilon_k\downarrow0$ so that
$E^{\varepsilon_k}$, $F^{\varepsilon_k}$, and $C^{\varepsilon_k}$
are $\mu$-continuity sets. Such a common sequence exists because, for
each compact $K$, the boundaries of $K^\varepsilon$ are contained in
the pairwise disjoint level sets
$\{x:\operatorname{dist}(x,K)=\varepsilon\}$, and only countably many
of these can have positive $\mu$-measure. Since
 $\theta E^{\varepsilon_k}+(1-\theta)F^{\varepsilon_k}
 \subset C^{\varepsilon_k}$,
weak convergence and log-concavity of $\mu_j$ yield
$ \mu(C^{\varepsilon_k})
 \geq \mu(E^{\varepsilon_k})^\theta
        \mu(F^{\varepsilon_k})^{1-\theta}$.
Letting $k\to\infty$ and using continuity from above proves the first claim. For the second claim, let $\mu$ be such a measure on $\mathbb R^m$,
let $\pi_j:\mathbb R^m\to\mathbb R$ be a coordinate projection, and
fix nonempty compact sets $E,F\subset\mathbb R$ and
$\theta\in(0,1)$. If both inverse images have positive $\mu$-measure,
inner regularity provides compact sets
$A_n\subset\pi_j^{-1}(E)$ and $B_n\subset\pi_j^{-1}(F)$ such that $\mu(A_n)\to\mu(\pi_j^{-1}(E))$ and $\mu(B_n)\to\mu(\pi_j^{-1}(F))$.
Since
$ \theta A_n+(1-\theta)B_n
 \subseteq
 \pi_j^{-1}\!\left(\theta E+(1-\theta)F\right)$,
log-concavity of $\mu$, followed by $n\to\infty$, proves the claim;
if either inverse image has zero measure, it is immediate.\end{proof}

\begin{proof}[Proof of Lemma~\ref{lem:log-concavity improve concentration}]
Let $X=(Z-a)/\sigma$, which has a log-concave law by Lemma~\ref{lem:log-concavity-properties}\textup{(v)}, and set
$u_0=2\sqrt C$.   Since $C\ge1$, the hypothesis applies at $u_0$ and gives
$p=\mathbb P(|X|\ge u_0)\le 1/4$.
Therefore, for $u\ge u_0$, we have 
\[
 \mathbb P(|X|\ge u)
 \le p^{(u/u_0+1)/2}
 \le \frac12\exp\left(-\frac{\log2}{2\sqrt C}u\right)
 \le \sqrt3\exp\left(-\frac{\log3}{4\sqrt C}u\right),
\]
where the first inequality is an equality if $u=u_0$, and if $u>u_0$, it follows from the log-concave dilation inequality
\cite[Corollary~4, specialized to a linear polynomial]{Fradelizi2009}, applied with dilation factor $u/u_0$ (and is trivial if the
left side vanishes); the second inequality uses $p\leq1/4$ and $u_0=2\sqrt C$; and the last inequality uses $2\log2\ge\log3$.  If $1\le u<u_0$, the right side of \eqref{eq:logconcave-upgrade} is at least one.
\end{proof}

\begin{lemma}\label{lem:density OU log concave}
For $d\in\mathbb N$, let $\mathsf W_d=\{x\in\mathbb R^d:x_1>\cdots>x_d\}$. Let $X=(X_1,\ldots,X_d)$ consist of independent Ornstein–Uhlenbeck processes, i.e., solutions of $dX_i(s)=-X_i(s)\,ds+dB_i(s)$, where $B_1,\ldots,B_d$ are independent standard Brownian motions. Write $\mathbb P_x$ for their law with
$X(0)=x$, and, for $x\in\mathsf W_d$, let $\sigma_{\mathsf W_d} =\inf\{s\geq0:X(s)\notin\mathsf W_d\}$, with $\inf\varnothing=\infty$. Then for every $t>0$, the killed transition measure $ \mathbb P_x (X(t)\in dy,\,t<\sigma_{\mathsf W_d} )$ has Lebesgue density
\[
 q_t^{\mathsf W_d}(x,y)
 =\det\left[
 \frac{1}{\sqrt{\pi(1-e^{-2t})}} 
 \exp\left(-\frac{(y_i-e^{-t}x_j)^2}{1-e^{-2t}}\right)
 \right]_{i,j=1}^d,
 \qquad x,y\in\mathsf W_d.
\]
Extended to be zero off $\mathsf W_d\times\mathsf W_d$, the function $(x,y)\mapsto q_t^{\mathsf W_d}(x,y)$ is continuous and log-concave on
$\mathbb R^{2d}$.  
\end{lemma}
\begin{proof}
    The scalar solution of $dX(s)=-X(s)\,ds+dB(s)$ from $X(0)=u$ satisfies
    \[
 X(s)=e^{-s}u+\int_0^s e^{-(s-r)}\,dB(r),
 \qquad
 \operatorname{Var}(X(s)\mid X(0)=u)
 =\frac{1-e^{-2s}}{2}.
\]
Hence the transition density of the Ornstein–Uhlenbeck process is 
\[
 p_s(u,v)=\frac{1}{\sqrt{\pi(1-e^{-2s})}}
 \exp\left(-\frac{(v-e^{-s}u)^2}{1-e^{-2s}}\right),
 \qquad s>0.
\]
 For fixed $t>0$ and $x,y\in\mathsf W_d$, we denote by $\mathbb P_{x,y}^{\mathrm{br}}$ the bridge law of independent Ornstein–Uhlenbeck processes from $x$ to $y$, i.e., the law of the process
\[
 X^{x,y,\mathrm{br}}(u) =\frac{\sinh(t-u)}{\sinh t}x
  +\frac{\sinh u}{\sinh t}y+H(u)-\frac{\sinh u}{\sinh t}H(t),
 \qquad u\in[0, t],
\]
where $H(u)=\int_0^u e^{-(u-r)}\,d\mathbf B(r)$ and $\mathbf B$ is the standard $d$-dimensional Brownian motion. This bridge process $X^{x,y,\mathrm{br}}$ has continuous paths and  has finite-dimensional density at times  $0<s_1<\cdots<s_k<t$ and locations  $(z_1,\ldots,z_k)$ given by
$ \rho_{s_1}(x,z_1)
 \rho_{s_2-s_1}(z_1,z_2)\cdots
 \rho_{t-s_k}(z_k,y) /\rho_t(x,y)$, where   $ \rho_s(u,v)=\prod_{i=1}^d p_s(u_i,v_i)$.

For clarity, we split the rest of the proof into three steps.\\

\noindent\textbf{Step 1. }  Define 
\[
\begin{aligned}
E_{\mathrm{op}} &=\{X(s)\in\mathsf W_d:0\leq s\leq t\},\qquad\qquad E_{\mathrm{cl}}=\{X(s)\in\overline{\mathsf W_d}:0\leq s\leq t\},\\
  E_m&=\{X(jt/2^m)\in\mathsf W_d:0\leq j\leq 2^m\},\qquad m\in\mathbb N. 
\end{aligned}
\]    In this step, we show that 
\begin{equation}\label{eq:first claim gq3we}
E_{\mathrm{op}}
 \subseteq\cap_{m\in\mathbb N}E_m
 \subseteq E_{\mathrm{cl}},\qquad\qquad
    \mathbb P_{x,y}^{\mathrm{br}}
 \left(E_{\mathrm{cl}}\setminus E_{\mathrm{op}}\right)=0  .
\end{equation}
The first assertion follows from the path continuity and the union of the dyadic grids.   
The second assertion is vacuous for $d=1$. We denote 
$ r(s)= (e^{2s}-1)/2$ and $T=r(t)>0$. 
In view of $ e^sX_i(s)=x_i+\int_0^s e^u\,dB_i(u)$, after the deterministic time change $s\mapsto r(s)$ and
conditioning on $X(t)=y$,  for every $i\in\llbracket1,d-1\rrbracket$, the adjacent gap
$e^s(X_i(s)-X_{i+1}(s))$, parametrized on $[0,T]$, is a Brownian bridge with variance parameter $2$ and positive endpoints
$a_i=x_i-x_{i+1}$ and $ b_i=e^t(y_i-y_{i+1})$. 
For a Brownian bridge $G$ of duration $T$, variance parameter $\sigma^2$, and endpoints $a,b>0$, the reflection principle gives, for
$c<\min(a,b)$,
\[
 \mathbb P\left(\min_{ u\in[0, T]}G(u)>c\right)
 =1- e^{- 2(a-c)(b-c)/(\sigma^2T) }.
\]
Letting $c\uparrow0$ through negative values and $c\downarrow0$ through positive values shows that the minimum has no atom at zero.
Therefore, an adjacent gap has probability zero of being nonnegative throughout $[0,t]$ while vanishing somewhere. Taking the finite union over $i\in\llbracket1,d-1\rrbracket$ concludes $ \mathbb P_{x,y}^{\mathrm{br}}
 \left(E_{\mathrm{cl}}\setminus E_{\mathrm{op}}\right)=0$. \\

\noindent\textbf{Step 2. } For $x,y\in\mathbb R^d$, $t>0$, $m\in\mathbb N$, with $z_0=x\in\mathbb R^d$ and $z_{2^m}=y\in\mathbb R^d$, we define 
\[\begin{aligned}
 K^m_t(x,y) &=
 \int_{(\mathbb R^d)^{2^m-1}}
 \mathbf1_{\mathsf W_d^{2^m+1}}(z_0,\ldots,z_{2^m})
 \prod_{r=1}^{2^m}\prod_{i=1}^d
 p_{t/2^m}(z_{r-1,i},z_{r,i})\,
 dz_1\cdots dz_{2^m-1}\\
& =
 \begin{cases}
 \rho_t(x,y)\mathbb P_{x,y}^{\mathrm{br}}(E_m),
     & x,y\in\mathsf W_d,\\
 0,  & \text{otherwise}.
 \end{cases} \end{aligned}
\]
 In this step, we show that $K^m_t(x,y)$ is log-concave on $(x,y)\in\mathbb R^{2d}$.   
Indeed, the integrand of $ K^m_t(x,y)$ is log-concave in $(z_0,\ldots,z_{2^m})$: the set $\mathsf W_d^{2^m+1}$ is convex, and the logarithm of each Gaussian factor is a constant minus a convex quadratic form. Moreover, we have   $ 0\leq K^m_t(x,y)
 \leq\prod_{i=1}^d p_t(x_i,y_i)<\infty$ by removing the chamber indicator and using Chapman--Kolmogorov.  
By Lemma~\ref{lem:log-concavity-properties}\textup{(i)} and \textup{(iii)}, $K^m_t$ is log-concave on $\mathbb R^{2d}$. \\

\noindent\textbf{Step 3. }
In this step, we conclude the proof of the lemma.
 When $x,y\in\mathsf W_d$, define
\[
 h_t(x,y)    = \rho_t(x,y)\mathbb P_{x,y}^{\mathrm{br}}(E_{\mathrm{op}}),
\]
and let $h_t(x,y) =0$ otherwise. Observe from \eqref{eq:first claim gq3we} that $\lim_{m\to\infty}K^m_t(x,y)=h_t(x,y)$ for each fixed $x,y\in\mathbb R^d$. By the claim of the last step and Lemma~\ref{lem:log-concavity-properties}\textup{(i)},   $h_t(x,y)$ is log-concave on $(x,y)\in\mathbb R^{2d}$. For every bounded Borel function $\varphi$,  bridge disintegration gives
\[
 \mathbb E_x\left[
 \varphi(X(t))\mathbf1_{\{t<\sigma_{\mathsf W_d}\}}
 \right]
 =\int_{\mathsf W_d}\varphi(y)
 h_t(x,y)\,dy.
\]
Thus $h_t(x,\cdot)$ is a density of the killed transition measure.

We recall that $q_t^{\mathsf W_d}(x,y)=\det[p_t(x_i,y_j)]_{i,j=1}^d$ for $x,y\in\mathsf W_d$. We next claim that
\begin{equation}\label{eq:tyserth}
 \mathbb P_x\left(
 X_i(t)\in E_i \mbox{ for  }1\leq i\leq d,\,
 t<\sigma_{\mathsf W_d}\right)
 =
 \det\left[\int_{E_i}p_t(x_j,v)\,dv\right]_{i,j=1}^d
 =
 \int_{E_1\times\cdots\times E_d}q_t^{\mathsf W_d}(x,y)\,dy,
\end{equation}
for each bounded open intervals $E_1,\ldots,E_d$ such that every point of $E_i$ is larger than every point of $E_{i+1}$.
Indeed,
the $d$-dimensional Ornstein–Uhlenbeck process is a time-homogeneous strong Markov process with continuous paths and continuous transition density $\rho_t$, so the
Karlin--McGregor formula \cite[Theorem~1]{KarlinMcGregor1959} applies.   Since a continuous noncolliding path that starts in $\mathsf W_d$ cannot leave that chamber, only the identity endpoint assignment in the Karlin--McGregor alternating sum can occur, and we conclude the proof of  \eqref{eq:tyserth}.

For fixed $x\in\mathsf W_d$, continuity of $p_t$ implies continuity of $q_t^{\mathsf W_d}(x,\cdot)$. In view of \eqref{eq:tyserth}, its integral over every ordered open rectangle is nonnegative. Hence $q_t^{\mathsf W_d}(x,y)\geq0$ on $\mathsf W_d$.
  Moreover, the Leibniz expansion gives $\int_{\mathsf W_d}|q_t^{\mathsf W_d}(x,y)|\,dy\leq d!$, because each of its $d!$ product terms has integral at most one.
Since the ordered open rectangles with rational endpoints form an intersection-stable determining basis of $\mathsf W_d$,  and the finite
measure with density $q_t^{\mathsf W_d}(x,\cdot)$ and the killed transition measure agree on this determining class, we know that they agree on every Borel subset of $\mathsf W_d$. By the uniqueness of densities, we have
\[
 h_t(x,y)=q_t^{\mathsf W_d}(x,y)
 \quad\text{for Lebesgue-almost every }y\in\mathsf W_d
\]
for every $x\in\mathsf W_d$. Therefore,  $h_t=q_t^{\mathsf W_d}$  Lebesgue-almost everywhere on $\mathbb R^{2d}$.

The function $q_t^{\mathsf W_d}$ is continuous on $\mathbb R^{2d}$. Indeed, it is continuous on $\mathsf W_d\times\mathsf W_d$. At a point of $\overline{\mathsf W_d}\times\overline{\mathsf W_d}$ outside the open product, at least one endpoint lies on $\partial\mathsf W_d$; two
coordinates of that endpoint coincide, so two columns (for $x$) or two rows (for $y$) of the limiting determinant coincide, and the
determinant tends to zero. If either endpoint lies outside $\overline{\mathsf W_d}$, the zero extension vanishes in a neighborhood.
The above argument also shows that it is nonnegative everywhere. By Lemma~\ref{lem:log-concavity-properties}\textup{(ii)}, we have that $q_t^{\mathsf W_d}$ is log-concave.  
\end{proof}

\begin{proof}[Proof of Proposition \ref{prop:log concave prelimit}] Recall $N=2n$ and $n\in\mathbb N$. Throughout the proof, we fix $\tau>0$ and let $u=e^{-\tau}$, $v=\sqrt{1-u^2}$. 

We first prove the simplicity of the spectrum. For every deterministic $H\in\mathsf{Herm}_N$, the matrix $uH+vB_N$ has a Lebesgue density on $\mathsf{Herm}_N$ because $v>0$. The set of Hermitian matrices having a repeated eigenvalue is the zero set of the discriminant of the characteristic polynomial. This is
a nonzero real polynomial in the $N^2$ real coordinates of a Hermitian matrix, and hence its zero set has Lebesgue measure zero. Therefore,
for every deterministic $H$,
$ \mathbb P  (
 uH+vB_N\ \text{has a repeated eigenvalue} )=0$.
Integrating this identity against the law of $A_N$ proves that $A_N(\tau)$ has simple spectrum almost surely.
We split the rest of the proof into five steps.\\

\noindent\textbf{Step 1. }
We write $\lambda^{\downarrow}(H)$ to denote the (weakly) decreasingly ordered eigenvalue vector of a matrix $H$.
  Suppose $a$ has probability density $\rho(a)\propto\mathbf 1_{\{a_1>\cdots>a_n\}}e^{-\sum_i a_i^2}\Delta_n(a)^4$, 
  let $F_0(a)=\mathrm{diag}[a_1,a_1,\dots,a_n,a_n]$,
  and let $B_N$ be an independent GUE random matrix. In this step, we show that 
  \begin{equation}\label{eq:claim step 1 prop log conc}
      \lambda^\downarrow(uF_0(a)+v\,B_N)\overset{d}{=}\lambda^\downarrow (A_N(\tau)).
  \end{equation}
Indeed, let $\nu_{\rm GUE}$ denote the law of $B_N$. For each bounded Borel function $\Phi:\mathbb R^N\to\mathbb R$, define 
\[
 K_\Phi(H)
 =
 \int_{\mathsf{Herm}_N}
 \Phi\!\left(
 \lambda^\downarrow(uH+vG)
 \right)\nu_{\rm GUE}(dG).
\]
 The map $(H,G)\mapsto \Phi  (\lambda^\downarrow(uH+vG) )$ is bounded and Borel, hence $K_\Phi$ is Borel. If a deterministic $H$ has ordered spectrum $(c_1,c_1,\ldots,c_n,c_n)$, choose a unitary $V$ such that $H=V(F_0(c))V^*$. Unitary invariance of $\nu_{\rm GUE}$ gives
$ K_\Phi(H)=K_\Phi(F_0(c))$. Independence of $A_N$ and $B_N$ gives
 \[ 
 \mathbb E\Phi\!\left(
 \lambda^\downarrow(A_N(\tau))
 \right) =
 \mathbb E K_\Phi(A_N)  =
 \int_{\mathbb R^n}K_\Phi(F_0(a))\rho(a)\,da  =
 \mathbb E\Phi\!\left(
\lambda^\downarrow(uF_0(a)+vB_N)
 \right),
\] 
and we conclude the proof of \eqref{eq:claim step 1 prop log conc}.\\

\noindent\textbf{Step 2. } For $\varepsilon>0$, define $ \mathsf D_\varepsilon
 =
 \{a\in\mathbb R^n:
   a_i-a_{i+1}>\varepsilon,\ 1\le i<n\}$. For $\varepsilon\geq0$, define
\[
 \iota_\varepsilon(a)
 =
 (a_1+\varepsilon/2,a_1-\varepsilon/2,\ldots,
  a_n+\varepsilon/2,a_n-\varepsilon/2).
\]
Then $\mathsf D_\varepsilon$ is nonempty and
$\iota_\varepsilon(\mathsf D_\varepsilon)\subset\mathsf W_N$.
Let $a^{(\varepsilon)}$ have probability density $\rho_{\varepsilon}$ on $\mathsf D_\varepsilon$ proportional to
$ \Delta_N(\iota_\varepsilon(a)) e^{-|\iota_\varepsilon(a)|^2/2}$.
This defines a probability law: the integrand is strictly positive on the nonempty open set $\mathsf D_\varepsilon$, while its zero extension is bounded by a polynomial times a Gaussian. 
 We define the map $F_\varepsilon(a)=\operatorname{diag}[\iota_\varepsilon(a)]$, which agrees with the map $F_0$ defined in Step 1 when $\varepsilon=0$. In this step, we show that 
\begin{equation}\label{eq:claim step2 log conc}
 F_\varepsilon(a^{(\varepsilon)})\Rightarrow F_0(a)\qquad\mbox{ as } \varepsilon\downarrow0,
\end{equation}
where $a$ has probability density $\rho$   as in Step 1.
Indeed, we can write $\rho=h/Z$ and $\rho_{\varepsilon}=h_{\varepsilon}/Z_{\varepsilon}$, where  
\[\begin{aligned}
 h_\varepsilon(a)
& =
 \mathbf 1_{\mathsf D_\varepsilon}(a)
 e^{-\sum_i a_i^2-n\varepsilon^2/4}
 \prod_{i<j}(a_i-a_j)^2
 \left((a_i-a_j)^2-\varepsilon^2\right),\\
 h(a)&= \mathbf 1_{\{a_1>\cdots>a_n\}}
 e^{-\sum_i a_i^2}\Delta_n(a)^4,
\end{aligned}\]
$Z=\int h$, and $Z_{\varepsilon}=\int h_{\varepsilon}$.
We have $0\le h_\varepsilon\le h$ and $ h_\varepsilon(a)\to h(a)$ pointwise as $\varepsilon\downarrow0$. 
Since $0<\int_{\mathbb R^n}h(a)\,da=Z <\infty$, dominated convergence gives $\|h_\varepsilon-h\|_{L^1}\to0$. Hence $Z_\varepsilon\to Z>0$ and
\[
 \|\rho_\varepsilon-\rho\|_{L^1}
 \le
 \frac{\|h_\varepsilon-h\|_{L^1}}{Z_\varepsilon}
 +\left|\frac1{Z_\varepsilon}-\frac1Z\right|Z
 \longrightarrow0.
\]
For every bounded function $\Psi:\mathsf{Herm}_N\to\mathbb R$ that is Lipschitz with respect to the operator norm,
\[
\begin{aligned}
 \left|\mathbb E\Psi (F_\varepsilon(a^{(\varepsilon)}) )-\mathbb E\Psi(F_0(a))\right|
 &\le
 \|\Psi\|_\infty\|\rho_\varepsilon-\rho\|_{L^1}
 +\operatorname{Lip}_{\rm op}(\Psi)
  \sup_{a\in\mathbb R^n}
  \|F_\varepsilon(a)-F_0(a)\|_{\rm op} \\
 &=
 \|\Psi\|_\infty\|\rho_\varepsilon-\rho\|_{L^1}
 +\frac{\varepsilon}{2}\operatorname{Lip}_{\rm op}(\Psi)
 \longrightarrow0,
\end{aligned}
\]
and we conclude the proof of \eqref{eq:claim step2 log conc}.\\

\noindent\textbf{Step 3. } 
For $\varepsilon>0$, define $\nu_\varepsilon$ to be the law of
\[
\lambda^{\downarrow}(u F_\varepsilon(a^{(\varepsilon)}) +v B_N),
\]
where $a^{(\varepsilon)}$ and $B_N$ are independent. 
In this step, we prove that $\nu_\varepsilon$ has a density
proportional to
\begin{equation}
 \mathbf 1_{\mathsf W_N}(y)\Delta_N(y)
 \int_{\mathsf D_\varepsilon}
 q_\tau^{\mathsf W_N}(\iota_\varepsilon(a),y)
 e^{-|\iota_\varepsilon(a)|^2/2}\,da.
 \label{eq:split-evolved-density new}
\end{equation}
Indeed, for each fixed $x\in\mathsf W_N$, let $Q_\tau(x,\cdot)$ denote the density on $\mathsf W_N$ of
\[
 \lambda^{\downarrow}(u\operatorname{diag}[x]+vB_N).
\]
The density of this matrix at $H\in\mathsf{Herm}_N$ is proportional to $\exp (-v^{-2}  \operatorname{Tr} ((H-u\operatorname{diag}[x])^2) )$. By  the Weyl integration formula
\cite[Section~1.3, (1.27)]{Forrester2010} and the HCIZ identity
\cite[(3.4)]{ItzyksonZuber1980}, applied with the distinct
coordinates $x,y\in\mathsf W_N$ and the parameter $2u/v^2>0$,
\begin{equation*} 
    \begin{split}
         Q_\tau(x,y)
 &=
 \frac{v^{-N^2 }}{Z_N^{\rm GUE}}
 e^{-(|y|^2+u^2|x|^2)/v^2}\Delta_N(y)^2
 \frac{\prod_{k=1}^{N-1}k!}{(2u/v^2)^{N(N-1)/2}}
 \frac{\det[e^{2ux_jy_i/v^2}]_{i,j=1}^N}
      {\Delta_N(x)\Delta_N(y)}\\
      &= u^{-N(N-1)/2}\frac{\Delta_N(y)}{\Delta_N(x)}
 q_\tau^{\mathsf W_N}(x,y).
    \end{split}
\end{equation*}
Mixing this against $\rho_{\varepsilon}\propto \Delta_N(\iota_\varepsilon(a)) e^{-|\iota_\varepsilon(a)|^2/2}$ gives \eqref{eq:split-evolved-density new}.\\

\noindent\textbf{Step 4. }  
In this step, we show that $\nu_{\varepsilon}$ is a log-concave probability measure. By Lemma~\ref{lem:density OU log concave}, the zero
extension of $q_\tau^{\mathsf W_N}$ is jointly log-concave. Therefore
\[
 H_\varepsilon(a,y)
 =
 \mathbf 1_{\mathsf D_\varepsilon}(a)
 q_\tau^{\mathsf W_N}(\iota_\varepsilon(a),y)
 e^{-|\iota_\varepsilon(a)|^2/2}
\]
is jointly log-concave on $\mathbb R^{n+N}$ by Lemma~\ref{lem:log-concavity-properties}\textup{(i)}. 
Moreover,
\[
\begin{aligned}
 0 <
 \int_{\mathbb R^{n+N}}H_\varepsilon(a,y)\,da\,dy  =
 \int_{\mathsf D_\varepsilon}
 e^{-|\iota_\varepsilon(a)|^2/2}
 \left(
 \int_{\mathsf W_N}
 q_\tau^{\mathsf W_N}(\iota_\varepsilon(a),y)\,dy
 \right)da  \le
 \int_{\mathsf D_\varepsilon}
 e^{-|\iota_\varepsilon(a)|^2/2}\,da
 <\infty,
\end{aligned}
\]
where the strict positivity  follows from Step~3: if
$\int_{\mathbb R^{n+N}}H_\varepsilon(a,y)\,da\,dy=0$, then Tonelli's theorem would force the integral over $a$ in
\eqref{eq:split-evolved-density new} to vanish for Lebesgue-almost every $y$; hence the function displayed there would vanish almost everywhere, contradicting that it is proportional to a probability density.

Lemma~\ref{lem:log-concavity-properties}\textup{(iii)} shows that
\[
 G_\varepsilon(y)
 =
 \int_{\mathbb R^n}H_\varepsilon(a,y)\,da
\]
is finite and log-concave on $\mathbb R^N$.
 Furthermore,
\[
 \log\Delta_N(y)=\sum_{i<j}\log(y_i-y_j),
 \qquad y\in\mathsf W_N,
\]
is concave. Since $\mathsf W_N$ is convex, the zero extension of
$\mathbf 1_{\mathsf W_N}\Delta_N$ is log-concave by
Lemma~\ref{lem:log-concavity-properties}\textup{(i)}. Thus the density
in \eqref{eq:split-evolved-density new} is log-concave, and
Lemma~\ref{lem:log-concavity-properties}\textup{(iv)} implies that
$\nu_\varepsilon$ is a log-concave probability measure.

\noindent\textbf{Step 5. } In this step, we conclude the proof of the proposition. 
Indeed, \eqref{eq:claim step2 log conc} gives 
\[
 (F_{\varepsilon}(a^{(\varepsilon)}),B_N)\Rightarrow(F_0(a),B_N).
\]
By  Weyl's inequality,
\[
 \max_{1\le j\le N}
 |\lambda_j^\downarrow(H)-\lambda_j^\downarrow(K)|
 \le \|H-K\|_{\rm op},
\]
and hence the map
$ (H,G)\longmapsto
 \lambda^\downarrow(uH+vG)$
is continuous. The continuous mapping theorem yields 
\[
 \nu_\varepsilon\Longrightarrow
 \mathrm{Law}\!\left(
 \lambda^\downarrow(uF_0(a)+vB_N)
 \right).
\]
By Step 4, each $ \nu_\varepsilon$ is log-concave. Using
Lemma~\ref{lem:log-concavity-properties}\textup{(vi)}, the limiting law is log-concave.
\end{proof}

\subsection{The scaling limit convergence: Proof of Proposition \ref{prop:scaling limit PP}} \label{subsec: proof of scaling limit}
In this subsection, we prove Proposition \ref{prop:scaling limit PP}, which states that for each fixed $t>0$ and $k\in\mathbb N$, the rescaled $k$th eigenvalue of $A_N(tN^{-1/3})$ (as in Definition \ref{def:prelimit}) converges weakly to the curve height $\hsai_k(t)$. We begin by considering the union of such points for all admissible $k$. We realize that the prelimit is a Pfaffian point process on $\mathbb R$ (with correlation kernel given by a rescaled version of \cite{ForresterNagaoHonner1999}), and the proposed limit is also a Pfaffian point process on $\mathbb R$ (as shown in Proposition \ref{Prop.HSAiryinfty}). We adopt the framework in \cite{dimitrov2026half}. We first prove the weak convergence of Pfaffian point processes in vague topology, which follows from showing the local uniform convergence on $(x,y)\in\mathbb R^2$ of the corresponding Pfaffian correlation kernels. Then to establish the weak convergence of the $k$-th ordered point in these point processes, we prove an upper tail bound which rules out the particles escaping to infinity.

\begin{definition}\label{def:kernel KN}
For $t>0$,  $N,m\in\mathbb N$, and $j,n\in\mathbb Z_{\geq0}$,
let $q_N=e^{-tN^{-1/3}}$, $a_N=(1-q_N^2)/(2(1+q_N^2))$,
$d_N=(\sqrt2N^{1/6})^{-1}$, and $h_N(x)=\sqrt{2N}+d_Nx$,
\begin{equation*}
    \begin{aligned}
    H_j(x)&=(-1)^je^{x^2}\frac{d^j}{dx^j}e^{-x^2},\\
         \phi_j(x)
& =
 (\sqrt\pi\,2^jj!)^{-1/2}e^{-x^2/2} H_j(x), \\
  K_m^{\mathrm{Her}}(u,v)
 &=\sum_{j=0}^{m-1}\phi_j(u)\phi_j(v),\\
 \psi_{j,N}(x)
 &=
 2^{-1/4}N^{1/12}\phi_j\left(h_N(x)\right),\\
   J_{n,N}(x)
& =
 \int_x^\infty
 e^{-a_N(h_N(v)^2-h_N(x)^2)}
 \psi_{n,N}(v)\,dv,
    \end{aligned}
\end{equation*}
where $H_j(x)$ is known as the physicists' Hermite polynomial.  
 
For $t>0$, $x,y\in\mathbb R$, and even $N\in\mathbb N$, we define the kernel 
    \begin{equation*}
        K^{t,N}(x,y)=\begin{bmatrix}
        K^{t,N}_{11}(x,y) & K^{t,N}_{12}(x,y)\\
        K^{t,N}_{21}(x,y) &K^{t,N}_{22}(x,y)
    \end{bmatrix},
    \end{equation*}
where

\begin{align}
    \label{eq:K11Ndef}K_{11}^{t,N}(x,y)=&
- \frac{N^{-1/3}}{1+q_N^2}
 \sum_{k=0}^{N/2-1}e^{-2t(2k+1)N^{-1/3}}  
 \left(J_{N-1-2k,N}(x)\psi_{N-1-2k,N}(y)
 -
 \psi_{N-1-2k,N}(x)J_{N-1-2k,N}(y)\right),\\
\label{eq:K12Ndef}         K_{12}^{t,N}(x,y)=&
d_NK_N^{\mathrm{Her}}(h_N(x),h_N(y))  -\frac{q_N^2}{1+q_N^2} 
 \psi_{N,N}(y)J_{N-1,N}(x),\\ 
\label{eq:K22Ndef}  K_{22}^{t,N}(x,y)
 =&  
 -2d_N^3(x-y)
 \sum_{m=N+1}^\infty
 q_N^{2m-2N}K_m^{\mathrm{Her}}(h_N(x),h_N(y)),
\end{align}
and $K_{21}^{t,N}(x,y)=- K_{12}^{t,N}(y,x)$.
\end{definition}

\begin{lemma}\label{lem:PPP structure}
Fix $t\in(0,\infty)$ and even $N\in\mathbb N$. For $1\leq k\leq N$, let
$Z_k^N(t)$ be as in Proposition~\ref{prop:scaling limit PP}, and define
\[
 M^{t,N}(\omega,A)
 =\sum_{k=1}^N\mathbf 1\{Z_k^N(t,\omega)\in A\}.
\]
Then $M^{t,N}$ is a Pfaffian point process on $\mathbb R$, in the
sense of \cite[Definition~5.6]{dimitrov2026half}, with reference
measure $\mathrm{Leb}$ and correlation kernel $K^{t,N}$ from
Definition~\ref{def:kernel KN}.
\end{lemma}
\begin{proof}
In this proof, we let $\tau_N=tN^{-1/3}$, $q=q_N=e^{-\tau_N}$, $d_N=(\sqrt2N^{1/6})^{-1}$, and $E_N=e^{N\tau_N}=q^{-N}$. Then 
$h_N(x)=\sqrt{2N}+d_Nx$. For clarity, we split the rest of the proof into three steps.\\

\noindent\textbf{Step 1. } In this step, we show that $M^{t,N}$ is a Pfaffian point process on $\mathbb R$, with correlation kernel 
\begin{equation}\label{eq:scaled-FNH-kernel}
 \mathsf K_N(x,y)
 =
 \begin{bmatrix}
 -E_N^{-2}D_N^{\mathrm F}(h_N(x),h_N(y))
 &d_NS_N(h_N(y),h_N(x))\\
 -d_NS_N(h_N(x),h_N(y))
 &d_N^2E_N^2\widetilde I_N(h_N(x),h_N(y))
 \end{bmatrix},
\end{equation} 
and with reference
measure $\mathrm{Leb}$,
where $S_N,\widetilde I_N,D_N^{\mathrm F}$ are the equal-time scalar kernels from \cite[Section~2]{ForresterNagaoHonner1999}, evaluated at
time $\tau_N$.  
Indeed, conditional
on $A_N=H$, the matrix $A_N(\tau_N)$ has density
\begin{equation}\label{eq:conditional-SU-density new}
 G\longmapsto
 c_{N,\tau_N}
 \exp\!\left(-\frac{\operatorname{Tr}((G-qH)^2)}{1-q^2}\right)
\end{equation}
with respect to Lebesgue measure on $\mathsf{Herm}_N$, where $c_{N,\tau_N}$ is independent of $H$. Take final symmetry $\beta=2$ and initial symmetry $\beta'=4$ in \cite[(1.3) and (2.7)]{ForresterNagaoHonner1999}. After eliminating the delta constraints in their $\beta'=4$ initial law, its distinct eigenvalues have density proportional to
 $e^{-\sum_{j=1}^{N/2}a_j^2}\Delta_{N/2}(a)^4$
and each eigenvalue has multiplicity two. This is precisely the spectral law in Definition~\ref{def:prelimit}. Moreover, \eqref{eq:conditional-SU-density new} is the conditional matrix density in \cite[(1.3)]{ForresterNagaoHonner1999} for $\beta=2$.  Since $B_N$ is unitarily invariant, the conditional law of the
eigenvalues in \eqref{eq:conditional-SU-density new} depends on $H$
only through its spectrum. Consequently, $H$ may be replaced by a
self-dual representative having the same doubled spectrum, and the
cited $\beta'=4$ transition formula applies.

Let 
$ J=\begin{bmatrix}0&1\\-1&0\end{bmatrix}$ 
and
$ R=\operatorname{diag}[\im,-\im]$. 
Taking the parameters denoted by $n,m$ in
\cite{ForresterNagaoHonner1999} to be $0,r$, respectively, and taking
$\tau_1=\tau_N/2$ and $\tau_2=\tau_N$ in
\cite[(2.8) and (2.24)--(2.26)]{ForresterNagaoHonner1999}
gives the equal-time correlation functions for
$1\leq r\leq N$. The matrix denoted $Z$ in \cite[(2.26)]{ForresterNagaoHonner1999} has diagonal blocks $-J$. With our convention $\operatorname{Pf}(J)=1$, the quaternion determinant is obtained using blocks $J$ instead, so the one-point function is $S_N(u,u)$.  Under our convention, the Pfaffian kernel is
\begin{equation} \label{eq:yethnm}
 K_N^0(u,v)
 =R\ J\begin{bmatrix}
  S_N(u,v)&\widetilde I_N(u,v)\\
  D_N^{\mathrm F}(u,v)&S_N(v,u)
 \end{bmatrix} R^{\mathsf T}
 =
 \begin{bmatrix}
  -D_N^{\mathrm F}(u,v)&S_N(v,u)\\
  -S_N(u,v)&\widetilde I_N(u,v)
 \end{bmatrix}.
\end{equation} 
We next show that this also gives the correct Pfaffian correlation kernel for $r>N$. Indeed, since the process
has  $N$ points, its $r$th  correlation function vanishes for $r>N$.  It suffices to show that the $2r\times2r$ block matrix
 $[K_N^0(u_a,u_b)]_{a,b=1}^r$
is singular, so that its Pfaffian vanishes.  Set $f_k(u)  = \sqrt{w(u)}e^{\gamma_k\tau_N}R_k(u;\tau_N)$ for $0\leq k<N$. By \cite[(2.20) and (2.22)]{ForresterNagaoHonner1999}, for fixed $u_1,\ldots,u_r$ there exist row vectors
$V_0,\ldots,V_{N-1}\in\mathbb C^{2r}$, independent of $i$, such that
\[
 \operatorname{row}_{2i-1}
 \left[K_N^0(u_a,u_b)\right]_{a,b=1}^r
 =
 \sum_{k=0}^{N-1}f_k(u_i)V_k,
 \qquad 1\leq i\leq r.
\]
Indeed, both $-D_N^{\mathrm F}(u_i,u_j)$ and $S_N(u_j,u_i)$ depend on $u_i$ only through the $N$ functions $f_k(u_i)$.  Hence, if
$r>N$, these $r$ rows are linearly dependent. We conclude that \eqref{eq:yethnm} gives the Pfaffian kernel for all $r\in\mathbb N$.

Under the change of variables $u=h_N(x)$, the $r$-point correlation function acquires  Jacobian $d_N^r$. Denote $\mathsf G_N=\operatorname{diag}[E_N^{-1},d_NE_N]$, and we have  $\det\mathsf G_N=d_N$. Therefore, the rescaled process $M^{t,N}$ formed by $Z_k^N(t)$,
$1\leq k\leq N$, has Pfaffian kernel
 $ \mathsf G_NK_N^0(h_N(x),h_N(y))\mathsf G_N^{\mathsf T}$,
which equals $\mathsf K_N(x,y)$ given by \eqref{eq:scaled-FNH-kernel}.\\

\noindent\textbf{Step 2. } Define 
\[ 
 F_{j,N}(u)
 =e^{a_Nu^2}\int_u^\infty e^{-a_Nz^2}\phi_j(z)\,dz.
\]
In this step, we show that 
\begin{align}
    \label{eq:FNH-D-exact}D_N^{\mathrm F}(u,v)
& =
 \frac1{1+q^2}
 \sum_{ 1\leq j<N,\, j\ {\rm odd}} 
 q^{-2j}
 \left(
 F_{j,N}(u)\phi_j(v)-\phi_j(u)F_{j,N}(v)
 \right),\\
 \label{eq:FNH-S-exact} S_N(u,v)
 &=
 K_N^{\mathrm{Her}}(u,v)
 -
 \frac{\sqrt{2N}\,q^2}{1+q^2}
 \phi_N(u)F_{N-1,N}(v),\\
 \label{eq:FNH-I-exact}\widetilde I_N(u,v)
 &=
-2(u-v) \sum_{j=N}^\infty
  q^{2(j+1)} K_{j+1}^{\mathrm{Her}}(u,v) .
\end{align}   Indeed, in the notation from \cite{ForresterNagaoHonner1999}, $ w(u)=e^{-u^2}$, $ C_j(u)=2^{-j}H_j(u)$, and $\gamma_j=j+\frac12$.
We use the monic Hermite norm $\mathfrak h_j=\int_{\mathbb R}w(u)C_j(u)^2\,du$ to fix the normalization in \cite[(4.1)]{ForresterNagaoHonner1999}.  Since $\int_{\mathbb R}e^{-u^2}H_j(u)^2\,du
 =2^jj!\sqrt{\pi}$,
we have
\begin{equation}\label{eq:correct-Hermite-norm}
 \mathfrak h_j=2^{-j}\sqrt{\pi}\,j!,
 \qquad\mbox{ and hence }\qquad
 \phi_j(u)= \sqrt{w(u)}\,C_j(u)/\sqrt{\mathfrak h_j }.
\end{equation} 
We begin by proving the formula \eqref{eq:FNH-D-exact} for $D_N^{\mathrm F}$.
For the $\beta'=4$ initial condition,
\cite[(4.5), (4.6)]{ForresterNagaoHonner1999} gives
\begin{equation}\label{eq:FNH-r-values}
 r_\ell(0)
 =2^{-2\ell}\sqrt{\pi}\,(2\ell+1)!,
 \qquad
 R_{2\ell+1}(u;\tau_N)=C_{2\ell+1}(u).
\end{equation}
Let $j=2\ell+1$ be odd. Using \cite[(4.37)]{ForresterNagaoHonner1999} with $2m=j-1$ and multiply it by $\sqrt{w(u)}$, we get
\begin{equation}\label{eq:FNH-even-R-F}
 \sqrt{w(u)}\,e^{\gamma_{j-1}\tau_N}R_{j-1}(u;\tau_N)
 =
 \frac{2q^{-j+1/2}\sqrt{\mathfrak h_j}}{1+q^2}F_{j,N}(u),
\end{equation}
where we use
$ a_N+ 1/2=1/(1+q^2)$
 and 
 $H_j(z)e^{-z^2/2}=2^j\sqrt{\mathfrak h_j}\,\phi_j(z)$.
On the other hand, \eqref{eq:FNH-r-values} and
\eqref{eq:correct-Hermite-norm} give
\begin{equation}\label{eq:FNH-odd-R-phi}
 \sqrt{w(u)}\,e^{\gamma_j\tau_N}R_j(u;\tau_N)
 =
 q^{-j-1/2}\sqrt{\mathfrak h_j}\,\phi_j(u),
 \qquad
 r_\ell(0)=2\mathfrak h_j.
\end{equation}
Inserting \eqref{eq:FNH-even-R-F} and
\eqref{eq:FNH-odd-R-phi} into
\cite[(2.22)]{ForresterNagaoHonner1999}, term by term for
$j=2\ell+1$, we conclude the proof of \eqref{eq:FNH-D-exact}. 
  
  We next prove the formula \eqref{eq:FNH-S-exact} for $S_N$.
By \cite[(2.27)--(2.29)]{ForresterNagaoHonner1999}, $S_N=S_{N,1}+S_{N,2}$. At equal times,
\[
 S_{N,1}(u,v)
 =
 \sqrt{w(u)w(v)}
 \sum_{j=0}^{N-1} C_j(u)C_j(v)/\mathfrak h_j 
 =
 K_N^{\mathrm{Her}}(u,v).
\]
Combining \cite[(4.35)]{ForresterNagaoHonner1999} with  $\gamma_N=N+ 1/2$, 
$\mathfrak h_{N-1}/\mathfrak h_N=2/N$,   and \eqref{eq:FNH-even-R-F} with $j=N-1$,
\[
 S_{N,2}(u,v)
 =
 -\frac N2
 \frac{\sqrt{w(u)}\,C_N(u)e^{-\gamma_N\tau_N}}{\mathfrak h_N}
 \left(
 \sqrt{w(v)}\,R_{N-2}(v;\tau_N)
 e^{\gamma_{N-2}\tau_N}
 \right) =
 -\frac{\sqrt{2N}\,q^2}{1+q^2}
 \phi_N(u)F_{N-1,N}(v),
\] and we conclude the proof of \eqref{eq:FNH-S-exact}.

We next prove the formula \eqref{eq:FNH-I-exact} for
$\widetilde I_N$. Put
 $\Uppsi_j(u) =e^{\gamma_j\tau_N}\Phi_j(u;\tau_N)$.
By \cite[(2.18), (4.8)]{ForresterNagaoHonner1999} and
\eqref{eq:correct-Hermite-norm},  
\begin{align*}
 \Uppsi_{2\ell}(u)
 &=
 \frac{r_\ell(0)q^{\gamma_{2\ell+1}}}
 {\sqrt{\mathfrak h_{2\ell+1}}}\phi_{2\ell+1}(u), \\
 \Uppsi_{2\ell+1}(u)
 &=
 -\frac{r_\ell(0)q^{\gamma_{2\ell}}}
 {\sqrt{\mathfrak h_{2\ell}}}\phi_{2\ell}(u)
 +
 \frac{(\ell+1)r_\ell(0)q^{\gamma_{2\ell+2}}}
 {\sqrt{\mathfrak h_{2\ell+2}}}\phi_{2\ell+2}(u). 
\end{align*}
The norm ratios are
\begin{equation*} 
 \frac{r_\ell(0)}
 {\sqrt{\mathfrak h_{2\ell}\mathfrak h_{2\ell+1}}}
 =\sqrt{2(2\ell+1)},
 \qquad
 \frac{(\ell+1)r_\ell(0)}
 {\sqrt{\mathfrak h_{2\ell+1}\mathfrak h_{2\ell+2}}}
 =\sqrt{2(2\ell+2)}.
\end{equation*}
The Hermite recurrence gives  
\begin{equation}\label{eq:Hermite-CD-wedge}
 \sqrt{2m}
 \left(
 \phi_{m-1}(u)\phi_m(v)
 -\phi_m(u)\phi_{m-1}(v)
 \right)
 =
 -2(u-v)K_m^{\mathrm{Her}}(u,v)\qquad\mbox{ for }m\in\mathbb N.
\end{equation}
For a finite truncation of  
\cite[(2.21)]{ForresterNagaoHonner1999}, using the above four displays, 
we know that its $\ell$th summand is
\begin{align*}  -2(u-v)q^{2(2\ell+1)}K_{2\ell+1}^{\mathrm{Her}}(u,v)  
 -2(u-v) q^{2(2\ell+2)} K_{2\ell+2}^{\mathrm{Her}}(u,v).
\end{align*}
To pass from finite partial sums to the infinite series, choose $\eta\in(q^2,1)$.  Mehler's formula
\cite[(18.18.28)]{DLMF}
gives
$ \sum_{j=0}^{\infty}\eta^j\phi_j(z)^2
 = 1/\sqrt{\pi(1-\eta^2)} \cdot
e^{- z^2 (1-\eta)/(1+\eta) }$.
Consequently, uniformly when $(u,v)$ ranges over a compact set, the
absolute value of the $j$th summand in \eqref{eq:FNH-I-exact} is bounded by
$ C\sqrt{j+1}\left( q^2/\eta \right)^j$.
This bound is summable because $q^2/\eta<1$.  Hence the resulting
series converges absolutely and locally uniformly, and we may sum the
finite identities above and pass to the limit, obtaining \eqref{eq:FNH-I-exact}.\\

\noindent\textbf{Step 3. } In this step, we conclude the proof of the lemma.  
By Definition~\ref{def:kernel KN},
\begin{equation*} 
 \phi_j(h_N(x))
 =2^{1/4}N^{-1/12}\psi_{j,N}(x),
 \qquad
 F_{j,N}(h_N(x))
 =2^{-1/4}N^{-1/4}J_{j,N}(x).
\end{equation*}
By a change of variables in the identities \eqref{eq:FNH-D-exact}, \eqref{eq:FNH-S-exact}, and \eqref{eq:FNH-I-exact}, and using the above display, we conclude that $\mathsf K_N=K^{t,N}$, where  $\mathsf K_N$ is given by \eqref{eq:scaled-FNH-kernel}. Note also that each entry of $K^{t,N}$ is continuous and locally bounded, in view of the local uniform convergence together with polynomial--Gaussian majorants uniform on compact $(x,y)$-sets for the remaining defining integrals. We conclude the proof.
\end{proof}

\begin{proposition}\label{prop:PPConvergence}
Assume the same notation as in Lemma \ref{lem:PPP structure}. Then for $t\in(0,\infty)$, as $N\to\infty$ through even $N$, $M^{t,N}\Rightarrow M^{t;\mathrm{hs};\infty}$ in the space of locally bounded measures on $\mathbb R$ with the vague topology, where  
\[
 M^{t;\mathrm{hs};\infty}(\omega,A)
  =\sum_{k=1}^{\infty}\mathbf 1\{\hsai_k(t,\omega)\in A\}.
\]  
\end{proposition}
\begin{proof} In this proof, we use the notation from Definition \ref{def:kernel KN}, and denote $\tau_N=tN^{-1/3}$, $r_N=q_N^2$,
\[ 
 \mu_m=\sqrt{2m+1},
 \qquad
 \delta_m=2^{-1/2}m^{-1/6},
 \qquad
 \overline\phi_m=2^{-1/4}m^{1/12}\phi_m,\qquad m\in\mathbb N.
\]
We want to use the convergence criterion of Pfaffian point processes given by \cite[Proposition~5.10]{dimitrov2026half}, and the major part is to prove  $K^{t,N}(x,y)\to K^{t;\mathrm{hs};\infty}(x,y)=\hski(t,x;t,y)$ locally uniformly in $(x,y)\in\mathbb R^2$.

For clarity, we split the proof into five steps.\\

\noindent\textbf{Step 1. } In this step, we prove the following limits and estimates: there exist an absolute constant  $C \in(0,\infty)$ and constants $C_L\in(0,\infty)$, $m_L\in\mathbb N$ depending only on $L>0$, such that
\begin{align}
\label{eq:PP-Hermite-Airy} \lim_{m\to\infty}\overline\phi_m(\mu_m+\delta_ms)&=\Ai(s) \qquad&\mbox{ uniformly on compact subsets of }s\in\mathbb R,\\
    \label{eq:PP-Hermite-bound}  
     \left|\overline\phi_m(\mu_m+\delta_ms)\right|
 &\leq C_Le^{-s}&\mbox{ for } m\geq m_L,\, s\geq-L,\,L>0,\\
 \label{eq:PP-global-Hermite}  \sup_{u\in\mathbb R}|\phi_m(u)|
 &\leq C(m+1)^{-1/12}\qquad&\mbox{ for each }m\in\mathbb Z_{\geq0},\end{align}
 and for each fixed $A,U\in(0,\infty)$, 
 \begin{align} 
 \label{eq:PP-moving-Hermite}   & \lim_{N\to\infty}\sup_{\substack{|x|\leq A,\ |m-N|\leq UN^{1/3}}}
 \left|
 \psi_{m,N}(x)-
 \Ai\left(x-\frac{m-N}{N^{1/3}}\right)
 \right|&=0 ,\\
 \label{eq:PP-moving-kernel} &\lim_{N\to\infty}\sup_{\substack{|x|,|y|\leq A,\ |m-N|\leq UN^{1/3}}}
 \left|
 d_NK_m^{\mathrm{Her}}(h_N(x),h_N(y))
 -K_{\Ai}\left(x-\frac{m-N}{N^{1/3}},
                 y-\frac{m-N}{N^{1/3}}\right)
 \right|&=0 ,\\ 
 \label{eq:PP-J-limit}&\lim_{N\to\infty} \sup_{\substack{|x|\leq A,\ 0\leq(N-m)N^{-1/3}\leq U}}
 \left|J_{m,N}(x)-
 h_t\left(x+(N-m)N^{-1/3}\right)\right|&=0,
\end{align}
where we recall that $h_t(z)=\int_0^\infty e^{-ts}\Ai(z+s)\,ds$. 

Indeed, \eqref{eq:PP-Hermite-Airy}, \eqref{eq:PP-Hermite-bound} follow from \cite[Proposition~1]{JohnstoneMa2012}, and \eqref{eq:PP-global-Hermite} follows from \cite[(61)]{Johansson2007Gumbel}. 
Observe 
\begin{equation}
     \label{eq:PP-psi-normalization} \psi_{m,N}(v)
 =\left( N/m\right)^{1/12}
 \overline\phi_m(\mu_m+\delta_m\sigma_{m,N}(v)),\qquad\mbox{ where }\sigma_{m,N}(v)= (h_N(v)-\mu_m)/\delta_m.
\end{equation}
Since $\sigma_{m,N}(x)=x-(m-N)N^{-1/3}+O_{A,U}(N^{-1/3})$ uniformly on $\{|x|\leq A,\ |m-N|\leq UN^{1/3}\}$, we conclude \eqref{eq:PP-moving-Hermite} from   \eqref{eq:PP-Hermite-Airy}, \eqref{eq:PP-psi-normalization}, and $m/N\to1$. 
Moreover, combining \eqref{eq:PP-Hermite-bound} and \eqref{eq:PP-psi-normalization}, we know that
there exist $C_{A,U},c_{A,U}\in(0,\infty)$ and $N_{A,U}\in\mathbb N$ depending only on $A,U\in(0,\infty)$, such that 
\begin{equation}
 |\psi_{m,N}(x+w)|+ |\psi_{m-1,N}(x+w)|\leq C_{A,U}e^{-c_{A,U}w}
 \qquad\mbox{ for } |x|\leq A,\,w\geq0, |m-N|\leq UN^{1/3},\, N\geq N_{A,U}. 
 \label{eq:PP-moving-bound}
\end{equation} 

We next prove \eqref{eq:PP-moving-kernel}.
 The Hermite projection identity
\cite[(19)]{JohnstoneMa2012} with
$z=d_Nw$, gives 
 \begin{align}
 &d_NK_m^{\mathrm{Her}}(h_N(x),h_N(y)) =\frac12\left(\frac{m}{N}\right)^{1/2}
 \int_0^\infty
 \left(
 \psi_{m,N}(x+w)\psi_{m-1,N}(y+w)
 +\psi_{m-1,N}(x+w)\psi_{m,N}(y+w)
 \right)\,dw.
 \label{eq:PP-Hermite-projection}
\end{align}
To take the $N\to\infty$ limit of this integral  on $\{|x|,|y|\leq A,\ |m-N|\leq UN^{1/3}\}$, we first restrict to $w\in[0,R]$. On this interval, \eqref{eq:PP-moving-Hermite} gives the uniform convergence of the integrand to the limit. The integral on $w\in(R,\infty)$ is bounded by $C_{A,U}e^{-c_{A,U}R}$ in view of  \eqref{eq:PP-moving-bound}, and the limiting Airy integral has the same type of tail bound. Letting first $N\to\infty$ and then $R\to\infty$ therefore proves \eqref{eq:PP-moving-kernel}.

We finally prove \eqref{eq:PP-J-limit}. Indeed, for each $R>0$ we have 
\begin{equation} \label{eq:PP-exponent-bound}
    \begin{split}
        \lim_{N\to\infty}a_N\left(h_N(x+s)^2-h_N(x)^2\right)=ts\qquad\mbox{ uniformly for }|x|\leq A,\, s\in[0,R].
        \end{split}
\end{equation} This follows from  $2a_N=\tanh(tN^{-1/3})$ and $h_N(x+s)^2-h_N(x)^2=d_Ns (2\sqrt{2N}+d_N(2x+s) )$. 
Equations \eqref{eq:PP-moving-Hermite} and \eqref{eq:PP-moving-bound} give convergence of the integrand in  $J_{m,N}$,
uniformly for  $s\in[0,R]$. For $N$ sufficiently large depending on $A$, we have $h_N(x)>0$ for
$|x|\leq A$. Since $h_N$ is increasing, for $s\geq0$, $\exp (-a_N(h_N(x+s)^2-h_N(x)^2)) \in(0,1]$.
Consequently, \eqref{eq:PP-moving-bound} bounds the part with $s\geq R$
by $C_{A,U}e^{-c_{A,U}R}$, uniformly in $m,N,x$ in $\{|x|\leq A,\ 0\leq(N-m)N^{-1/3}\leq U\}$;
Lemma~\ref{lem:Airy decay} gives the same type of tail bound for the
limiting integral. Letting first $N\to\infty$ and then $R\to\infty$ proves \eqref{eq:PP-J-limit}.\\

\noindent\textbf{Step 2. } In this step, we show that, locally uniformly in $(x,y)$,
\begin{align}
\lim_{N\to\infty} K_{11}^{t,N}(x,y)
 &=\frac14\int_0^\infty e^{-2tu}
 \left(
 \Ai(x+u)h_t(y+u)-h_t(x+u)\Ai(y+u)
 \right)\,du=K_{11}^{t;\mathrm{hs};\infty}(x,y),
 \label{eq:PP-K11-limit}
\end{align}
where the second equality follows from $h_t'=th_t-\Ai$ and \eqref{eq:alt11} in Lemma \ref{lem:alternative kernels}.
For $0\leq k\leq N/2-1$, denote $u_{k,N}=(2k+1)N^{-1/3}$ and $m_{k,N}=N-1-2k$, and hence $u_{k,N}=(N-m_{k,N})N^{-1/3}$. Rewriting
\eqref{eq:K11Ndef},  
\begin{align*}
 K_{11}^{t,N}(x,y)
 =\frac{N^{-1/3}}{1+r_N}\sum_{k=0}^{N/2-1}e^{-2tu_{k,N}}
 \left( \psi_{m_{k,N},N}(x)J_{m_{k,N},N}(y) 
  -J_{m_{k,N},N}(x)\psi_{m_{k,N},N}(y)\right).
\end{align*}
Fix $A,U>0$. By \eqref{eq:PP-moving-Hermite},
\eqref{eq:PP-J-limit}, uniformly for $|x|,|y|\leq A$,
$u_{k,N}\leq U$, the expression in parentheses converges to
\[
 \Ai(x+u_{k,N})h_t(y+u_{k,N})
 -h_t(x+u_{k,N})\Ai(y+u_{k,N}).
\]
The points $u_{k,N}$ form a mesh of size $2N^{-1/3}$,
and $r_N\to1$. Hence, locally uniformly
in $(x,y)$,
\begin{equation} \label{eq:PP-K11-truncated}
    \begin{split}
         &\lim_{N\to\infty}\frac{N^{-1/3}}{1+r_N}
 \sum_{\substack{0\leq k<N/2,\, u_{k,N}\leq U}}e^{-2tu_{k,N}}
 \left(\psi_{m_{k,N},N}(x)J_{m_{k,N},N}(y)
 -J_{m_{k,N},N}(x)\psi_{m_{k,N},N}(y)\right) \\
 & =\frac14\int_0^Ue^{-2tu}
 \left(\Ai(x+u)h_t(y+u)-h_t(x+u)\Ai(y+u)\right)\,du.
    \end{split}
\end{equation} 
It remains to control the remaining indices uniformly. By
\eqref{eq:PP-global-Hermite}, for $m\geq1$ and $z\in\mathbb R$,
\begin{equation}
 |\psi_{m,N}(z)|
 \leq C\left( N/(m+1)\right)^{1/12}.
 \label{eq:PP-K11-global-psi}
\end{equation}
For $|z|\leq A$ and $s\geq0$, we have
 $a_N (h_N(z+s)^2-h_N(z)^2 )
 =a_Nd_Ns (2\sqrt{2N}+d_N(2z+s) )$.
Note also that $a_Nd_N(2\sqrt{2N}-2Ad_N)\to t$, using \eqref{eq:PP-K11-global-psi}, for all sufficiently large $N$, we have
\begin{equation}
 \sup_{|z|\leq A}|J_{m,N}(z)|
 \leq C_{A,t}\left( N/(m+1)\right)^{1/12},
 \qquad 1\leq m<N.
 \label{eq:PP-K11-global-J}
\end{equation}
Indeed, we can bound the defining integral for $J_{m,N}(z)$ by  $C\left( N/(m+1)\right)^{1/12}
 \int_0^\infty e^{-ts/2}\,ds$.
Combining \eqref{eq:PP-K11-global-psi} and
\eqref{eq:PP-K11-global-J}, and splitting the remaining indices
according as $m_{k,N}\geq N/2$ or $m_{k,N}<N/2$, yields
\begin{align*}
 &\sup_{|x|,|y|\leq A}\frac{N^{-1/3}}{1+r_N}
 \sum_{\substack{0\leq k<N/2,\,u_{k,N}>U}}e^{-2tu_{k,N}}
 \left|\psi_{m_{k,N},N}(x)J_{m_{k,N},N}(y)
 -J_{m_{k,N},N}(x)\psi_{m_{k,N},N}(y)\right|\\
 & \leq C_{A,t}N^{-1/3}
 \sum_{\substack{0\leq k<N/2,\,u_{k,N}>U,\ m_{k,N}\geq N/2}}
 e^{-2tu_{k,N}}
 +C_{A,t}N^{5/6}e^{-tN^{2/3}} \leq C_{A,t}e^{-2tU}+o(1).
\end{align*}
Indeed, the first sum is bounded by the geometric series with ratio
$e^{-4tN^{-1/3}}$, while $m_{k,N}<N/2$ implies
$u_{k,N}>N^{2/3}/2$. Finally, Lemma~\ref{lem:Airy decay} shows that, uniformly for $|x|,|y|\leq A$, the
absolute value of the integral over $u\in(U,\infty)$ in
\eqref{eq:PP-K11-limit} is at most $C_{A,t}e^{-2tU}$.
Letting first $N\to\infty$ in \eqref{eq:PP-K11-truncated} and then
$U\to\infty$ proves \eqref{eq:PP-K11-limit}.\\

\noindent\textbf{Step 3. } In this step, we observe that, locally uniformly in $(x,y)$, 
\begin{equation} \label{eq:PP-K12-limit}
\lim_{N\to\infty} K_{12}^{t,N}(x,y)=
 K_{\Ai}(x,y)-\frac12h_t(x)\Ai(y)
 =K_{12}^{t;\mathrm{hs};\infty}(x,y).
\end{equation}
Indeed, recall from \eqref{eq:K12Ndef} that 
\begin{equation*}
 K_{12}^{t,N}(x,y)
 =d_NK_N^{\mathrm{Her}}(h_N(x),h_N(y))
 -\frac{q_N^2}{1+q_N^2}\psi_{N,N}(y)J_{N-1,N}(x). 
\end{equation*}
Then \eqref{eq:PP-K12-limit} follows from combining \eqref{eq:PP-moving-Hermite}, \eqref{eq:PP-moving-kernel}, \eqref{eq:PP-J-limit},  
$q_N^2/(1+q_N^2)\to1/2$, and \eqref{eq:alt12} in Lemma \ref{lem:alternative kernels}.\\

\noindent\textbf{Step 4. } In this step, we show that, locally uniformly in $(x,y)$,
\begin{align}
\lim_{N\to\infty} K_{22}^{t,N}(x,y)=
-(x-y)\int_0^\infty
 e^{-2tu}K_{\Ai}(x-u,y-u)\,du=K_{22}^{t;\mathrm{hs};\infty}(x,y).
 \label{eq:PP-K22-limit}
\end{align}
Indeed, we first verify the second equality above. Lemma~\ref{lem:Airy decay} gives
$K_{\Ai}(x-u,y-u)=O_A(1+u^{1/2})$ uniformly for $|x|,|y|\leq A$, while $\partial_uK_{\Ai}(x-u,y-u)=\Ai(x-u)\Ai(y-u)$.
Thus integration by parts gives
\[
 \int_0^\infty e^{-2tu}K_{\Ai}(x-u,y-u)\,du
 =
 \frac1{2t}\left(
 K_{\Ai}(x,y)
 +\int_0^\infty e^{-2tu}\Ai(x-u)\Ai(y-u)\,du
 \right).
\]
By \eqref{eq:Darboux KAi} in Definition~\ref{def:Airy kernel},
\[
 \Ai(x-u)\partial_u\Ai(y-u)
 -\partial_u\Ai(x-u)\Ai(y-u)
 =-(x-y)K_{\Ai}(x-u,y-u).
\]
In view of \eqref{eq:alt22} in Lemma \ref{lem:alternative kernels}, we conclude the second equality in \eqref{eq:PP-K22-limit}.

We next prove the limit in \eqref{eq:PP-K22-limit}.
Since $2d_N^2=N^{-1/3}$, setting $\ell=m-N$ in
\eqref{eq:K22Ndef} gives
\begin{equation*}
 K_{22}^{t,N}(x,y)
 =
 -(x-y)N^{-1/3}\sum_{\ell=1}^{\infty}
 e^{-2t\ell N^{-1/3}}
 d_NK_{N+\ell}^{\mathrm{Her}}(h_N(x),h_N(y)).
\end{equation*}
Fix $A>0$. For every fixed $U>0$, \eqref{eq:PP-moving-kernel}
implies, uniformly for $|x|,|y|\leq A$ and
$1\leq\ell\leq UN^{1/3}$,
\[
 d_NK_{N+\ell}^{\mathrm{Her}}(h_N(x),h_N(y))
 =
 K_{\Ai}\left(x-\ell N^{-1/3},y-\ell N^{-1/3}\right)+o(1).
\]
Consequently, the truncated sums converge locally
uniformly in $(x,y)\in\mathbb R^2$:
\begin{align*}
 \lim_{N\to\infty} N^{-1/3}\sum_{\ell=1}^{\lfloor UN^{1/3}\rfloor}
 e^{-2t\ell N^{-1/3}}
 d_NK_{N+\ell}^{\mathrm{Her}}(h_N(x),h_N(y))=
 \int_0^U e^{-2tu}K_{\Ai}(x-u,y-u)\,du .
\end{align*}
It remains to justify the truncation. 
By
\eqref{eq:PP-moving-kernel} with $m=N$,
$d_NK_N^{\mathrm{Her}}(h_N(x),h_N(y))$ is bounded uniformly for
$|x|,|y|\leq A$ and all sufficiently large $N$. Using
\eqref{eq:PP-global-Hermite} and
  $K_{N+\ell}^{\mathrm{Her}}(u,v)
 =
 K_N^{\mathrm{Her}}(u,v)
 +\sum_{j=N}^{N+\ell-1}\phi_j(u)\phi_j(v)$, uniformly for $|x|,|y|\leq A$ and $\ell\geq1$,
$ d_N\left|
 K_{N+\ell}^{\mathrm{Her}}(h_N(x),h_N(y))
 \right|
 \leq C_A+C\ell N^{-1/3}$.
Hence, for $U\geq1$,
\[
 \sup_{|x|,|y|\leq A}
 N^{-1/3}\sum_{\ell>UN^{1/3}}
 e^{-2t\ell N^{-1/3}}
 d_N\left|
 K_{N+\ell}^{\mathrm{Her}}(h_N(x),h_N(y))
 \right|
 \leq C_{A,t}e^{-tU}.
\]
Indeed, the left side is bounded by the corresponding geometric
sums with summands
$e^{-2t\ell N^{-1/3}}(1+\ell N^{-1/3})$.
Letting first $N\to\infty$ and then $U\to\infty$ proves the first equality in \eqref{eq:PP-K22-limit}.\\

\noindent\textbf{Step 5. } In this step, we conclude the proof of the proposition.  Combining \eqref{eq:PP-K11-limit}, \eqref{eq:PP-K12-limit},   
\eqref{eq:PP-K22-limit},  and the fact that
$ K_{21}^{t,N}(x,y)=-K_{12}^{t,N}(y,x)$ and 
$K_{21}^{t;\mathrm{hs};\infty}(x,y)=-K_{12}^{t;\mathrm{hs};\infty}(y,x)$,  we have
\[
 \lim_{N\to\infty,\,N\in2\mathbb N}
 \max_{i,j\in\{1,2\}}
 \sup_{x,y\in \mathcal V}
 \left|K_{ij}^{t,N}(x,y)-K^{t;\mathrm{hs};\infty}_{ij}(x,y)\right|=0\qquad\mbox{ for every compact }\mathcal V\subset \mathbb R.
\]
By the integral formulas on the right sides of \eqref{eq:PP-K11-limit}, \eqref{eq:PP-K12-limit}, and \eqref{eq:PP-K22-limit}, together with
Lemma~\ref{lem:Airy decay} and dominated convergence, we know that $K^{t;\mathrm{hs};\infty}:\mathbb R^2\to\operatorname{Mat}_2(\mathbb C)$ is continuous.
Using \cite[Proposition~5.10]{dimitrov2026half},
there exists a Pfaffian point process $M$ on $\mathbb R$ with
reference measure $\mathrm{Leb}$ and correlation kernel $K^{t;\mathrm{hs};\infty}$ such
that $M^{t,2n}\Rightarrow M$.
On the other hand, Proposition~\ref{Prop.HSAiryinfty} with
$\mathsf S=\{t\}$, and \cite[Lemma~5.13]{dimitrov2026half} show that the random measure
$M^{t;\mathrm{hs};\infty}$ is a Pfaffian point
process on $\mathbb R$ with the same reference measure
$\mathrm{Leb}$ and the same correlation kernel $K^{t;\mathrm{hs};\infty}$. By
\cite[Proposition~5.8(3)]{dimitrov2026half},
$M\overset{d}=M^{t;\mathrm{hs};\infty}$, which proves the proposition.
\end{proof}

\begin{lemma}\label{lem:PPP upper tail}
Fix $t\in(0,\infty)$. For each even $N\in\mathbb N$, let $M^{t,N}$ 
be the point process from Lemma~\ref{lem:PPP structure}. There exist
constants $C,c\in(0,\infty)$ and $N_0\in\mathbb N$ such that, for every
$R\geq0$ and every even $N\geq N_0$,
\[
 \mathbb E\left[M^{t,N}((R,\infty))\right]
 \leq Ce^{-cR}.
\] 
\end{lemma}
\begin{proof}
    By \cite[(5.3) and (5.12)]{dimitrov2026half} and monotone convergence,
    \begin{equation}\label{eq:mdyrtu}
    \mathbb E\left[M^{t,N}((R,\infty))\right] =
 \int_R^{\infty} K_{12}^{t,N}(x,x)\,dx.
     \end{equation}
Recall from \eqref{eq:K12Ndef} that
\begin{equation}\label{eq:rnmtu}
    K_{12}^{t,N}(x,y)= 
d_NK_N^{\mathrm{Her}}(h_N(x),h_N(y))  -\frac{q_N^2}{1+q_N^2} 
 \psi_{N,N}(y)J_{N-1,N}(x).\end{equation}
In view of \eqref{eq:PP-moving-bound} with $A=U=1$, $m=N$, $x=0$, there exist constants $C_0,c_0\in(0,\infty)$ and $N_0\in\mathbb N$ such that,
\begin{equation}\label{eq:PPP-upper-tail-psi}
 |\psi_{N,N}(u)|+|\psi_{N-1,N}(u)|
 \leq C_0e^{-c_0u}\qquad\mbox{ for each }u\geq0,\,N\geq N_0.
\end{equation}
Using \eqref{eq:PP-Hermite-projection} with $m=N$ and $y=x$, and combining with \eqref{eq:PPP-upper-tail-psi}, for each $x\geq0$,
\[
 \left|d_NK_N^{\mathrm{Her}}(h_N(x),h_N(x))\right|
 =
\left| \int_0^\infty
 \psi_{N,N}(x+w)\psi_{N-1,N}(x+w)\,dw\right|\leq
 C_0^2\int_0^\infty e^{-2c_0(x+w)}\,dw
 =
 \frac{C_0^2}{2c_0}e^{-2c_0x}.
\]
 For $v\ge x\ge0$, we have $h_N(v)\ge h_N(x)>0$ and $a_N>0$.
Thus the exponential weight in $J_{N-1,N}(x)$ is at most one.
By \eqref{eq:PPP-upper-tail-psi},
 \[
 |J_{N-1,N}(x)|
 \leq
 \int_x^\infty|\psi_{N-1,N}(v)|\,dv
 \leq C_0/c_0\cdot e^{-c_0x}.
\]
Because $q_N\in(0,1)$, we have
$q_N^2/(1+q_N^2)<1/2$. Combining with the last display and \eqref{eq:rnmtu}, 
\[
 \left|K_{12}^{t,N}(x,x)\right|
 \leq C_0^2/c_0\cdot e^{-2c_0x}\qquad\mbox{ for each }x\geq0,\,N\geq N_0.
\]
Plugging in \eqref{eq:mdyrtu}, we have $ \mathbb E [M^{t,N}((R,\infty)) ]
 \leq C_0^2/(2c_0^2)\cdot e^{-2c_0R}$. We conclude the proof.
\end{proof}

\begin{proof}[Proof of Proposition \ref{prop:scaling limit PP}] 
Fix $t>0$. For each even $N$, let $ b_N(t)=\min(Z_N^N(t),-N)$ and define $\mathsf Z^N(t)=(\mathsf Z_i^N(t))_{i\geq1}\in\mathbb R^\infty$
by $\mathsf Z_i^N(t)=Z_i^N(t)$ if $1\leq i\leq N$ and $\mathsf Z_i^N(t)=b_N(t)-(i-N)$ if $i>N$. Define  
\[
 \mathsf M^{t,N}
 =\sum_{i=1}^\infty\delta_{\mathsf Z_i^N(t)}
 =M^{t,N}+\sum_{j=1}^\infty\delta_{b_N(t)-j}.
\]
This random measure on $\mathbb R$ is locally finite. The second summand is supported on $(-\infty,-N-1]$, and hence it converges
vaguely to the zero measure. By Proposition~\ref{prop:PPConvergence}
and Slutsky's theorem,
\begin{equation}\label{eq:augmented-PP-convergence}
 \mathsf M^{t,N}\Longrightarrow M^{t;\mathrm{hs};\infty}
 \qquad\text{as $N\to\infty$ through even integers}.
\end{equation}
Since $\mathsf Z_1^N(t)=Z_1^N(t)$, we have
$ \mathbb P\left(\mathsf Z_1^N(t)\geq a\right)
 \leq
 \mathbb P\left(M^{t,N}((a-1,\infty))\geq1\right)
 \leq
 \mathbb E M^{t,N}((a-1,\infty))$ for every $a>1$.
Thus Lemma~\ref{lem:PPP upper tail} yields
\begin{equation}\label{eq:rightmost-upper-tightness}
 \lim_{a\to\infty}\limsup_{N\to\infty,\,N\in2\mathbb N } 
 \mathbb P\left(\mathsf Z_1^N(t)\geq a\right)=0.
\end{equation}
Since $\hsai$ is an $\mathbb N$-indexed  line ensemble,
$ M^{t;\mathrm{hs};\infty}(\mathbb R) 
 =\infty$
almost surely.
Combining this with \eqref{eq:augmented-PP-convergence} and
\eqref{eq:rightmost-upper-tightness}, after reindexing the sequence
by the even integers, all the hypotheses of
\cite[Proposition~2.21]{dimitrov2024airy} are verified. Hence
$\{\mathsf Z^N(t):N\ {\rm even}\}$ is tight in $\mathbb R^\infty$,
and in particular each coordinate family is tight. Applying
\cite[Proposition~2.19]{dimitrov2024airy} with $r=1$ to
\eqref{eq:augmented-PP-convergence}, there exists an ordered random
vector $X=(X_i)_{i\geq1}$ such that
$\mathsf Z^N(t)\Rightarrow X$ in $\mathbb R^\infty$, and
$ \sum_{i=1}^\infty\delta_{X_i}
 \overset{d}{=}M^{t;\mathrm{hs};\infty}$.
  In view of \cite[Corollary~2.20]{dimitrov2024airy}, $X\overset{d}{=} (\hsai_i(t))_{i\geq1}$. Therefore $ \mathsf Z^N(t)\Rightarrow (\hsai_i(t) )_{i\geq1}$.  
For fixed $k$,   $\mathsf Z_k^N(t)=Z_k^N(t)$ once $N\geq k$. We conclude the proof.
\end{proof}

\bibliographystyle{goodbibtexstyle.bst}
\bibliography{Concentration}

\end{document}